\documentclass[a4paper,12pt]{article}

\usepackage{amsthm,mathrsfs,amssymb,amsmath,bbm}
\usepackage{fullpage,enumerate}
\usepackage{makeidx}
\usepackage{graphicx}
\usepackage{xr,color}
\usepackage{tcolorbox}
\usepackage[pdftex]{hyperref}

\newenvironment{thnumber}{

\vspace{-.05in} \begin{enumerate}}{\end{enumerate}} % enumeration within Theorems, Lemmas etc.

\newtheorem{theorem}{Theorem}
\newtheorem{lemma}[theorem]{Lemma}
\newtheorem{proposition}[theorem]{Proposition}
\newtheorem{corollary}[theorem]{Corollary}

\theoremstyle{remark}
\newtheorem{remark}[theorem]{Remark}
\newtheorem{example}[theorem]{Example}

\numberwithin{theorem}{section}
\numberwithin{equation}{section}

\newcommand{\R}{\mathbb{R}}
\newcommand{\C}{\mathbb{C}}

\newcommand{\E}{{\mathbf E}}
\renewcommand{\P}{{\mathbf P}}
\newcommand{\Q}{{\mathbf Q}}

\newcommand{\defeq}{\stackrel{\mathrm{def}}{=}}

\newcommand{\sE}{\mathcal{E}}
\newcommand{\sF}{\mathcal{F}}

\newcommand{\sG}{\mathcal{G}}

\newcommand{\sT}{\mathcal{T}}
\newcommand{\sS}{\mathcal{S}}
\newcommand{\sJ}{\mathcal{J}}
\newcommand{\sP}{\mathcal{P}}

\newcommand{\sL}{\mathcal{L}}

\newcommand{\scrA}{\mathscr{A}}

\newcommand{\sM}{\mathscr{M}}

\newcommand{\real}{{\rm Re}}

\newcommand{\sB}{\mathcal{B}}
\newcommand{\one}{\mathbbm{1}}

\DeclareMathOperator{\sign}{sign}

\newcommand{\mc}{\ast}

\newcommand{\mcs}{\circledast}

\begin{document}

\thispagestyle{empty}

\begin{center}
{\Large \bf Multiplicative convolution with symmetries\\[.3em] in Euclidean space and on the sphere}\\

\vspace{.5cm}

{\large \scshape Felix~Nagel}\footnote{Email: research@felixnagel.org}\\

\vspace{.2cm}

{\it Institute of Mathematical Statistics and Actuarial Science\\
University of Berne, Switzerland}

\vspace{.8cm}

\begin{minipage}{.85\textwidth}
\begin{center}
{\bf Abstract}
\end{center}
\vspace{-.3cm}
Multiplicative convolution $\mu \mc \nu$ of two finite signed measures $\mu$ and $\nu$ on $\R^n$ and a related product $\mu \mcs \nu$ on the sphere $S^{n-1}$ are studied. For fixed $\mu$ the injectivity in $\nu$ of both operations is characterised given an arbitrary group of reflections along the coordinate axes. The results for the sphere yield generalised versions of the theorems in \cite{mol:nag:21} about convex bodies.
\vspace{.2cm}


\vspace{.2cm}

\end{minipage}
\end{center}

\vspace{.5cm}

\section{Introduction}

The space $\R^n$ where $n \geq 1$ can be equipped with componentwise multiplication, and for two finite signed measures $\mu$ and $\nu$ on $\R^n$ the corresponding convolution $\mu \mc \nu$ can be defined. The question arises for which $\mu$ the homogeneous equation $\mu \mc \nu = 0$ implies $\nu = 0$. Such a measure $\mu$ is called universal in this paper.

Similar convolution equations of measures on locally compact Abelian groups or semigroups have been studied for a long time~\cite{den60,dav:sha:87}. Results are applied in numerous papers to characterise various probability laws, see \cite[Section~4]{dav:sha:87} for an overview. More recently they have been considered in relation to inverse problems for regular variation~\cite{jac:mic:ros:sam09,dam:mik:ros:14}.

The main problem in the case of $\R^n$ originates from the fact that it is not a group. Clearly vectors in $(0, \infty)^n$ with componentwise multiplication form a group. Also real vectors with non-zero components, $(\R \setminus \{0\})^n$, with componentwise multiplication are a group. The same is true for any subset $(\R \setminus \{0\})^E$ in the coordinate subspace $\R^E$ where $E \subset \{1, \ldots, n\}$. We may think of $\R^n$ as partitioned into such sets, and the induced decomposition of measures is used in the present study.

Signed measures are considered here (as opposed to non-negative ones) because these are needed in geometric applications, for example in the analytic representation of generalised zonoids, an important class of convex bodies, see for example \cite[p.~195]{schn2}. Moreover signed measures can be identified with continuous linear functionals on standard function spaces.

Since the product of two negative coordinates is positive, all orthants $(\R \setminus \{0\})^E$ for given $E$ should be analysed together. This fact also motivates to consider reflections along coordinate axes, i.e.\ the change of sign of some non-zero coordinates. Symmetries with respect to reflections often appear naturally in applications; e.g.\ in the representation of generalised zonoids the signed measure can always assumed to be symmetric about the origin. It is intuitive that, if $\mu$ is preserved under a reflection, then also $\mu \mc \nu$ has the same property. Thus $\mu$ is not fully universal but it may be universal when $\nu$ is confined to the same class of symmetries. Since we consider signed measures, it makes sense to include both even and odd terms under an arbitrary group of reflections, for example the space of measures that change their sign under reflection about the origin.

Measures and signed measures on the sphere are used in various contexts, for instance in convex geometry as mentioned above and in the representation of the characteristic function of multivariate stable distributions, see~\cite[Theorem~2.3.1]{sam:taq94}. This motivates the study of the product $\mu \mcs \nu = P_S(\mu \mc \nu)$ of two finite signed measures $\mu$ and $\nu$ on $S^{n-1}$. It consists of multiplicative convolution and subsequent projection $P_S$ on the sphere. A similar spherical projection of a measure is used in~\cite{weil82}. As for $\R^n$ we pose the question: For which $\mu$ does $\mu \mcs \nu = 0$ imply $\nu = 0$, given an arbitrary reflection group?

For each of the two products ($\mc$ on $\R^n$, and $\mcs$ on $S^{n-1}$) and an arbitrary choice of reflection symmetries a characterisation of universal signed measures is derived in Sections \ref{sec - universality Rn} and \ref{sec - universality sphere}. Several special cases are stated, for example origin-symmetric and unconditional measures.

Finally multiplicative convolution on the sphere is applied to the geometric problem considered in~\cite[Section~3]{mol:nag:21}. Their goal was to determine which generalised zonoids are "universal" in a sense defined in that paper. In the present study their main results are generalised, weaker conditions are stated, and they are proven to be necessary.

Section~\ref{sec - multiplicative convolutions} starts with the basic definitions and some properties of multiplicative convolution on $\R^n$. A decomposition into zero and non-zero coordinates as mentioned above is defined in Section~\ref{sec - decomp measures}. In Section~\ref{sec - reflections} the role of reflections along coordinate axes is extensively discussed. Here another decomposition into symmetry components is introduced, and an appropriate transformation of such components. In Section~\ref{sec - lifting} a map that connects the problem for $\R^n$ with that for $S^n$ in $\R^{n+1}$ is specified. Section~\ref{sec - uniqueness} contains the basic uniqueness results for measures using harmonic analysis, both for Euclidean space and for the sphere. Sections \ref{sec - universality Rn} and \ref{sec - universality sphere} contain the proofs of the main theorems, viz.\ the characterisations of universality for both situations. Various useful special cases are stated as well. It is then shown in Section~\ref{sec - application} how the results about generalised zonoids in \cite{mol:nag:21} can be generalised using the new characterisations. Section~\ref{sec - measures of degree n} gives a criterion (expressed in terms of orthogonal projections) in order to satisfy some conditions for universality. It is finally demonstrated in Section~\ref{sec - universality lifting} how the results for $\R^n$ can be derived from those for $S^n$. In spite of this fact all proofs leading to the main theorems are separated for both cases into Sections \ref{sec - universality Rn} and \ref{sec - universality sphere} in order to make the two lines of thought more independent.

\section{Multiplicative convolution}
\label{sec - multiplicative convolutions}

\subsection{Preliminaries}

For any set $A$ the family of subsets of $A$ is denoted by $\sP(A)$. We also use the notation $\sP^2(A) = \sP(\sP(A))$. For $n \geq 1$ we set $[n] = \{1, 2, \ldots, n \}$, $\sP_n = \sP([n])$, and
\[
\R^n_+ \,=\, \big\{ x \in \R^n \,:\, \forall j \in [n] \;\; x_j \geq 0 \big\}
\]
Whenever we use any of the sets $[n]$, $\sP_n $, $\R^n$, $\R^n_+$, or $(0,\infty)^n$ it is tacitly understood that $n \geq 1$. The $\sigma$-field of Borel measurable sets in $\R^n$ is denoted by~$\mathcal{B}(\R^n)$. For $A \in \mathcal{B}(\R^n)$ the sub-$\sigma$-field of $\mathcal{B}(\R^n)$ on~$A$ is denoted by~$\mathcal{B}(A)$. A $\sigma$-additive set function $\mu : \mathcal{B}(A) \longrightarrow \R$ with $\mu(\O) = 0$ is called {\em measure on $A$}\footnote{In the literature (see for example \cite[Chapter~4]{cohn} or \cite[Definition~VII.1.1]{elstrodt}) our "measure" is mostly called "finite signed measure". We use this convention to avoid endless repetitions because in the present context most $\sigma$-additive set functions have range~$\R$.}. If the range of $\mu$ is $\R_+$, then we say that $\mu$ is {\em non-negative} and write $\mu \geq 0$. For $A \in \mathcal{B}(\R^n)$ the vector space of measures on $A$ is denoted by $\sM(A)$. We say that $\mu \in \sM(\R^n)$ is {\em concentrated on $A$} if $\mu(B) = 0$ for any $B \in \mathcal{B}(\R^n)$ with $B \subset \R^n \setminus A$, see for example \cite[p.~130]{cohn}. Moreover $\sM(A)$ can be identified with a vector subspace of $\sM(\R^n)$, viz.\ the vector subspace of measures concentrated on~$A$.

A function $g : \R^n \longrightarrow \R$ is called {\em positively one-homogeneous} if $g(a x) = a g(x)$ holds for every $a \geq 0$ and $x \in \R^n$, see \cite[p.~24]{schn2}. Given a function $f : S^{n-1} \longrightarrow \R$ there always exists a positively one-homogeneous extension to $\R^n$, viz.
\[
g : \R^n \longrightarrow \R,\quad g(x) = \left\{ \begin{array}{ll}
\|x\|_2 \, f(x / \|x\|_2), & x \in \R^n \setminus \{0\} \\[1em]
0, & x = 0
\end{array}\right.
\]
Clearly if $f$ is measurable (continuous), then $g$ is measurable (continuous). The set of positively one-homogeneous measurable functions from $\R^n$ to $\R$ bounded on $S^{n-1}$ is denoted by $M_b^{\rm ph}$.

For any two vectors $x, y \in \R^n$ the Euclidean inner product is denoted by $\langle x, y \rangle$, and the componentwise product by $xy$. The latter product is sometimes also called Hadamard product. For a set $A \subset \R^n$ and $x \in \R^n$ we write
\begin{equation}
\label{eq - comp prod vec set}
xA = Ax = \{ xy : y \in A \}
\end{equation}

\subsection{Basic properties}

For $\mu, \nu \in \sM(\R^n)$ we define $\mu \mc \nu \in \sM(\R^n)$, called the {\em multiplicative convolution}, as follows:
\[
(\mu \mc \nu)(A) = \int_{\R^n} \mu(dx) \int_{\R^n} \nu(dy) \; \one_A(xy), \quad A \in \mathcal{B}(\R^n)
\]
Note that $(\mu \mc \nu)(\R^n) = \mu(\R^n) \, \nu(\R^n)$. If $|\mu|(\{ x_j = 0 \}) = 0$ for each $j \in [n]$, then we can write, for each $A \in \mathcal{B}(\R^n)$,
\begin{eqnarray*}
(\mu \mc \nu)(A) & = & \int_{(\R \setminus \{0\})^n} \mu(dx) \int_{\R^n} \nu(dy) \; \one_{x^{-1}A} (y) \\
 & = & \int_{(\R \setminus \{0\})^n} \mu(dx) \, \nu(x^{-1}A)
\end{eqnarray*}

\begin{remark}
If $\xi$ and $\eta$ are independent random vectors in $\R^n$ with $\mu = \mathcal{L}(\xi)$ and $\nu = \mathcal{L}(\eta)$, then $\mu \mc \nu = \mathcal{L}(\xi \eta)$.
\end{remark}

\begin{example}
\label{ex - delta prod}
Let $\mu \in \sM(\R^n)$ and $x, y \in \R^n$. Then
\[
\delta_0 \mc \mu =  \mu(\R^n) \, \delta_0,\quad\quad \delta_x \mc \delta_y = \delta_{xy}
\]
If $x \neq 0$ and $A \in \mathcal{B}(\R^n)$, then
\[
\delta_x \mc \mu \, (A) = \mu(x^{-1}A)
\]
\end{example}

\begin{proposition}
\label{prop - multconv properties}
Let $\mu, \nu, \sigma \in \sM(\R^n)$ and $a \in \R$. We have
\begin{thnumber}
\item \label{prop - multconv properties commut} $\mu \mc \nu = \nu \mc \mu$
\item \label{prop - multconv properties const} $(a \mu) \mc \nu \,=\, \nu \mc (a \mu) \,=\, a (\mu \mc \nu)$
\item \label{prop - multconv properties distr} $(\mu + \nu) \mc \sigma = \mu \mc \sigma + \nu \mc \sigma$
\item \label{prop - multconv properties assoc} $(\mu \mc \nu) \mc \sigma = \mu \mc (\nu \mc \sigma)$
\end{thnumber}
\end{proposition}

\begin{remark}
\label{rema - mult conv jord decomp}
For $\mu, \nu \in \sM(\R^n)$ the minimality property of the Jordan decomposition (see e.g.~\cite[Satz~VII.1.12]{elstrodt}) yields 
\begin{eqnarray*}
(\mu \mc \nu)_+ & \leq & \mu_+ \mc \nu_+ + \mu_- \mc \nu_- \\
(\mu \mc \nu)_- & \leq & \mu_+ \mc \nu_- + \mu_- \mc \nu_+
\end{eqnarray*}
It follows that 
\[
| \mu \mc \nu | \,\leq\, |\mu| \mc |\nu|,
\]
and therefore
\[
\| \mu \mc \nu \| \,\leq\, \| \mu \| \, \| \nu \|
\]
Thus $\sM(\R^n)$ with multiplicative convolution as multiplication is a commutative Banach algebra over $\R$ with unit $\delta_{1_n}$. In particular multiplicative convolution is continuous (simultaneously in both arguments) with respect to the total variation norm.
\end{remark}

\begin{proposition}
\label{prop - multconv prod}
Let $m, n \geq 1$, $\mu, \nu \in \sM(\R^n)$, and $\rho, \sigma \in \sM(\R^m)$. Then
\[
(\mu \otimes \rho) \mc (\nu \otimes \sigma) \,=\, (\mu \mc \nu) \otimes (\rho \mc \sigma)
\]
\end{proposition}

\begin{remark}
\label{rema - prod multconv rv}
Notice that in the case of probability measures Proposition~\ref{prop - multconv prod} says that, for independent random vectors $\xi$ and $\eta$ in $\R^n$ and $\chi$ and $\zeta$ in $\R^m$ we have $(\xi, \chi) (\eta, \zeta) = (\xi \eta, \chi \zeta)$. This equation holds even almost surely (or everywhere), and it holds even if the set of all four vectors is not independent.
\end{remark}

\subsection{Orthogonal projections}

Next we introduce our notation for orthogonal projections on coordinate subspaces of $\R^n$. The standard basis vectors of $\R^n$ are denoted by $e_i$ for $i \in [n]$. For $E \subset [n]$ we define the vector subspace $H_E = \text{lin} \{ e_i : i \in E \}$ and the linear projection on $H_E$ by
\[
P_E : \R^n \longrightarrow \R^n,\quad P_E(e_j) = \left\{ \begin{array}{ll}
e_j, & j \in E \\[1em]
0, & j \notin E
\end{array}\right.
\]
When $E$ is a singleton, say $E = \{ i \}$ for some $i \in [n]$, we also write $H_i$ and $P_i$, respectively. For the projection of a set $A \subset \R^n$ we write
\[
P_E A \,=\, P_E(A) \,=\, \{ P_E(x) \,:\, x \in A \}
\]
For $E, F \subset [n]$ we clearly have $P_E \circ P_F = P_{E \cap F}$. Note also that for a positively one-homogeneous function $f$ on $\R^n$ the composition $f \circ P_E$ is positively one-homogeneous too.

For $\mu \in \sM(\R^n)$ and $E \subset [n]$ the marginal measure of $\mu$ on $H_E$ is given by $P_E(\mu) = \mu \circ P_E^{-1}$. Clearly $P_E$ is a linear operator on $\sM(\R^n)$. $\mu$ is concentrated on $H_E$ if and only if $P_E(\mu) = \mu$. Note also that
\begin{equation}
\label{eq - proj E pm}
(P_E(\mu))_+ \leq P_E(\mu_+),\quad (P_E(\mu))_- \leq P_E(\mu_-)
\end{equation}
holds for $E \subset [n]$ by the minimality property of the Jordan decomposition (see e.g.~\cite[Satz~VII.1.12]{elstrodt}), but equality does not hold in general in~\eqref{eq - proj E pm}.

The interaction of projection on a coordinate subspace with multiplicative convolution can be stated as follows.

\begin{proposition}
\label{prop - multconv proj}
Let $\mu, \nu \in \sM(\R^n)$ and let $E \subset [n]$. We have
\[
P_E(\mu \mc \nu) \,=\, P_E(\mu) \mc \nu \,=\, \mu \mc P_E(\nu) \,=\, P_E(\mu) \mc P_E(\nu)
\]
\end{proposition}

\begin{proof}
Given $x, y \in \R^n$, we have
\[
P_E(xy) \,=\, P_E(x) \,y \,=\, x P_E(y) \,=\, P_E(x) P_E(y)
\]
This implies the assertion.
\end{proof}

\begin{remark}
\label{rema - orth proj random vector}
Given a random vector $\xi$ in $\R^n$ and $E \subset [n]$ we can apply the projection also to $\xi$ so that $P_E(\xi)$ denotes the random vector with components $\xi_i$ for $i \in E$ and $0$ for $i \notin E$.
\end{remark}

\subsection{Integrability}

A measure $\mu \in \sM(\R^n)$ is called {\em integrable} if
\[
\int_{\R^n} \|x\|_2 \, |\mu|(dx) < \infty
\]
For $A \in \mathcal{B}(\R^n)$ the subspace of integrable members of $\sM(A)$ is denoted by $\sM^{\rm int}(A)$.

Clearly if $\mu$ is integrable, then $\mu_+$ and $\mu_-$ are integrable, as well are $P_E(\mu)$ (by~\eqref{eq - proj E pm}) for every $E \subset [n]$. Notice also that every $\mu \in \sM(\R^n)$ concentrated on a bounded set is integrable. In particular this is true if $\mu$ is concentrated on $S^{n-1}$.

\begin{proposition}
\label{prop - integr mult conv}
Let $\mu, \nu \in \sM(\R^n)$. If $\mu$ and $\nu$ are integrable, then $\mu \mc \nu$ is integrable.
\end{proposition}

\begin{proof}
By Remark~\ref{rema - mult conv jord decomp} we may assume that $\mu$ and $\nu$ are non-negative. In this case $\mu \mc \nu$ is non-negative. Noting that
\[
\|xy\|_2 \leq \|x\|_2 \, \|y\|_2,\quad x, y \in \R^n,
\]
we obtain
\begin{eqnarray*}
\int_{\R^n} \mu \mc \nu (dx) \, \|x\|_2 & = & \int_{\R^n} \mu(dx) \int_{\R^n} \nu(dy) \, \|xy\|_2 \\
 & \leq & \int_{\R^n} \mu(dx) \int_{\R^n} \nu(dy) \, \|x\|_2 \, \|y\|_2
\end{eqnarray*}
\end{proof}

\begin{proposition}
\label{prop - integr prod}
Let $\mu \in \sM(\R^m)$ and $\nu \in \sM(\R^n)$ where $m, n \geq 1$. If $\mu$ and $\nu$ are integrable, then $\mu \otimes \nu$ is integrable.
\end{proposition}

\begin{proof}
First note that $\mu \otimes \nu = \rho - \sigma$ where
\[
\rho = \mu_+ \otimes \nu_+ + \mu_- \otimes \nu_-,\quad \sigma = \mu_+ \otimes \nu_- + \mu_- \otimes \nu_+
\]
So $\rho$ and $\sigma$ are non-negative measures on $\R^{m + n}$. By the minimality of the Jordan decomposition
\[
(\mu \otimes \nu)_+ \leq \rho,\quad (\mu \otimes \nu)_- \leq \sigma
\]
Hence we may assume that $\mu$ and $\nu$ are non-negative. In this case $\mu \otimes \nu$ is non-negative and we have
\begin{eqnarray*}
\int_{\R^{m + n}} \mu \otimes \nu (dx, dy) \, \|(x,y)\|_2 & \leq & \int_{\R^m} \mu(dx) \int_{\R^n} \nu(dy) \, (\|x\|_2 + \|y\|_2)
\end{eqnarray*}
The right-hand side is finite if both $\mu$ and $\nu$ are integrable.
\end{proof}

\subsection{Radial projection}

Set $\R^n_0 = \R^n \setminus \{0\}$ for $n \geq 1$. Given $\mu \in \sM^{\rm int}(\R^n)$, we define $\tilde{\mu} \in \sM(\R^n)$ by
\[
\tilde{\mu}(A) = \int_{\R^n} \mu(dx) \, \|x\|_2 \, \one_A(x),\quad A \in \mathcal{B}(\R^n)
\]
Note that $\tilde{\mu}$ is well-defined because of the integrability of $\mu$. Moreover $\tilde{\mu}$ is always concentrated on $\R^n_0$ even if $\mu$ has an atom at the origin. 

Now define the map
\[
R : \R^n_0 \longrightarrow S^{n-1},\quad R(x) = \frac{x}{\|x\|_2}
\]
We set $\hat{\mu} = R(\tilde{\mu})$. Thus $\hat{\mu} \in \sM(S^{n-1})$ and $\hat{\mu}(B) = \tilde{\mu}(R^{-1}(B))$ for every $B \in \mathcal{B}(S^{n-1})$. The composite map $\mu \mapsto \hat{\mu}$ is called {\em radial projection} and denoted by
\[
P_S : \sM^{\rm int}(\R^n) \longrightarrow \sM(S^{n-1}),\quad P_S(\mu) = \hat{\mu}
\]
To summarise, the radial projection is, for $B \in \mathcal{B}(S^{n-1})$,
\[
 P_S(\mu)(B) \,=\, \int_{S^{n-1}} P_S(\mu)(du) \, \one_B(u) \,=\, \int_{\R^n_0} \mu(dx) \, \|x\|_2 \, \one_B \left(\frac{x}{\|x\|_2}\right)
\]
Clearly $P_S$ is a linear operator and $P_S(\mu) \in \sM^{\rm int}(\R^n)$. The relation $P_S(\mu) = \mu$ holds if and only if $\mu$ is concentrated on~$S^{n-1}$. In particular we always have $P_S(P_S(\mu)) = P_S(\mu)$. This property justifies the name "radial projection" and the notation $P_S$. For $\mu \geq 0$ we have $P_S(\mu) \geq 0$. It follows that for any $\mu, \nu \in \sM^{\rm int}(\R^n)$ with $\mu \geq \nu$, we have $P_S(\mu) \geq P_S(\nu)$. For a Dirac measure $\delta_x$ with any $x \in \R^n_0$ we have
\begin{equation}
\label{eq - ps delta}
P_S(\delta_x) \,=\, \|x\|_2 \; \delta_{x / \|x\|_2}
\end{equation}
It is easy to see that $P_S$ is not injective.

\begin{proposition}
\label{prop - multconv poshom}
Let $\mu \in \sM^{\rm int}(\R^n)$ and $f \in M_b^{\rm ph}$. Then
\[
\int_{\R^n} \mu(dx) \, f(x) = \int_{S^{n-1}} P_S(\mu)(du) \, f(u)
\]
and both integrals are finite.
\end{proposition}

\begin{proof}
By the assumptions the integral
\[
\int_{\R^n_0} \mu(dx) \, \|x\|_2 \, f\left(\frac{x}{\|x\|_2}\right)
\]
exists and is finite. Now, on the one hand
\[
\int_{\R^n_0} \mu(dx) \, \|x\|_2 \, f\left(\frac{x}{\|x\|_2}\right) \,=\, \int_{\R^n_0} \mu(dx) \, f(x) \,=\, \int_{\R^n} \mu(dx) \, f(x),
\]
and on the other hand
\begin{eqnarray*}
\int_{\R^n_0} \mu(dx) \, \|x\|_2 \, f\left(\frac{x}{\|x\|_2}\right) & = & \int_{\R^n_0} \tilde{\mu}(dx) \, (f \circ R)(x)\\
 & = & \int_{S^{n-1}} \hat{\mu}(du) \, f(u)\\
 & = & \int_{S^{n-1}} P_S(\mu)(du) \, f(u)
\end{eqnarray*}
\end{proof}

\begin{corollary}
\label{coro - ps total mass}
Let $\mu \in \sM^{\rm int}(\R^n)$. Then
\[
\int_{\R^n} \mu(dx) \, \|x\|_2 \,=\, P_S(\mu)(S^{n-1})
\]
In particular, $\mu = 0$ implies $P_S(\mu) = 0$, and $\mu > 0$ on $\R^n_0$ implies $P_S(\mu) > 0$.
\end{corollary}

\begin{proof}
Apply Proposition~\ref{prop - multconv poshom} with $f(x) = \|x\|_2$.
\end{proof}

\begin{remark}
Consider the special case of the radial projection $P_S(\mu)$ of a probability measure $\mu$. Let $\xi$ be a random vector in $\R^n$ that has distribution $\mu$ under a basic probability measure $\P$. Expectations with respect to this probability measure are denoted by $\E_{\P}$. Suppose that $\mu$ is integrable, in other words $\E_{\P} \| \xi \|_2 < \infty$. We may define another basic probability measure $\Q$ through the derivative
\[
\frac{d \Q}{d \P} = c^{-1} \| \xi \|_2,\quad c = \E_{\P} \| \xi \|_2
\]
Then
\[
(P_S(\mu))(B) = c \, \Q ( \xi / \| \xi \|_2 \in B),\quad B \in \mathcal{B}(S^{n-1})
\]
For given  $f \in M_b^{\rm ph}$ the formula in Proposition~\ref{prop - multconv poshom} now reads
\[
\E_{\P} f(\xi) = c \, \E_{\Q} f( \xi / \| \xi \|_2 )
\]
\end{remark}

\begin{remark}
\label{rema - prop jordan ps}
For $\mu \in \sM^{\rm int}(\R^n)$ the minimality of the Jordan decomposition yields
\[
(P_S(\mu))_+ \leq P_S(\mu_+),\quad (P_S(\mu))_- \leq P_S(\mu_-)
\]
Consequently
\[
|P_S(\mu)| \leq P_S(|\mu|)
\]

\end{remark}

For $\mu \in \sM(\R^n)$ and a bounded measurable function $g : \R^n \longrightarrow \R$ we denote
\[
g \odot \mu (A) \,=\, \int_{\R^n} \mu(dx) \, g(x) \one_A(x),\quad A \in \sB(\R^n)
\]
Clearly $g \odot \mu \in \sM(\R^n)$. If $\mu$ is integrable, then $g \odot \mu$ is integrable.

A function $g : \R^n \longrightarrow \R$ is called {\em positively zero-homogeneous} if $g(a x) = g(x)$ holds for every $a > 0$ and $x \in \R^n_0$, 

\begin{corollary}
\label{coro - multconv dens poshom}
Let $\mu \in \sM^{\rm int}(\R^n)$, and let $g$ be a positively zero-homogeneous bounded measurable function on $\R^n$. Then $P_S(g \odot \mu) = g \odot P_S(\mu)$.
\end{corollary}

\begin{proof}
Let $f \in M_b^{\rm ph}$. Then
\begin{eqnarray*}
\int_{S^{n-1}} P_S(g \odot \mu)(du) \, f(u) & = & \int_{\R^n} g \odot \mu (dx) \, f(x) \\
 & = & \int_{\R^n} \mu (dx) \, g(x) f(x) \\
 & = & \int_{S^{n-1}} P_S(\mu) (du) \, g(u) f(u) \\
 & = & \int_{S^{n-1}} g \odot P_S(\mu) (du) \, f(u)
\end{eqnarray*}
\end{proof}

The relation of $P_S$ and projections on coordinate subspaces is specified in the following proposition.

\begin{proposition}
\label{prop - proj-s proj-e}
Let $\mu \in \sM^{\rm int}(\R^n)$ and $E \subset [n]$. Then
\[
P_S \circ P_E \circ P_S \, (\mu) \,=\, P_S \circ P_E \, (\mu)
\]
\end{proposition}

\begin{proof}
Let $f \in M_b^{\rm ph}$. We have
\begin{eqnarray}
\int_{S^{n-1}} d(P_S \circ P_E \circ P_S) (\mu) \; f & = & \int_{\R^n} d(P_E \circ P_S)(\mu) \; f \label{prop - proj-s proj-e 1}\\
 & = & \int_{\R^n} d P_S (\mu) \; f \circ P_E \nonumber\\
 & = & \int_{\R^n} d \mu \; f \circ P_E \label{prop - proj-s proj-e 3}\\
 & = & \int_{\R^n} d P_E(\mu) \; f \nonumber\\
 & = & \int_{S^{n-1}} d(P_S \circ P_E)(\mu) \; f \label{prop - proj-s proj-e 5}
\end{eqnarray}
where in \eqref{prop - proj-s proj-e 1} and \eqref{prop - proj-s proj-e 5} Proposition~\ref{prop - multconv poshom} is applied to $f$, and in \eqref{prop - proj-s proj-e 3} the same Proposition is applied to $f \circ P_E$.
\end{proof}

Multiplicative convolution of measures can be combined with the map $P_S$ leading to the following results.

\begin{proposition}
\label{prop - multconv ps}
For $\mu, \nu \in \sM^{\rm int}(\R^n)$ we have
\[
P_S(\mu \mc \nu) \,=\, P_S \big(P_S(\mu) \mc \nu \big) \,=\, P_S \big(\mu \mc P_S(\nu) \big) \,=\, P_S \big(P_S(\mu) \mc P_S(\nu) \big)
\]
\end{proposition}

\begin{proof}
We show the first identity. The other relations then easily follow.

Let $f \in M_b^{\rm ph}$. For each $x \in \R^n$ the map $y \mapsto f(xy)$ is positively one-homogeneous, measurable, and bounded on~$S^{n-1}$. Hence the function
\[
g : \R^n \longrightarrow \R,\quad g(x) = \int_{\R^n} \nu(dy) \, f(xy)
\]
is well-defined. Moreover $g$ is positively one-homogeneous and measurable. Since $\mu \mc \nu$ is integrable by Proposition~\ref{prop - integr mult conv}, we have
\begin{eqnarray}
\int_{S^{n-1}} P_S(\mu \mc \nu)(du) \, f(u) & = & \int_{\R^n} (\mu \mc \nu)(dx) \, f(x) \label{prop - multconv ps dom}\\
 & = & \int_{\R^n} \mu(dx) \int_{\R^n} \nu(dy) \, f(xy) \nonumber\\
 & = & \int_{\R^n} \mu(dx) \, g(x) \nonumber\\
 & = & \int_{S^{n-1}} P_S(\mu)(du) \, g(u) \label{prop - multconv ps g}\\
 & = & \int_{S^{n-1}} P_S(\mu)(du) \int_{\R^n} \nu(dy) \, f(uy) \nonumber\\
 & = & \int_{\R^n} \left(P_S(\mu) \mc \nu \right)(dx) \, f(x) \nonumber\\
 & = & \int_{S^{n-1}} P_S \left( P_S(\mu) \mc \nu \right)(du) \, f(u) \label{prop - multconv ps f}
\end{eqnarray}
where \eqref{prop - multconv ps dom}, \eqref{prop - multconv ps g}, and \eqref{prop - multconv ps f} follow from Proposition~\ref{prop - multconv poshom}.
\end{proof}

For $\mu, \nu \in \sM^{\rm int}(\R^n)$ we define a measure on $S^{n-1}$ by
\[
\mu \mcs \nu = P_S(\mu \mc \nu)
\]
In particular note that the restriction of $\mcs$ to $\sM(S^{n-1})$ is a binary operation on $\sM(S^{n-1})$, which is a motivation for this definition on its own.

\begin{example}
\label{ex - delta prod s}
Let $x, y \in \R^n$. Using Example~\ref{ex - delta prod} and formula~\eqref{eq - ps delta} we have
\[
\delta_x \mcs \delta_y \,=\, \|xy\|_2 \; \delta_{xy / \|xy\|_2}
\]
\end{example}

\begin{corollary}
\label{coro - radial proj prod}
For $\mu, \nu \in \sM^{\rm int}(\R^n)$ we have
\[
\mu \mcs \nu \,=\, P_S(\mu \mcs \nu) \,=\, P_S(\mu) \mcs \nu \,=\, \mu \mcs P_S(\nu) \,=\, P_S(\mu) \mcs P_S(\nu)
\]
\end{corollary}

\begin{proof}
This follows from Proposition \ref{prop - multconv ps}.
\end{proof}

\begin{proposition}
\label{prop - star prod properties}
Let $\mu, \nu, \sigma \in \sM^{\rm int}(\R^n)$ and $a \in \R$. Then
\begin{thnumber}
\item $\mu \mcs \nu = \nu \mcs \mu$
\item $(a \mu) \mcs \nu \,=\, \mu \mcs (a \nu) \,=\, a (\mu \mcs \nu)$
\item $(\mu + \nu) \mcs \sigma = \mu \mcs \sigma + \nu \mcs \sigma$
\item $(\mu \mcs \nu) \mcs \sigma = \mu \mcs (\nu \mcs \sigma)$
\end{thnumber}
\end{proposition}

The identities listed in Proposition~\ref{prop - star prod properties} are interesting specifically for $\sM(S^{n-1})$.

\begin{proposition}
\label{prop - star prod banach}
Let $\mu, \nu \in \sM(S^{n-1})$. Then
\[
\| \mu \mcs \nu \| \leq \| \mu \| \, \| \nu \|
\]
\end{proposition}

\begin{proof}
Recall from Remark~\ref{rema - mult conv jord decomp} that $|\mu \mc \nu| \leq |\mu| \mc |\nu|$. It follows that
\begin{equation}
\label{eq - inequ ps prod}
P_S(|\mu \mc \nu|) \leq P_S(|\mu| \mc |\nu|)
\end{equation}
Consequently
\begin{eqnarray}
\| \mu \mcs \nu \| & = & |P_S(\mu \mc \nu)| \, (S^{n-1}) \nonumber\\
 & \leq & P_S(|\mu \mc \nu|) \, (S^{n-1}) \label{eq - ps abs}\\
 & \leq & P_S(|\mu| \mc |\nu|) \, (S^{n-1}) \label{eq - ps prod abs}\\
 & = & \int_{\R^n} |\mu| \mc |\nu| (dx) \, \|x\|_2 \label{eq - ps rn}\\
 & = & \int_{S^{n-1}} |\mu| (dx) \int_{S^{n-1}} |\nu| (dy) \, \|xy\|_2 \nonumber\\
 & \leq & \int_{S^{n-1}} |\mu| (dx) \int_{S^{n-1}} |\nu| (dy) \;\leq\; \| \mu \| \, \| \nu \| \nonumber
\end{eqnarray}
where \eqref{eq - ps abs} follows from Remark~\ref{rema - prop jordan ps}, \eqref{eq - ps prod abs} from \eqref{eq - inequ ps prod}, and \eqref{eq - ps rn} from Corollary~\ref{coro - ps total mass}.
\end{proof}

\begin{remark}
Propositions \ref{prop - star prod properties} and \ref{prop - star prod banach} imply that $\sM(S^{n-1})$ with the operation $\mcs$ as multiplication is a commutative Banach algebra with unit $P_S(\delta_{1_n})$. Therefore as is the case with the operation~$\mc$ on $\sM(\R^n)$ the operation $\mcs$ is continuous (simultaneously in both arguments) with respect to the total variation norm.
\end{remark}

We also would like to define the projection of any measure on $S^{n-1}$ onto a measure on a {\em coordinate subsphere} $S_E = H_E \cap S^{n-1}$ where $E \subset [n]$. Note that $S_E$ is the unit sphere in the coordinate subspace $H_E$. To this end we define, for $E \subset [n]$, the operator $P^S_E = P_S \circ P_E$, mapping each member of $\sM^{\rm int}(\R^n)$ to an element of $\sM(S_E)$. So the most interesting subset of the domain of $P^S_E$ is $\sM(S^{n-1})$. The same map is used in~\cite[Section~3]{weil82} for the study of generalised zonoids. Linearity of $P^S_E$ follows from the linearity of $P_E$ and $P_S$. Proposition~\ref{prop - proj-s proj-e} yields $P^S_E \circ P^S_E = P^S_E$, confirming the projection property. Moreover for $E, F \subset [n]$ Proposition~\ref{prop - proj-s proj-e} yields
\[
P^S_E \circ P^S_F \,=\, P^S_F \circ P^S_E \,=\, P^S_{E \cap F}
\]

\begin{proposition}
\label{prop - proj s prod}
For $\mu, \nu \in \sM^{\rm int}(\R^n)$ we have
\[
P^S_E(\mu \mcs \nu) \,=\, P^S_E(\mu) \mcs \nu \,=\, \mu \mcs  P^S_E(\nu) \,=\, P^S_E(\mu) \mcs P^S_E(\nu)
\]
\end{proposition}

\begin{proof}
The first equation is shown as follows:
\begin{eqnarray}
P^S_E (\mu \mcs \nu) & = & P_S P_E P_S (\mu \mc \nu) \nonumber\\
 & = & P_S P_E (\mu \mc \nu) \label{prop - proj s prod 2}\\
 & = & P_S (P_E(\mu) \mc \nu) \label{prop - proj s prod 3}\\
 & = & P_S (P^S_E(\mu) \mc \nu) \label{prop - proj s prod 4}\\
 & = & P^S_E(\mu) \mcs \nu \nonumber
\end{eqnarray}
where \eqref{prop - proj s prod 2} follows from Proposition~\ref{prop - proj-s proj-e}, \eqref{prop - proj s prod 3} follows from Proposition~\ref{prop - multconv proj}, and \eqref{prop - proj s prod 4} follows from Proposition~\ref{prop - multconv ps}. The other equations can be shown similarly.
\end{proof}

Proposition~\ref{prop - proj s prod} is interesting mainly for $\mu$ and $\nu$ concentrated on~$S^{n-1}$.

\section{Decomposition of measures}
\label{sec - decomp measures}

We now introduce partitions of $\R^n$ and $S^{n-1}$ and corresponding decompositions of elements of $\sM(\R^n)$ and $\sM(S^{n-1})$ by restricting them to the members of those partitions. For $n \geq 1$ and $E \subset [n]$ define
\begin{eqnarray*}
A_E & = & \big\{x \in \R^n \;:\; \forall i \in [n] \;\; x_i \neq 0 \,\Longleftrightarrow\, i \in E \big\},\\
B_E & = & A_E \cap S^{n-1}
\end{eqnarray*}
In particular $A_{\O} = \{0\}$ and $B_{\O} = \O$. The family $\{ A_E : E \subset [n] \}$ is a partition of $\R^n$ and $\{ B_E : E \subset [n] \}$ is a partition of $S^{n-1}$. We have, for $E \subset [n]$,
\[
H_E = \bigcup_{ F \subset E} A_F,\quad S_E = \bigcup_{ F \subset E} B_F
\]
For $n \geq 1$ we define
\begin{equation}
\label{def - m}
M^n = \big\{ x \in \R_+^n \;:\;  \|x\|_1 \leq 1 \big\}
\end{equation}

For $\mu \in \sM(\R^n)$ we introduce the {\em coordinate decomposition}
\begin{equation}
\label{eq - sign meas decomp}
\mu = \sum_{E \subset [n]} \mu_E
\end{equation}
where $\mu_E(A) = \mu(A \cap A_E)$ for $A \in \mathcal{B}(\R^n)$.  In the special case where $\mu$ is concentrated on $S^{n-1}$ we obtain $\mu_E(A) = \mu(A \cap B_E)$. 

For $\sE \subset \sP_n$ and $A = \bigcup_{E \in \sE} A_E$ we also write $\sM(\sE)$ for $\sM(A)$ and $\sM^{\rm int}(\sE)$ for $\sM^{\rm int}(A)$. Moreover for $\sE \subset \sP_n$ with $\O \notin \sE$ and $B = \bigcup_{E \in \sE} B_E$, we write $\sM^{\rm sph}(\sE)$ for $\sM(B)$. Notice that the superscript "sph" distinguishes the measures concentrated on the sphere from that on the entire~$\R^n$.

We say that $\mu \in \sM(\R^n)$ is {\em of order} $E$ if $|\mu|(A_E) > 0$ and $|\mu|(A_F) = 0$ for all $F \subset [n]$ with $F \neq E$. Specifically if $\mu$ is of order $\O$, then it is an atom at the origin. Obviously the order is not defined for all $\mu \in \sM(\R^n)$. Rather \eqref{eq - sign meas decomp} is the decomposition of $\mu$ into terms of various orders; for each $E \subset [n]$ (including $E = \O$) the term $\mu_E$ is either zero or of order $E$. Moreover the {\em degree of} $\mu$ is the maximum integer $k \in \{0, 1, \ldots, n\}$ such that there is $E \subset [n]$ with $|E| = k$ and $\mu_E \neq 0$.

The decomposition induces mappings $R_E(\mu) = \mu_E$ for $E \subset [n]$. These are linear projections, that is, $R_E$ is linear on $\sM(\R^n)$ and $R_E \circ R_E = R_E$. Clearly we have
\[
(\mu_E)_+ = (\mu_+)_E,\quad (\mu_E)_- = (\mu_-)_E,
\]
and therefore also
\begin{equation}
\label{eq - decomp meas pm}
|\mu_E| = (|\mu|)_E
\end{equation}
Consequently if $\mu$ is integrable, then also $\mu_E$ is integrable.

\begin{remark}
Given a probability measure $\mu$ on $\R^n$, let $\xi$ be a random vector such that $\mathcal{L}(\xi) = \mu$. For  $E \subset [n]$ and $A \in \mathcal{B}(\R^n)$ we have
\[
\mu_E(A) \,=\, \mu(A \cap A_E) \,=\, \P(\xi \in A \cap A_E) \,=\, \P(\xi \in A \,|\, \xi \in A_E) \, \P(\xi \in A_E)
\]
So the decomposition~\eqref{eq - sign meas decomp} can be written as
\[
\P(\xi \in A) = \sum_{E \subset [n]} \P(\xi \in A \,|\, \xi \in A_E) \, \P(\xi \in A_E)
\]
\end{remark}

\begin{remark}
\label{rema - coord decomposition random vectors}
Let $\mu$ be a probability measure on $\R^n$ and $\xi$ a random vector with $\mathcal{L}(\xi) = \mu$. We may decompose $\xi$ as follows:
\begin{equation}
\label{eq - random vector coord decomp}
\xi = \sum_{E \subset [n]} \xi_E, \quad \xi_E = \xi \, \one_{A_E}(\xi)
\end{equation}
Now suppose that $\mu$ has degree $k$. Then there is $E \subset [n]$ such that $k = |E|$ and $P(\xi \in A_E) > 0$; however $P(\xi \in A_F) = 0$ holds for all $F \subset [n]$ with $k < |F|$. However \eqref{eq - random vector coord decomp} does not exactly coincide with the coordinate decomposition of measures. In fact for $E \subset [n]$ we have
\[
\mathcal{L}(\xi_E) \,=\, \mu_E + \P(\xi \notin A_E) \delta_0
\]
The decomposition~\eqref{eq - random vector coord decomp} is used in Lemmas 3.4 and 3.8 in~\cite{mol:nag:21}.
\end{remark}

We next consider the connection between coordinate decomposition and multiplicative convolution.

\begin{proposition}
\label{prop - decomp prod}
Let $E, F \subset [n]$ and $\mu, \nu \in \sM(\R^n)$ such that $\mu$ is of order $E$ and $\nu$ is of order $F$. Then $\mu \mc \nu$ is of order $E \cap F$ or equal to $0$. If both $\mu, \nu \geq 0$, then $\mu \mc \nu$ is of order $E \cap F$.
\end{proposition}

\begin{proof}
First assume that $\mu, \nu \geq 0$. Note that if $x \in A_E$ and $y \in A_F$, then $xy \in A_{E \cap F}$. Thus
\[
\one_{A_E}(x) \one_{A_F}(y) \leq \one_{A_{E \cap F}}(xy),\quad x,y \in \R^n
\]
Consequently
\begin{eqnarray*}
(\mu \mc \nu)(\R^n) & = & \int_{\R^n} \mu(dx) \int_{\R^n} \nu(dy) \; \one_{A_E}(x) \one_{A_F}(y) \\
 & \leq & \int_{\R^n} \mu(dx) \int_{\R^n} \nu(dy) \; \one_{A_{E \cap F}}(xy) \\
 & = & (\mu \mc \nu)(A_{E \cap F})
\end{eqnarray*}
Hence
\[
(\mu \mc \nu)(A_{E \cap F}) \,=\, (\mu \mc \nu)(\R^n) \,=\,  \mu(\R^n) \nu(\R^n) > 0
\]
This shows the second statement.

For general $\mu$ and $\nu$ set $\rho = \mu \mc \nu$. Then $\rho = \sigma - \tau$ where
\[
\sigma = \mu_+ \mc \nu_+ + \mu_- \mc \nu_-,\quad \tau = \mu_+ \mc \nu_- + \mu_- \mc \nu_+
\]
By the assumptions $\mu_+$ and $\mu_-$ each is of order $E$ or equal to zero, whereas $\nu_+$ and $\nu_-$ each is of order $F$ or equal to zero. It follows from the first part of the proof that $\sigma$ is of order $E \cap F$ or equal to zero, and the same holds for $\tau$.
The minimality property of the Jordan decomposition yields
\[
\rho_+ \leq \sigma,\quad \rho_- \leq \tau
\]
Therefore $\rho_+$ and $\rho_-$ are of order $E \cap F$ or equal to $0$, and the same follows for $\rho$.
\end{proof}

In particular Proposition~\ref{prop - decomp prod} yields the formula $(\mu \mc \nu)(A_{E\cap F}) = \mu(A_E) \, \nu(A_F)$.

\begin{corollary}
\label{cor - multconv decomp}
Given $\mu, \nu \in \sM(\R^n)$ we have
\[
\mu \mc \nu = \sum_{G \subset [n]} (\mu \mc \nu)_G,\quad (\mu \mc \nu)_G = \sum_{E, F \subset [n] \atop E \cap F = G} \mu_E \mc \nu_F
\]
If $k$ is the degree of $\mu$ and $l$ is the degree of $\nu$, then the degree of $\mu \mc \nu$ is smaller or equal than the minimum of $k$ and $l$.
\end{corollary}

\begin{corollary}
\label{coro - mult conv decomp note}
Let $n \geq 1$, and for $j \in \{1, 2\}$ let $\sE_j \subset \sP_n$ and $\mu_j \in \sM(\sE_j)$. Then $\mu_1 \mc \mu_2 \in \sM(\sE)$ where
\[
\sE = \big\{ E_1 \cap E_2 \,:\, E_1 \in \sE_1,\; E_2 \in \sE_2 \big\}
\]
\end{corollary}

We turn to the relation between $P_S$ and the decomposition of measures. Also recall that $R_{\O}(\nu) = 0$ for every $\nu \in \sM(S^{n-1})$.

\begin{proposition}
\label{prop - proj S decomp E}
For every $\mu \in \sM^{\rm int}(\R^n)$ and $E \subset [n]$ we have
\[
(P_S(\mu))_E = P_S(\mu_E)
\]
\end{proposition}

\begin{proof}
Let $f \in M_b^{\rm ph}$. Then
\begin{eqnarray}
\int_{S^{n-1}} R_E P_S(\mu) (du) \, f(u) & = & \int_{S^{n-1}} P_S(\mu) (du) \, \one_{A_E}(u) \, f(u) \nonumber \\
 & = & \int_{\R^n} \mu (dx) \, \one_{A_E}(x) \, f(x) \label{prop - proj S decomp E 1} \\
 & = & \int_{\R^n} \mu_E (dx) \, f(x) \nonumber \\
 & = & \int_{S^{n-1}} (P_S(\mu_E))(du) \, f(u) \label{prop - proj S decomp E 2}
\end{eqnarray}
where \eqref{prop - proj S decomp E 1} and \eqref{prop - proj S decomp E 2} follow from Proposition~\ref{prop - multconv poshom} since $\one_{A_E} f \in M_b^{\rm ph}$.
\end{proof}

\begin{remark}
Note that if $\mu \in \sM^{\rm int}(\R^n)$ is of order $E$ for some $E \subset [n]$, then $P_S(\mu)$ is of order $E$ or equal to zero. Consequently the degree of $P_S(\mu)$ is always smaller or equal to the degree of $\mu$. It also follows that, if $\mu$ is concentrated on $H_E$, then $P_S(\mu)$ is concentrated on $S_E$.
\end{remark}

\begin{proposition}
\label{prop - decomp star prod}
Let $E, F \subset [n]$ and $\mu, \nu \in \sM^{\rm int}(\R^n)$ such that $\mu$ is of order $E$ and $\nu$ is of order $F$. Then $\mu \mcs \nu$ is of order $E \cap F$ or equal to $0$. If $\mu, \nu \geq 0$ and $E \cap F \neq \O$, then $\mu \mcs \nu$ is of order $E \cap F$.
\end{proposition}

\begin{proof}
This is a consequence of Propositions \ref{prop - decomp prod} and \ref{prop - proj S decomp E} and Corollary~\ref{coro - ps total mass}.
\end{proof}

\begin{corollary}
\label{cor - decomp star prod}
Given $\mu, \nu \in \sM^{\rm int}(\R^n)$ we have
\[
\mu \mcs \nu = \sum_{G \subset [n] \atop G \neq \O} (\mu \mcs \nu)_G,\quad (\mu \mcs \nu)_G = \sum_{E, F \subset [n] \atop E \cap F = G} \mu_E \mcs \nu_F
\]
If $k$ is the degree of $\mu$ and $l$ is the degree of $\nu$, then the degree of $\mu \mcs \nu$ is smaller or equal than the minimum of $k$ and $l$.
\end{corollary}

\begin{proof}
This follows from Proposition~\ref{prop - decomp star prod}.
\end{proof}

\begin{corollary}
\label{coro - mult conv decomp note sphere}
Let $n \geq 1$, and for $j \in \{1, 2\}$ let $\sE_j \subset \sP_n \setminus \{\O\}$ and $\mu_j \in \sM^{\rm int}(\sE_j)$. Then $\mu_1 \mcs \mu_2 \in \sM^{\rm sph}(\sE)$ where
\[
\sE = \big\{ E_1 \cap E_2 \,:\, E_1 \in \sE_1,\; E_2 \in \sE_2 \big\} \setminus \{\O\}
\]
\end{corollary}

In order to see the interrelation between coordinate decomposition and projections on coordinate subspaces first note that $P_E(x) = x$ holds for $x \in H_E$, in particular for $x \in A_E$. In fact we have $P_F(A_E) = A_{E \cap F}$ for every $E, F \subset [n]$.

\begin{proposition}
\label{prop - proj decomp measure}
Let $\mu \in \sM(\R^n)$ and $E, F \subset [n]$. Then
\begin{thnumber}
\item \label{prop - proj decomp measure proj} $P_{E}(\mu_E) = \mu_E$
\item \label{prop - proj decomp measure simpl} $P_{E \cap F}(\mu_E) = P_F(\mu_E)$
\item \label{prop - proj decomp measure term} $P_F(\mu_E)$ either equals zero or is of order $E \cap F$.
\item \label{prop - proj decomp measure nonneg} If $\mu \geq 0$ and $\mu_E > 0$, then $P_F(\mu_E)$ is of order $E \cap F$.
\end{thnumber}
\end{proposition}

\begin{proof}
To see \eqref{prop - proj decomp measure proj} note that, for $A \in \mathcal{B}(\R^n)$,
\[
(P_E(\mu_E))(A) \,=\, \mu(P_E^{-1}(A) \cap A_E) \,=\, \mu(A \cap A_E) \,=\, \mu_E(A)
\]

To see \eqref{prop - proj decomp measure simpl} let $A \in \mathcal{B}(\R^n)$ and notice that
\begin{eqnarray}
(P_{E \cap F}(\mu_E))(A) & = & \mu_E(P_{E \cap F}^{-1}(A)) \nonumber\\
 & = & \mu_E(P_E^{-1} P_F^{-1}(A)) \nonumber\\
 & = & (P_E(\mu_E))(P_F^{-1}(A)) \nonumber\\
 & = & \mu_E(P_F^{-1}(A)) \label{prop - proj decomp measure simpl first}\\
 & = & (P_F(\mu_E))(A) \nonumber
\end{eqnarray}
where \eqref{prop - proj decomp measure simpl first} follows from~\eqref{prop - proj decomp measure proj}.

We next prove \eqref{prop - proj decomp measure nonneg}. So assume that $\mu$ is non-negative. Let $G \subset [n]$ with $G \neq E \cap F$. Then
\begin{eqnarray*}
(P_F^{-1}(A_G)) \cap A_E & \subset & (P_F^{-1}(A_G)) \cap P_F^{-1}(P_F(A_E))\\
 & = & P_F^{-1}(A_G \cap A_{E \cap F})\\
 & = & \O
\end{eqnarray*}
Hence
\[
(P_F(\mu_E))_G (\R^n) \,=\, \mu(P_F^{-1}(A_G)) \cap A_E) \,=\, 0
\]
If $\mu_E(\R^n) > 0$, then $(P_F(\mu_E))(\R^n) = \mu(A_E) > 0$. Consequently
 \[
(P_F(\mu_E))_{E \cap F} (\R^n) \,>\, 0
\]
Therefore $P_F(\mu_E)$ is of order $E \cap F$.

Finally to prove \eqref{prop - proj decomp measure term} note that
\begin{equation}
\label{eq - proj signed jordan}
P_F(\mu_E) = P_F((\mu_+)_E) - P_F((\mu_-)_E)
\end{equation}
From statement \eqref{prop - proj decomp measure nonneg} we know that each term on the right-hand side of \eqref{eq - proj signed jordan} is either zero or of order $E \cap F$.
\end{proof}

Loosely speaking, in applying $P_F$ to $\mu_E$ the coordinates $E \setminus F$ are integrated over, the coordinates $E \cap F$ remain unchanged, the coordinates $F \setminus E$ have been zero already before, and the coordinates $[n] \setminus (E \cup F)$ are integrated over though having been zero.

Now we define
\[
\sigma_0^n = \sum_{s \in \{1, 2 \}^n} \prod_{j = 1}^n (-1)^{s_j} \; \delta_s = \sum_{s \in \{1, 2 \}^n} (-1)^{\sum_{j = 1}^n s_j} \; \delta_s
\]
Clearly
\[
\sigma_0^n \,\in\, \sM^{\rm int}((0,\infty)^n) \,\subset\, \sM^{\rm int}(A_{[n]})
\]

\begin{remark}
Note that $\sigma_0^n$ is a transformation of $\delta_{[n]}$; the latter notation is defined in Section~\ref{subsec - symmetry decomposition}. In this transformation the location of atoms is translated and scaled and the weight of every atom multiplied with $2^n$. However $\delta_{[n]}$ itself does not have the properties of $\sigma_0^n$ in Lemma~\ref{le - deltaj even odd class} with respect to the index transformation. Thus we need both measures.
\end{remark}

\begin{lemma}
\label{le - proj sigma0}
Let $n \geq 1$ and $E \subsetneqq [n]$. Then $P_E(\sigma_0^n) = 0$.
\end{lemma}

\begin{proof}
It is enough to show the claim for $|E| = n-1$. Without loss of generality assume that $E = [n - 1]$. Note that
\begin{equation}
\label{eq - sigma0n repr}
\sigma_0^n \,=\, \sum_{s_j \in \{1, 2 \} \atop j = 1,\ldots n-1} \prod_{j = 1}^{n - 1} (-1)^{s_j} \sum_{s_n \in \{1, 2 \}} (-1)^{s_n} \, \delta_s
\end{equation}
It follows that $P_{[n - 1]}(\sigma_0^n) = 0$ by linearity.
\end{proof}

\begin{corollary}
\label{coro - sigma0 n term}
Let $\nu \in \sM(\R^n)$. Then:
\begin{thnumber}
\item \label{coro - sigma0 n term rn} $\nu \mc \sigma_0^n = \nu_{[n]} \mc \sigma_0^n$ 
\item \label{coro - sigma0 n term sphere} If $\nu$ is integrable, then $\nu \mcs \sigma_0^n = \nu_{[n]} \mcs \sigma_0^n$ 
\end{thnumber}
\end{corollary}

\begin{proof}
To see \eqref{coro - sigma0 n term rn} note that
\begin{eqnarray*}
\nu \mc \sigma_0^n & = & \sum_{E \subset [n]} (\nu_E \mc \sigma_0^n)\\
 & = & \sum_{E \subset [n]} (P_E(\nu_E) \mc \sigma_0^n)\\
 & = & \sum_{E \subset [n]} (\nu_E \mc P_E(\sigma_0^n))\\
 & = & \nu_{[n]} \mc \sigma_0^n
\end{eqnarray*}

Statement~\eqref{coro - sigma0 n term sphere} follows by applying $P_S$ on both sides.
\end{proof}

\section{Reflections}
\label{sec - reflections}

\subsection{Basic definitions}

For $E \subset [n]$ we define the {\em multiple reflection}
\[
T_E : \R^n \longrightarrow \R^n,\quad (T_E(x))_j = \left\{ \begin{array}{ll}
-x_j, & j \in E \\[1em]
x_j, & j \notin E
\end{array}\right.
\]
In particular $T_{\O}$ is the identity map on~$\R^n$. For $j \in [n]$ the map $T_{\{j\}}$ is the reflection at the hyperplane perpendicular to~$e_j$. We define $\sT = \{ T_E : E \subset [n] \}$ which contains all linear maps on $\R^n$ represented by diagonal matrices with diagonal entries in $\{ 1, -1\}$. Notice that $(\sT, \circ)$ where $\circ$ denotes the composition of linear maps is a finite Abelian group with identity element $T_{\O}$.

For $E, F \subset [n]$ the {\em symmetric difference} is defined as usual by
\[
E \Delta F \,=\, (E \setminus F) \cup (F \setminus E)
\]
For $E, F, G \subset [n]$ we have
\begin{thnumber}
\item $E \Delta \O = E$,
\item \label{eq - boolean group idem} $E \Delta E = \O$,
\item $(E \Delta F) \Delta G = E \Delta (F \Delta G)$,
\item $E \Delta F = F \Delta E$
\end{thnumber}
This means that $(\sP_n, \Delta)$ is a finite Abelian group with identity element $\O$. It is a Boolean group since each element is the inverse of itself as stated in~\eqref{eq - boolean group idem}. It is clearly isomorphic to the group $(\sT, \circ)$, a group isomorphism being given by $E \mapsto T_E$. In the context of symmetries of measures on $\R^n$ some important subgroups are: $\{ \O \}$ (no non-trivial symmetry), $\{\O, [n]\}$ (only reflection at the origin), and $\sP_n$ (all multiple reflections). Given any $\sF \subset \sP_n$, the smallest subgroup $\sG$ of $\sP_n$ such that $\sF \subset \sG$ is called the {\it group generated by}~$\sF$, that is the intersection of all subgroups of $\sP_n$ containing~$\sF$. For definiteness we agree that the group generated by $\O$ is $\{ \O \}$.

\begin{proposition}
\label{prop - f and g generated}
Let $\sF \subset \sP_n$ such that $\sF \neq \O$, and $\sG$ the group generated by $\sF$. Then for each $G \in \sG$ there are $m \geq 1$ and $F_1, \ldots, F_m \in \sF$ such that $G = F_1 \Delta \cdots \Delta F_m$. If $G \neq \O$, then $F_1, \ldots, F_m$ can chosen to be distinct.
\end{proposition}

\begin{proof}
This follows from the fact that $\sG$ is a finite subgroup in which each element can be written as a finite product of the generating set $\sF$, see for example \cite[Proposition~3.3 and Corollary~3.4]{grillet}. Since the group is Abelian and every element is the inverse of itself, we can assume that $F_1, \ldots, F_m$ are distinct, provided $G \neq \O$.
\end{proof}

For $E, F, G \subset [n]$ we also have the distribution law
\begin{equation}
\label{eq - distr law sets}
(E \Delta F) \cap G \,=\, (E \cap G) \Delta (F \cap G)
\end{equation}
For $\sF \subset \sP_n$ and $E \subset [n]$ we use the notation
\[
\sF |_E \,=\, \{ F \cap E \,:\, F \in \sF \}
\]
Consequently $\sF |_E \subset \sP(E)$ .

\begin{remark}
\label{rema - proj gen commute}
Let $\sF \subset \sP_n$ and $E \subset [n]$. Let $\sG$ be the group generated by $\sF$. Then the distribution law \eqref{eq - distr law sets} implies that the group generated by $\sF |_E$ is $\sG |_E$. This also implies that, if $\sF$ is a group, then $\sF |_E$ is a group.
\end{remark}

\subsection{Invariant measures}

In this section we consider the action of multiple reflections on measures.

\begin{proposition}
\label{rema - mult refl commute}
For $\mu \in \sM(\R^n)$ and $E, F \subset [n]$ we have
\begin{thnumber}
\item $T_F P_E (\mu) = P_E T_F (\mu)$
\item $T_F R_E (\mu) = R_E T_F (\mu)$
\item If $\mu$ is of order $E$, then $T_F(\mu)$ is of order $E$.
\item The degrees of $\mu$ and $T_F(\mu)$ are equal.
\item If $\mu$ is integrable, then $T_F(\mu)$ is integrable.
\item If $\mu$ is integrable, then $T_F P_S (\mu) = P_S T_F (\mu)$ and $T_F P^S_E (\mu) = P^S_E T_F (\mu)$.
\end{thnumber}
\end{proposition}

\begin{proposition}
\label{prop - refl conv rn}
For $\mu, \nu \in \sM(\R^n)$ and $E \subset [n]$ we have
\begin{equation}
\label{eq - refl conv}
T_E(\mu \mc \nu) \,=\, T_E(\mu) \mc \nu \,=\, \mu \mc T_E(\nu)
\end{equation}
If $\mu$ and $\nu$ are integrable, then
\begin{equation}
\label{eq - refl conv sphere}
T_E(\mu \mcs \nu) \,=\, T_E(\mu) \mcs \nu \,=\, \mu \mcs T_E(\nu)
\end{equation}
\end{proposition}

\begin{proof}
\eqref{eq - refl conv} is straightforward to show. Then \eqref{eq - refl conv sphere} is a consequence of \eqref{eq - refl conv} and Proposition~\ref{rema - mult refl commute}.
\end{proof}

Notice that \eqref{eq - refl conv} is not equal to $T_E(\mu) \mc T_E(\nu)$ in general.

For given $F \subset [n]$ a set $A \in \sB(\R^n)$ is called {\em invariant under} $T_F$ if $T_F(A) = A$. Given $F \subset [n]$, and $A \in \sB(\R^n)$ invariant under $T_F$, $\mu \in \sM(A)$ is called {\em even under $T_F$} if $T_F(\mu) = \mu$, and it is called {\em odd under $T_F$} if $T_F(\mu) = -\mu$.

\begin{example}
Given $E, F \subset [n]$, the sets $A_E$ and $B_E$ are invariant under $T_F$.
\end{example}

\begin{remark}
\label{rema - prod even odd}
Given $E \subset [n]$ and $\mu, \nu \in \sM(\R^n)$, it follows from Proposition~\ref{prop - refl conv rn} that if $\mu$ or $\nu$ or both are even under $T_E$, then $\mu \mc \nu$ is even under $T_E$; similarly if $\mu$ or $\nu$ or both are odd under $T_E$, then $\mu \mc \nu$ is odd under $T_E$.
\end{remark}

\begin{proposition}
\label{prop - mult conv even odd support}
For $j \in \{1, 2\}$ let $\sE_j, \sF_e^j, \sF_o^j \subset \sP_n$ and $\mu_j \in \sM(\sE_j; \sF_e^j, \sF_o^j)$. Then $\mu_1 \mc \mu_2 \in \sM(\sE; \sF_e, \sF_o)$ where
\begin{eqnarray*}
\sE & = & \big\{ E_1 \cap E_2 \,:\, E_1 \in \sE_1,\; E_2 \in \sE_2 \big\}, \\
\sF_e & = & \sF_e^1 \cup \sF_e^2, \\
\sF_o & = & \sF_o^1 \cup \sF_o^2
\end{eqnarray*}
\end{proposition}

\begin{proof}
This follows from Corollary~\ref{coro - mult conv decomp note} and Remark~\ref{rema - prod even odd}.
\end{proof}

\begin{proposition}
\label{prop - mult conv even odd support sphere}
For $j \in \{1, 2\}$ let $\sE_j, \sF_e^j, \sF_o^j \subset \sP_n$ with $\O \notin \sE_j$, and let $\mu_j \in \sM^{\rm int}(\sE_j; \sF_e^j, \sF_o^j)$. Then $\mu_1 \mcs \mu_2 \in \sM^{\rm sph}(\sE; \sF_e, \sF_o)$ where
\begin{eqnarray*}
\sE & = & \big\{ E_1 \cap E_2 \,:\, E_1 \in \sE_1,\; E_2 \in \sE_2 \big\} \setminus \{ \O \}, \\
\sF_e & = & \sF_e^1 \cup \sF_e^2, \\
\sF_o & = & \sF_o^1 \cup \sF_o^2
\end{eqnarray*}
\end{proposition}

\begin{proof}
This follows from Corollary~\ref{coro - mult conv decomp note sphere} and Remark~\ref{rema - prod even odd}.
\end{proof}

Now let $\sF \subset \sP_n$. A set $A \in \sB(\R^n)$ is called {\em invariant under} $\sF$ if $T_F(A) = A$ holds for every $F \in \sF$. Furthermore $\mu \in \sM(\R^n)$ is called {\em even (odd) under} $\sF$ if $\mu$ is even (odd) under $T_F$ for every $F \in \sF$. For $\sF_e, \sF_o \subset \sP_n$ and $A \in \sB(\R^n)$ invariant under $\sF_e \cup \sF_o$ we denote by $\sM(A; \sF_e, \sF_o)$ the subset of those $\mu \in \sM(A)$ that are even under $\sF_e$ and odd under $\sF_o$. The integrable members of $\sM(A; \sF_e, \sF_o)$ are denoted by $\sM^{\rm int}(A; \sF_e, \sF_o)$.

Let $\sF_e, \sF_o \subset \sP_n$ and set $\sF = \sF_e \cup \sF_o$. Let $\sG$ be the group generated by $\sF$. Note that $\sF = \O$ is possible, yielding $\sG = \{ \O \}$. If $\sF \neq \O$, then for each $G \in \sG$ Proposition~\ref{prop - f and g generated} yields a (non-unique) representation $G = F_1 \Delta \ldots \Delta F_m$ where $m \geq 1$ and $F_1, \ldots, F_m \in \sF$. Denote by $N_o$ the number of sets from $\sF_o$ in this representation. Obviously $N_o$ can be zero.  The family of sets $G \in \sG$ such that there exists a representation with even (odd) $N_o$ is denoted by $\sG_e$ ($\sG_o$). In general it is possible that $\sG_e \cap \sG_o \neq \O$. In the case $\sF = \O$ we set $\sG_e = \{ \O \}$ and $\sG_o = \O$. We call $\sG_e$ and $\sG_o$ the {\em even part} and {\em odd part of} $\sG$, respectively. Clearly $\sF_e \subset \sG_e$, $\sF_o \subset \sG_o$, and $\sG_e \cup \sG_o = \sG$. Also note that $\sG_e$ is always a subgroup, though not necessarily the one generated by $\sF_e$. The procedure just mentioned defines a map
\[
\gamma : \sP_n^2 \longrightarrow \sP_n^2,\quad \gamma(\sF_e, \sF_o) = (\sG_e, \sG_o)
\]
An element $(\sF_e, \sF_o)$ of the domain $\sP_n^2$ of $\gamma$ is also called a {\em generating pair}. An element $(\sG_e, \sG_o)$ of its range is also called a {\em symmetry pair}. If $\sG_e \cap \sG_o = \O$, then $(\sF_e, \sF_o)$ and $(\sG_e, \sG_o)$ are called {\em proper}. In other words, for a proper symmetry pair and a fixed element $G \in \sG$ the number of terms from $\sF_o$ is either even in all representations in terms of sets from $\sF$ or it is odd in all such representations. In general the condition $\sF_e \cap \sF_o = \O$ is not sufficient for a generating pair to be proper.

\begin{remark}
\label{rema - repeated generation}
The procedure of first generating a group from a generating pair and then identifying the even and odd parts of the group can be performed repeatedly, which finally yields the same symmetry pair. In other words for any $\sF_e, \sF_o \subset \sP_n$ we always have
\[
(\gamma \circ \gamma)(\sF_e, \sF_o) \,=\, \gamma(\sF_e, \sF_o) 
\]
\end{remark}

\begin{remark}
\label{rema - characterisation sym pair}
Note that any generating pair $(\sF_e, \sF_o)$ (whether proper or not) is a symmetry pair if and only if all of the following conditions are satisfied:
\begin{thnumber}
\item $\sF_e$ is a subgroup of $\sP_n$.
\item $\forall E, F \in \sF_o \quad E \Delta F \in \sF_e$
\item $\forall E \in \sF_e, F \in \sF_o \quad E \Delta F \in \sF_o$
\end{thnumber}
\end{remark}

\begin{remark}
\label{rema - feneq0 multi refl}
Notice that in the context of multiple reflections we can always assume that $\O \in \sF_e$ as this corresponds to the identity transformation.
\end{remark}

\begin{remark}
\label{rema - measures even odd deg}
Let $\sF_e, \sF_o \subset \sP_n$, and let $\sG$ be the group generated by $\sF_e \cup \sF_o$. Further let $A \in \sB(\R^n)$ be invariant under $\sG$ and $\mu \in \sM(A; \sF_e, \sF_o)$. It follows that $\mu$ is even under $\sG_e$ and odd under $\sG_o$.  Consequently $\sM(A; \sF_e, \sF_o) = \sM(A; \sG_e, \sG_o)$. Notice that if $(\sF_e, \sF_o)$ is not proper, then there is $G \in \sG_e \cap \sG_o$. In this case $\mu$ is both even and odd under $T_G$ and therefore $\mu = 0$. So for $\sM(A; \sF_e, \sF_o)$ in order to be non-trivial we need to require that $(\sF_e, \sF_o)$ is proper.
\end{remark}

For the following analysis it is useful to introduce a specific notation for the restriction of a generating pair on $[n]$ to a subset $E$ of $[n]$. For any $E \subset [n]$ we denote
\[
\rho_E : \sP_n^2 \longrightarrow \sP(E)^2,\quad \rho_E(\sF_e, \sF_o ) = (\sF_e |_E, \sF_o |_E)
\]
In Remark~\eqref{rema - proj gen commute} we mentioned that the generation of a group commutes with the restriction to a subset $E$ of~$[n]$. This fact is extended by the following observation.

\begin{remark}
\label{rema - proj even odd}
For $E \subset [n]$ and $\sF_e, \sF_o \subset \sP_n$ we have
\[
\gamma \circ \rho_E \, (\sF_e, \sF_o ) = \rho_E \circ \gamma \, (\sF_e, \sF_o )
\]
So the generation of the group together with the subsequent extraction of the even and odd parts commutes with the restriction to~$E$.
\end{remark}

\begin{remark}
Note that it can occur that the generating pair $(\sF_e |_E, \sF_o |_E)$ is not proper even if $(\sF_e, \sF_o)$ is proper. For instance if there is $F \in \sF_o$ such that $E \cap F = \O$, then $\O \in \sF_o |_E$.
\end{remark}

\begin{example}
Let $n = 2$, $\sF_e = \{ \{1\} \}$, and $\sF_o = \{ \{ 2 \} \}$. Then $\sG_e = \{ \O, \{1\} \}$, and $\sG_o = \{ \{ 2 \}, \{ 1, 2 \} \}$. Hence the generating pair $(\sF_e, \sF_o)$ is proper. Setting $E = \{ \{ 1 \} \}$, the restrictions are $\sF_e |_E = \{ \{1\} \}$ and $\sF_o |_E = \{ \O \}$. Thus $(\sF_e |_E, \sF_o |_E)$ is not proper.
\end{example}

For $\sE, \sF_e, \sF_o \subset \sP_n$ and $A = \bigcup_{E \in \sE} A_E$ we also write $\sM(\sE; \sF_e, \sF_o)$ for $\sM(A; \sF_e, \sF_o)$ and $\sM^{\rm int}(\sE; \sF_e, \sF_o)$ for $\sM^{\rm int}(A; \sF_e, \sF_o)$. Further note that if $\mu \in \sM(\R^n)$ is concentrated on $S^{n-1}$, then also $T_E(\mu)$ is concentrated on $S^{n-1}$ for all $E \subset [n]$. If $\O \notin \sE$ and $B = \bigcup_{E \in \sE} B_E$, we write $\sM^{\rm sph}(\sE; \sF_e, \sF_o)$ for $\sM(B; \sF_e, \sF_o)$.

\subsection{Symmetrisation of measures}

For $\sF \subset \sP_n$ we define the following operators on $\sM(\R^n)$:
\begin{eqnarray}
\label{eq - def mf}
M_{\sF}^+ & = & 2^{-|\sF|} \, \prod_{F \in \sF} (T_{\O} + T_F), \\
M_{\sF}^- & = & 2^{-|\sF|}  \, \prod_{F\in \sF} (T_{\O} - T_F) \nonumber
\end{eqnarray}
For $\sF = \O$ both operators are defined to be the identity map. For any $\mu \in \sM(\R^n)$ it is evident that $M_{\sF}^+(\mu)$ is even under $\sF$ and $M_{\sF}^-(\mu)$ is odd under $\sF$. Moreover we have
\begin{equation}
\label{eq - mf projection property}
M_{\sF}^+ M_{\sF}^+ = M_{\sF}^+,\quad M_{\sF}^- M_{\sF}^- = M_{\sF}^-
\end{equation}
The normalisation factors in~\eqref{eq - def mf} are chosen such that \eqref{eq - mf projection property} hold. Note that $\mu$ is even under $\sF$ if and only if $M_{\sF}^+(\mu) = \mu$. Further $\mu$ is odd under $\sF$ if and only if $M_{\sF}^-(\mu) = \mu$. Finally, for $E \subset [n]$, if $\mu$ is even (odd) under~$T_E$, then both $M_{\sF}^+(\mu)$ and $M_{\sF}^-(\mu)$ are even (odd) under~$T_E$.

Furthermore for any group $\sG \subset \sP_n$ we define the following linear operator on $\sM(\R^n)$:
\begin{equation}
\label{eq - group average}
M_{\sG} = |\sG|^{-1} \sum_{G \subset \sG} T_G
\end{equation}
For any $\mu \in \sM(\R^n)$, $M_{\sG}(\mu)$ is even under $\sG$. Moreover $M_{\sG}(\mu) = \mu$ holds if and only if $\mu$ is even under $\sG$.

\begin{proposition}
\label{prop - equ sym group av}
Let $\sG \subset \sP_n$ be a group. Then $M_{\sG}^+ = M_{\sG}$ holds on $\sM(\R^n)$.
\end{proposition}

\begin{proof}
Let $\mu \in \sM(\R^n)$. Then
\[
M_{\sG}^+(\mu) \,=\, M_{\sG} M_{\sG}^+ (\mu) \,=\, M_{\sG}^+ M_{\sG} (\mu) \,=\, M_{\sG} (\mu)
\]
where the first and last equalities hold because $M_{\sG}^+(\mu)$ and $M_{\sG}(\mu)$ are even under $\sG$ respectively.
\end{proof}

In order to combine the two types of operators we define for $\sF_e, \sF_o \subset \sP_n$ the operator
\[
M_{\sF_e, \sF_o} = M_{\sF_e}^+ M_{\sF_o}^-
\]
Whenever $\sF_e$ and $\sF_o$ are evident from the context we may drop the subscripts and simply write $M$ instead of $M_{\sF_e, \sF_o}$. From \eqref{eq - mf projection property} it follows that
\[
M_{\sF_e, \sF_o} M_{\sF_e, \sF_o} = M_{\sF_e, \sF_o}
\]
We denote by $\sG$ the group generated by $\sF_e \cup \sF_o$ and by $\sG_e$ and $\sG_o$ the even and odd parts of $\sG$ respectively. For any $\mu \in \sM(\R^n)$ it is clear that $M_{\sF_e, \sF_o}(\mu)$ is even under $\sG_e$ and odd under $\sG_o$. Moreover $\mu$ is even under $\sF_e$ and odd under $\sF_o$ if and only if $M_{\sF_e, \sF_o}(\mu) = \mu$, and this holds if and only if $\mu$ is even under $\sG_e$ and odd under $\sG_o$. Similarly as in the proof of Proposition~\ref{prop - equ sym group av} we can show that
\[
M_{\sF_e, \sF_o} = M_{\sG_e, \sG_o}
\]
Note that the definitions imply that $M_{\sF_e, \sF_o}$ is a linear combination of $T_G$, $G \in \sG$. If $\mu$ is concentrated on $S^{n-1}$, then also $M_{\sF_e, \sF_o}(\mu)$ is concentrated on~$S^{n-1}$.

\begin{remark}
\label{rema - prop sym measures}
Letting $\sF_e, \sF_o \subset \sP_n$ some properties of $M = M_{\sF_e, \sF_o}$ are derived from Proposition~\ref{rema - mult refl commute} as follows. For any $E \subset [n]$, $M$ commutes with $P_E$ and $R_E$ on $\sM(\R^n)$, and it commutes with $P_S$ and $P^S_E$ on $\sM^{\rm int}(\R^n)$. It follows that if $\mu \in \sM(\R^n)$ is even under $\sF_e$ and odd under $\sF_o$, then $P_E(\mu)$ and $R_E(\mu)$ are even under $\sF_e$ and odd under $\sF_o$. Similarly if $\mu \in \sM^{\rm int}(\R^n)$ is even under $\sF_e$ and odd under $\sF_o$, then $P_S(\mu)$ and $P^S_E(\mu)$ are even under $\sF_e$ and odd under $\sF_o$. Moreover if $\mu$ is of order $E$, then $M(\mu)$ is of order $E$ or zero. If $\mu$ is integrable, then $M(\mu)$ is integrable.
\end{remark}

\begin{remark}
\label{rema - group conv}
Let $\sF_e, \sF_o \subset \sP_n$ and $\mu, \nu \in \sM(\R^n)$. Set $M = M_{\sF_e, \sF_o}$. It follows from Proposition~\ref{prop - refl conv rn} that
\[
M(\mu \mc \nu) \,=\, M(\mu) \mc \nu \,=\, \mu \mc M(\nu) \,=\, M(\mu) \mc M(\nu)
\]
where the last equality is a consequence of the fact that $M M = M$. If $\mu$ and $\nu$ are integrable, then also
\[
M(\mu \mcs \nu) \,=\, M(\mu) \mcs \nu \,=\, \mu \mcs M(\nu) \,=\, M(\mu) \mcs  M(\nu)
\]
\end{remark}

For $E \subset [n]$ we define $1_E \in \R^n$ by
\[
(1_E)_j = \left\{ \begin{array}{ll}
-1, & j \in E \\[1em]
1, & j \notin E
\end{array}\right.
\]
and write also $1_n = 1_{\O}$. Thus $1_n = (1, \ldots, 1)$. Multiple reflections applied to vectors yield $T_F(1_E) = 1_{E \Delta F}$ for $E, F \subset [n]$, in particular $T_F(1_n) = 1_F$, and applied to Dirac measures we have
\[
T_F(\delta_{1_E}) = \delta_{T_F(1_E)} = \delta_{1_{E \Delta F}},\quad T_F(\delta_{1_n}) = \delta_{1_F}
\]

\begin{remark}
Let $\mu \in \sM(\R^n)$, and let $\sG \subset \sP_n$ be a group. Then Remark~\ref{rema - group conv} implies
\[
M_{\sG}(\mu) \,=\, M_{\sG}(\mu \mc \delta_{1_n}) \,=\, \mu \mc M_{\sG}(\delta_{1_n})
\]
This shows that the action of $M_{\sG}$ on $\mu$ can be expressed as the multiplicative convolution of $\mu$ with a purely atomic probability measure that does not depend on $\mu$.
\end{remark}

The next proposition specifies a measure with a finite number of atoms that has precisely the reflection symmetries of a given proper symmetry pair.

\begin{proposition}
\label{prop - sym pair symmetrisation}
Let $\sF_e, \sF_o \subset \sP_n$. Set $\sF = \sF_e \cup \sF_o$ and $\rho =  M_{\sF_e, \sF_o}(\delta_{1_n})$. Denote by $\sG$ the group generated by $\sF$ and its even and odd parts by $\sG_e$ and $\sG_o$, respectively. Then:
\begin{thnumber}
\item \label{prop - sym pair symmetrisation 1} $\rho \in \sM^{\rm int}(A_{[n]}; \sG_e, \sG_o)$
\item \label{prop - sym pair symmetrisation 1a} There are unique $a_G \in \R$, $G \in \sG$, such that
\begin{equation}
\label{eq - rho decomp delta}
\rho \,=\, \sum_{G \in \sG} a_G \, \delta_{1_G}
\end{equation}
\item \label{prop - sym pair symmetrisation 2} $(\sF_e, \sF_o)$ is proper. $\;\Longleftrightarrow\;$ $M_{\sF_e, \sF_o} \neq 0$ $\;\Longleftrightarrow\;$ $\rho \neq 0$ $\;\Longleftrightarrow\;$ $a_{\O} > 0$
\item \label{prop - sym pair symmetrisation 3} If $E \subset [n]$ with $E \notin \sG$, then $\rho \perp T_E(\rho)$.
\item \label{prop - sym pair symmetrisation 4} Suppose $(\sF_e, \sF_o)$ is proper and let $E \subset [n]$. Then $\rho$ is even (odd) under $T_E$ if and only if $E \in \sG_e$ ($E \in \sG_o$).
\end{thnumber}
\end{proposition}

In the sum on the right-hand side of \eqref{eq - rho decomp delta} each term is an atom at a different location on the scaled unit sphere $n^{1/2} S^{n-1}$, one of which is $1_n$.

\begin{proof}[Proof of Proposition~\ref{prop - sym pair symmetrisation}]
Statements \eqref{prop - sym pair symmetrisation 1} and \eqref{prop - sym pair symmetrisation 1a} are clear.

We turn to statement~\eqref{prop - sym pair symmetrisation 2}. All implications from the right to the left are easy to see. It remains to show that the first condition implies that $a_{\O} > 0$. It follows from the definitions~\eqref{eq - def mf} that the coefficient $a_{\O}$ always has a strictly positive contribution from a power of $T_{\O}$, and it has further contributions if $\O = F_1 \Delta \ldots \Delta F_m$ holds for any $m \geq 1$ and $F_1, \ldots, F_m \in \sF$. Since $\O \notin \sG_o$, the number of terms from $\sF_o$ in such a representation must be even. Hence all possible other contributions to $a_{\O}$ are positive too. Consequently $a_{\O} > 0$.

To see statement~\eqref{prop - sym pair symmetrisation 3} let $E \subset [n]$ with $E \notin \sG$ and note that
\[
T_E(\rho) \,=\, \sum_{G \in \sG} a_G \, \delta_{1_{E \Delta G}}
\]
Since $E \Delta G \notin \sG$ for all $G \in \sG$, all the locations of the atoms in $T_E(\rho)$ are distinct from all the locations of the atoms in $\rho$.

In statement~\eqref{prop - sym pair symmetrisation 4} we show necessity since sufficiency is clear. Assume that $(\sF_e, \sF_o)$ is proper and that $\rho$ is either even or odd under $T_E$. It follows from \eqref{prop - sym pair symmetrisation 2} that $\rho \neq 0$ and from \eqref{prop - sym pair symmetrisation 3} that $E \in \sG$. Suppose $E \notin \sG_e$ ($E \notin \sG_o$). Then $E \in \sG_o$ ($E \in \sG_e$). Hence $\rho$ is odd (even) under $T_E$, and thus $\rho$ is not even (not odd) under $T_E$.
\end{proof}

\begin{remark}
Statements similar to that in Proposition~\ref{prop - sym pair symmetrisation} are true for $P_S(\rho)$, in particular $P_S(\rho) \in \sM(B_{[n]}; \sG_e, \sG_o)$. If $(\sF_e, \sF_o)$ is proper, then $P_S(\rho) \neq 0$ and if, in addition, $P_S(\rho)$ is even (odd) under $T_E$ for some $E \subset [n]$, then $E \in \sG_e$ ($E \in \sG_o$).
\end{remark}

\begin{remark}
\label{rema - symetr jacobi}
Let $\mu \in \sM(\R^n)$ and $\sF_e, \sF_o \subset \sP_n$, and set $M = M_{\sF_e, \sF_o}$. Then
\[
M(\mu) \,=\, M(\mu_+) - M(\mu_-) \,=\, (M(\mu))_+ - (M(\mu))_-,
\]
the last expression being the Jordan decomposition of $M(\mu)$ with both terms concentrated on disjoint sets, and the middle expression following from the linearity of $M$. Note however that $M(\mu_+)$ and $M(\mu_-)$ are not necessarily non-negative; so the minimality property of the Jordan decomposition does {\em not} yield
\[
(M(\mu))_+ \leq M(\mu_+),\quad (M(\mu))_- \leq M(\mu_-)
\]
in this case!
\end{remark}

\subsection{Specific symmetries}

A set $A \in \mathcal{B}(\R^n)$ is called {\em origin-symmetric} or {\em symmetric} if it is invariant under $T_{[n]}$, that is $A = -A$, and it is called {\em unconditional} if it is invariant under $\sP_n$. Moreover $\mu \in \sM(\R^n)$ where $n \geq 1$ is called {\em origin-symmetric} or {\em symmetric} if $\mu$ is even under $T_{[n]}$. Clearly this is the case if and only if $\mu(-A) = \mu(A)$ holds for every $A \in \mathcal{B}(\R^n)$. Further $\mu$ is called {\em anti-symmetric} if it is odd under $T_{[n]}$. This is equivalent to $\mu(-A) = -\mu(A)$ for every $A \in \mathcal{B}(\R^n)$. Finally $\mu$ is called {\em unconditional} if it is even under $\sP_n$.

It follows from the definitions that for any symmetric $A \in \mathcal{B}(\R^n)$ the subspace of symmetric members of $\sM(A)$ is given by $\sM(A; \{[n]\}, \O)$ and that of $\sM^{\rm int}(A)$ is given by $\sM^{\rm int}(A; \{[n]\}, \O)$. We also denote these two spaces by $\sM_{\rm sym}(A)$ and $\sM^{\rm int}_{\rm sym}(A)$. Similarly for any symmetric $A \in \mathcal{B}(\R^n)$ the subspace of anti-symmetric members of $\sM(A)$ is given by $\sM(A; \O, \{[n]\})$ and that of $\sM^{\rm int}(A)$ is given by $\sM^{\rm int}(A; \O, \{[n]\})$. These two spaces are also denoted by $\sM_{\rm asym}(A)$ and $\sM^{\rm int}_{\rm asym}(A)$. Moreover for any unconditional $A \in \mathcal{B}(\R^n)$ the subspace of unconditional members of $\sM(A)$ is given by $\sM(A; \sP_n, \O)$ and that of $\sM^{\rm int}(A)$ is given by $\sM^{\rm int}(A; \sP_n, \O)$. We also denote these two spaces by $\sM_{\rm unc}(A)$ and $\sM^{\rm int}_{\rm unc}(A)$.

For $\sE \subset \sP_n$ and $A = \bigcup_{E \in \sE} A_E$ we also write $\sM_{\rm sym}(\sE)$ for $\sM_{\rm sym}(A)$ and $\sM^{\rm int}_{\rm sym}(\sE)$ for $\sM^{\rm int}_{\rm sym}(A)$. If $\O \notin \sE$ and $B = \bigcup_{E \in \sE} B_E$, we write $\sM^{\rm sph}_{\rm sym}(\sE)$ for $\sM_{\rm sym}(B)$. Similar short notations are defined for unconditional and asymmetric measures in an obvious way.

\begin{proposition}
\label{prop - sym anti-sym}
Let $\mu, \nu \in \sM(\R^n)$. If $\mu$ is symmetric and $\nu$ is anti-symmetric, then $\mu \mc \nu = 0$.
\end{proposition}

\begin{proof}
This follows from Remark~\ref{rema - prod even odd}.
\end{proof}

In the two special cases $\sG = \{ \O, [n] \}$ and $\sG = \sP_n$ the operator~\eqref{eq - group average} is denoted by $M_{\rm sym}$ and $M_{\rm unc}$, respectively. Thus
\begin{equation}
\label{eq - msym measure def}
M_{\rm sym} = \frac{1}{2} (T_{\O} + T_{[n]}),\quad M_{\rm unc} = 2^{-n} \sum_{E \subset [n]} T_E
\end{equation}
More explicitly this yields, for $\mu \in \sM(\R^n)$ and $A \in \mathcal{B}(\R^n)$,
\[
M_{\rm sym}(\mu)(A) = \frac{1}{2}(\mu(A) + \mu(-A)),\quad M_{\rm unc}(\mu)(A) = 2^{-n} \sum_{s \in \sS^n} \mu(sA)
\]

If $\mu$ has a density with respect to the Lebesgue measure on $\R^n$, say $\mu = f \odot \lambda$, we set, for $x \in \R^n$,
\[
f_{\rm sym}(x) = \frac{1}{2}(f(x) + f(-x)),\quad f_{\rm unc}(x) = 2^{-n} \, \sum_{s \in \sS^n} f(sx)
\]
and obtain $M_{\rm sym}(\mu) = f_{\rm sym} \odot \lambda$ and $M_{\rm unc}(\mu) = f_{\rm unc} \odot \lambda$.

\begin{remark}
\label{rema - spec sym prop}
For $\mu \in \sM(\R^n)$ it follows from the general case that  $M_{\rm sym}(\mu)$ is symmetric and $M_{\rm unc}(\mu)$ is unconditional. Moreover $M_{\rm sym} M_{\rm sym} = M_{\rm sym}$ and $M_{\rm unc} M_{\rm unc} = M_{\rm unc}$. Further $\mu$ is symmetric if and only if $M_{\rm sym}(\mu) = \mu$, and $\mu$ is unconditional if and only if $M_{\rm unc}(\mu) = \mu$. Note also that
\[
\mu(\R^n) \,=\, M_{\rm sym}(\mu)(\R^n) \,=\, M_{\rm unc}(\mu)(\R^n)
\]
\end{remark}

For any unconditional $A \in \mathcal{B}(\R^n)$ the map $M_{\rm unc}$ yields a canonical bijection between $\sM(A \cap \R^n_+)$ and $\sM_{\rm unc}(A)$, which is very intuitive. This is formalised in Proposition~\ref{prop - bijection unc positive}. The following notation proves useful:
\begin{align*}
 & R_+ : \sM(\R^n) \longrightarrow \sM(\R^n_+),\\
 & R_+(\mu)(B) = \mu(B \cap \R^n_+),\quad B \in \mathcal{B}(\R^n)
\end{align*}
Note that $R_+ R_E = R_E R_+$ for every $E \subset [n]$.

\begin{proposition}
\label{prop - bijection unc positive}
Let $A \in \mathcal{B}(\R^n)$ be unconditional. Define
\[
\widetilde{M} : \sM(A \cap \R^n_+) \longrightarrow \sM_{\rm unc}(A),\quad \widetilde{M}(\mu) = M_{\rm unc}(\mu)
\]
and
\[
g : \R^n \longrightarrow \R^n_+,\quad g(x) = (|x_1|, \ldots, |x_n|)
\]
Then $\widetilde{M}$ is bijective and we have, for every $\mu \in \sM_{\rm unc}(A)$,
\begin{equation}
\widetilde{M}^{-1}(\mu) =  g(\mu) = \sum_{E \subset [n]} 2^{|E|} R_+(\mu_E) \label{eq - bijection unc positive}
\end{equation}
\end{proposition}

Observe the coefficients varying from term to term in the coordinate decomposition \eqref{eq - bijection unc positive}. Specifically the cases $A = \R^n$ and $A = A_{[n]}$ are interesting. The latter choice yields a bijection from $\sM((0,\infty)^n)$ to $\sM_{\rm unc}(A_{[n]})$. Note also that the restriction of $\widetilde{M}$ to integrable elements is a bijection between $\sM^{\rm int}(A \cap \R^n_+)$ and $\sM_{\rm unc}^{\rm int}(A)$ for each unconditional $A \in \mathcal{B}(\R^n)$.

\begin{proof}[Proof of Proposition~\ref{prop - bijection unc positive}]
First note that if $\mu \in \sM(\R^n)$ is concentrated on $A \cap \R^n_+$, then $M_{\rm unc}(\mu)$ is indeed concentrated on~$A$. Conversely if $\mu \in \sM_{\rm unc}(\R^n)$ is concentrated on~$A$, then the application of the operator on the right-hand side of \eqref{eq - bijection unc positive} yields a measure concentrated on~$\sM(A \cap \R^n_+)$. We may therefore assume that $A = \R^n$, the general case arising from restriction to the appropriate subspaces of measures.

Now let $E \subset [n]$, $\mu \in \sM_{\rm unc}(\R^n)$, and $B \in \mathcal{B}(\R^n)$. Then
\begin{eqnarray*}
g(\mu_E)(B) & = & \mu_E(g^{-1}(B)) \\
 & = & \mu_E(g^{-1}(B \cap \R^n_+)) \\
 & = & \mu_E \Big( \bigcup\nolimits_{F \subset [n]} T_F(B \cap \R^n_+) \Big) \\
 & = & \mu \Big( \bigcup\nolimits_{F \subset [n]} T_F(B \cap \R^n_+ \cap A_E) \Big) \\
 & = & \mu \Big( \bigcup\nolimits_{F \subset E} T_F(B \cap \R^n_+ \cap A_E) \Big) \\
 & = & \sum\nolimits_{F \subset E} \mu (T_F(B \cap \R^n_+ \cap A_E)) \\
 & = & \sum\nolimits_{F \subset E} \mu (B \cap \R^n_+ \cap A_E) \\
 & = & 2^{|E|} \mu_E(B \cap \R^n_+) \\
 & = & 2^{|E|} R_+(\mu_E)(B) 
\end{eqnarray*}
By linearity this shows the second equation in~\eqref{eq - bijection unc positive}.

Furthermore, for $\mu \in \sM(\R^n_+)$ and $B \in \mathcal{B}(\R^n_+)$, we have
\begin{eqnarray*}
(g \circ \widetilde{M})(\mu)(B) & = & M_{\rm unc}(\mu)(g^{-1}(B)) \\
 & = & 2^{-n} \sum\nolimits_{E \subset [n]} (T_E \mu)(g^{-1}(B)) \\
 & = & 2^{-n} \sum\nolimits_{E \subset [n]} \mu(T_E g^{-1}(B)) \\
 & = & 2^{-n} \sum\nolimits_{E \subset [n]} \mu(g^{-1}(B)) \\
 & = & \mu(g^{-1}(B)) \\
 & = & \mu(B)
\end{eqnarray*}

It remains to show that $(\widetilde{M} \circ g)(\mu) = \mu$ for $\mu \in \sM_{\rm unc}(\R^n)$. So let $\mu \in \sM_{\rm unc}(\R^n)$ and fix $F \subset [n]$. First note that for $j \notin F$ we have $T_{\{j\}} R_+ R_F = R_+ R_F$ and therefore
\[
(T_{\O} + T_{\{j\}}) R_+ (\mu_F) = 2 R_+ (\mu_F)
\]
Moreover, for $B \in \mathcal{B}(\R^n)$, we have
\begin{eqnarray*}
\prod\nolimits_{j \in F} (T_{\O} + T_{\{j\}}) R_+ (\mu_F) (B) & = & \sum\nolimits_{G \subset F} T_G R_+ (\mu_F) (B) \\
 & = & \sum\nolimits_{G \subset F} \mu_F ((T_G B) \cap \R^n_+) \\
 & = & \sum\nolimits_{G \subset F} \mu_F (T_G (B \cap (T_G \R^n_+)) \\
 & = & \sum\nolimits_{G \subset F} T_G(\mu_F) (B \cap (T_G \R^n_+)) \\
 & = & \sum\nolimits_{G \subset F} \mu_F (B \cap (T_G \R^n_+)) \\
 & = & \mu_F(B),
\end{eqnarray*}
the last equation being a consequence of the fact that the sets $ (T_G \R^n_+) \cap A_F$ for $G \subset F$ form a partition of $A_F$. Hence
\[
(\widetilde{M} \circ g)(\mu_F) \,=\, 2^{-n} \prod\nolimits_{j \in [n]} (T_{\O} + T_{\{j\}}) \, 2^{|F|} R_+(\mu_F) \,=\, \mu_F
\]
Summing over $F \subset [n]$ yields $(\widetilde{M} \circ g)(\mu) = \mu$.

\end{proof}

\begin{remark}
\label{rema - spec sym comm}
As special cases of Remark~\ref{rema - prop sym measures} note that the operators $M_{\rm sym}$ and $M_{\rm unc}$ commute with $P_E$ and $R_E$ on $\sM(\R^n)$ for every $n \geq 1$ and $E \subset [n]$. Moreover they commute with $P_S$ and $P^S_E$ on $\sM^{\rm int}(\R^n)$ for every $n \geq 1$ and $E \subset [n]$. It follows that if $\mu \in \sM(\R^n)$ is symmetric (unconditional), then $P_E(\mu)$ and $R_E(\mu)$ are symmetric (unconditional). Similarly if $\mu \in \sM^{\rm int}(\R^n)$ is symmetric (unconditional), then $P_S(\mu)$ and $P^S_E(\mu)$ are symmetric (unconditional). Moreover if $\mu$ is of order $E$, then $M(\mu)$ is of order $E$ or zero. Also note that if $\mu$ is integrable, then also $M_{\rm sym}(\mu)$ and $M_{\rm unc}(\mu)$ are integrable.
\end{remark}

\begin{remark}
\label{rema - sym conv prop}
Let $\mu, \nu \in \sM(\R^n)$. As special cases of Remark~\ref{rema - group conv} we obtain
\begin{thnumber}
\item $M_{\rm sym}(\mu \mc \nu) \,=\, M_{\rm sym}(\mu) \mc \nu \,=\, \mu \mc M_{\rm sym}(\nu) \,=\, M_{\rm sym}(\mu) \mc M_{\rm sym}(\nu)$
\item $M_{\rm unc}(\mu \mc \nu) \,=\, M_{\rm unc}(\mu) \mc \nu \,=\, \mu \mc M_{\rm unc}(\nu) \,=\, M_{\rm unc}(\mu) \mc M_{\rm unc}(\nu)$
\end{thnumber}
If $\mu$ and $\nu$ are integrable, then
\begin{thnumber}
\item $M_{\rm sym}(\mu \mcs \nu) \,=\, M_{\rm sym}(\mu) \mcs \nu \,=\, \mu \mcs M_{\rm sym}(\nu) \,=\, M_{\rm sym}(\mu) \mcs M_{\rm sym}(\nu)$
\item $M_{\rm unc}(\mu \mcs \nu) \,=\, M_{\rm unc}(\mu) \mcs \nu \,=\, \mu \mcs M_{\rm unc}(\nu) \,=\, M_{\rm unc}(\mu) \mcs M_{\rm unc}(\nu)$
\end{thnumber}
\end{remark}

In the case $\sG = \{ \O, [n] \}$, corresponding to symmetry with respect to the origin, we define
\[
\sigma^n_{\rm sym} \,=\, M_{\sG}(\delta_{1_n}) \,=\, 2^{-1} \big(\delta_{1_n} + \delta_{-1_n}\big)
\]
In the case $\sG = \sP_n$, corresponding to unconditionality, we define
\[
\sigma^n_{\rm unc} \,=\, M_{\sG}(\delta_{1_n}) \,=\, 2^{-n} \, \sum_{s \in \sS^n} \delta_s
\]

\begin{remark}
\label{rema - sigma-sym sigma-unc}
The measures $\sigma^n_{\rm sym}$ and $\sigma^n_{\rm unc}$ are purely atomic probability measures for every $n \geq 1$ as in the case of general $\sG$, and they are concentrated in $n^{1/2} S^{n-1}$, and so they are clearly integrable. Moreover $\sigma^n_{\rm sym}$ is symmetric, and $\sigma^n_{\rm unc}$ is unconditional. For any $\mu \in \sM(\R^n)$ we have
\begin{equation}
\label{eq - m sym unc conv}
M_{\rm sym}(\mu) = \mu \mc \sigma^n_{\rm sym},\quad M_{\rm unc}(\mu) = \mu \mc \sigma^n_{\rm unc}
\end{equation}
\end{remark}

\begin{remark}
Let $\xi$ be a random vector in $\R^n$. Denote by $\zeta_{\rm sym}$ and $\zeta_{\rm unc}$ two independent random vectors that are defined on the same probability space and distributed as $\sigma^n_{\rm sym}$ and $\sigma^n_{\rm unc}$, respectively. Then we may define random vectors
\begin{equation*}
\label{eq - constr sym unc rv}
\xi_{\rm sym} = \xi \zeta_{\rm sym},\quad \xi_{\rm unc} = \xi \zeta_{\rm unc}
\end{equation*}
Setting $\mu = \mathcal{L}(\xi)$, then \eqref{eq - m sym unc conv} yields $M_{\rm sym}(\mu) = \mathcal{L}(\xi_{\rm sym})$ and $M_{\rm unc}(\mu) = \mathcal{L}(\xi_{\rm unc})$. Hence the distribution of $\xi_{\rm sym}$ is symmetric and that of $\xi_{\rm unc}$ is unconditional.
\end{remark}

\begin{remark}
Recall that if $\mu \in \sM(S^{n-1})$, then $M_{\rm sym}(\mu), M_{\rm unc}(\mu) \in \sM(S^{n-1})$. In this case the action of both operators can also be expressed as the convolution with certain measures on the sphere as follows:
\[
M_{\rm sym}(\mu) = \mu \mcs P_S(\sigma^n_{\rm sym}),\quad M_{\rm unc}(\mu) = \mu \mcs P_S(\sigma^n_{\rm unc})
\]
These two relations do not hold for general $\mu$ of course.
\end{remark}

\subsection{Symmetry decomposition}
\label{subsec - symmetry decomposition}

For $J \subset [n]$ with $J \neq \O$ we define
\[
\sigma_J(x) \,=\, \prod_{j \in J} \sign(x_j),\quad x \in \R^n,
\]
where $\sign(a) \in \{ -1, 0, 1 \}$ is the usual sign of $a \in \R$, and $\sigma_{\O}(x) = 1$. Further we set $\sS^n = \{ -1, 1\}^n$. For $E \subset [n]$ we denote by $\sS^E$ the set of all $s \in \{ -1, 0, 1 \}^n$ such that $s_j \in \{ -1, 1\}$ for $j \in E$ and $s_j = 0$ for $j \notin E$. For $J \subset E \subset [n]$ with $J \neq \O$ we therefore have
\[
\sigma_J(s) = \prod_{j \in J} s_j,\quad s \in \sS^E
\]
For $J \subset E \subset [n]$ with $E \neq \O$ we define the purely atomic measure
\[
\delta_{E, J} \,=\, 2^{-|E|} \sum_{s \in \sS^E} \sigma_J(s) \, \delta_s,
\]
Moreover we set $\delta_{\O, \O} = \delta_0$. In the special case $E = [n]$ we also use the short notation $\delta_J = \delta_{[n], J}$.

\begin{remark}
\label{rema - sigma J sigma K sym difference}
For $J, K \subset [n]$ we have
\[
\sigma_J(s) \, \sigma_K(s) = \sigma_{J \Delta K}(s),\quad s \in \sS^n
\]
\end{remark}

\begin{remark}
\label{rema - trafo deltaj te}
For $E, F \subset [n]$, and $J \subset E$ we have
\[
T_F(\delta_{E, J}) = (-1)^{|J \cap F|} \, \delta_{E, J}
\]
In particular note that $\delta_{E, J}$ is either even or odd under $T_F$.
\end{remark}

\begin{remark}
\label{re - e j sym}
Notice that $\delta_{E, J} \in \sM^{\rm int}(A_E)$ and $\delta_{E, J} \neq 0$ hold for every $J \subset E \subset [n]$. Further $\delta_{E, J}$ is origin-symmetric if and only if $|J|$ is even, and it is anti-symmetric if and only if $|J|$ is odd. Moreover $\delta_{E, J}$ is unconditional if and only if $J = \O$. Furthermore
\[
\delta_{E, J} = \sigma_J \odot \delta_{E, \O},\quad\quad \delta_{J} = \sigma_J \odot \delta_{\O},\quad\quad \delta_{\O} = \sigma^n_{\rm unc},
\]
and more generally
\[
\sigma_J \odot \delta_K = \sigma_{J \Delta K} \odot \delta_{\O},\quad J, K \subset [n]
\]
Now let $\mu \in \sM(\R^n)$. If $\mu$ is symmetric and $|J|$ is odd, then $\delta_J \mc \mu = 0$ by Proposition~\ref{prop - sym anti-sym}. The same identity holds if $\mu$ is anti-symmetric and $|J|$ is even, or if $\mu$ is unconditional and $J \neq \O$ by Remark~\ref{rema - prod even odd}.
\end{remark}

\begin{remark}
\label{rema - pf deltaej}
For $J \subset F \subset E \subset [n]$ it is easy to see that
\[
P_F(\delta_{E, J}) \;=\; \delta_{F, J}
\]
On the other hand, for $J \subset E \subset [n]$ and $F \subset E$ such that $J \not\subset F$ we have $P_F(\delta_{E, J}) = 0$. Consequently, for $\mu \in \sM(\R^n)$, if $\mu$ is concentrated on $H_E$, then
\[
\mu \mc \delta_{[n],J} \;=\; \mu \mc \delta_{E,J}
\]
\end{remark}

\begin{proposition}
\label{prop - sum sigma J}
Let $J \subset [n]$. Then 
\[
\sum_{s \in \sS^n} \sigma_J(s) \,=\, 2^n \, \one \{ J = \O \}
\]
\end{proposition}

\begin{remark}
\label{rema - sum sigma j sigma k}
For $J, K \subset [n]$ Remark~\ref{rema - sigma J sigma K sym difference} and Proposition~\ref{prop - sum sigma J} imply 
\[
\sum_{s \in \sS^n} \sigma_J(s) \sigma_K(s) \,=\, 2^n \, \one \{ J = K \}
\]
\end{remark}

\begin{proposition}
\label{prop - sigma s t identities}
Let $n \geq 1$ and $s, t \in \sS^n$. Then
\begin{thnumber}
\item \label{prop - sigma s t identities 1} $\sum_{J \subset [n]} \sigma_J(s) \,=\, 2^n \, \one \{s = 1_n \}$
\item \label{prop - sigma s t identities 2} $\sum_{J \subset [n]} \sigma_J(s) \sigma_J(t) \,=\, 2^n \, \one \{s = t \}$
\item \label{prop - sigma s t identities 3} $\sum_{J \subset [n]} \delta_J \,=\, \delta_{1_n}$
\end{thnumber}
\end{proposition}

\begin{proof}
We first show~\eqref{prop - sigma s t identities 1}. The assertion is clear for $s = 1_n$. So assume that $s \neq 1_n$. Without loss of generality $s_1 = -1$. Then
\[
\sum_{J \subset [n]} \sigma_J(s) \,=\, \sum_{K \subset \{1\}} \sigma_K((s_1)) \sum_{L \subset \{2, 3, \ldots, n\}} \sigma_L((s_2, \ldots, s_n)) \,=\, 0
\]

Statement~\eqref{prop - sigma s t identities 2} follows from \eqref{prop - sigma s t identities 1} and the fact that
$\sigma_J(s) \sigma_J(t) = \sigma_J(st)$ for $J \subset [n]$.

Statement~\eqref{prop - sigma s t identities 3} is a consequence of~\eqref{prop - sigma s t identities 1}.
\end{proof}

\begin{proposition}
\label{prop - sign meas decomp signs rn}
Let $\nu \in \sM(\R^n)$. Then
\begin{thnumber}
\item \label{prop - sign meas decomp signs rn decomp} $\nu = \sum_{K \subset [n]} (\delta_K \mc \nu)$
\item \label{prop - sign meas decomp signs rn zero} $\nu = 0$ if and only if $\, \delta_J \mc \nu = 0$ for every $J \subset [n]$.
\end{thnumber}
\end{proposition}

\begin{proof}
Statement~\eqref{prop - sign meas decomp signs rn decomp} follows from Proposition~\ref{prop - sigma s t identities}~\eqref{prop - sigma s t identities 3}.

In statement~\eqref{prop - sign meas decomp signs rn zero} sufficiency is a consequence of statement~\eqref{prop - sign meas decomp signs rn decomp} and necessity is trivial.
\end{proof}

\begin{remark}
As an important (and obvious) consequence of Proposition~\ref{prop - sign meas decomp signs rn}~\eqref{prop - sign meas decomp signs rn zero} note that for $\mu, \nu \in \sM(\R^n)$ we have $\mu = \nu$ if and only if $\delta_J \mc \mu = \delta_J \mc \nu$ for every $J \subset [n]$. This also implies that in the expansion of Proposition~\ref{prop - sign meas decomp signs rn}~\eqref{prop - sign meas decomp signs rn decomp} each term is unique for given~$\nu$.
\end{remark}

\begin{corollary}
\label{coro - sign meas bn decomp signs}
Let $\nu \in \sM(S^{n-1})$. Then
\begin{thnumber}
\item \label{coro - sign meas bn decomp signs decomp} $\nu = \sum_{K \subset [n]} \delta_K \mcs \nu$
\item \label{coro - sign meas bn decomp signs zero} $\nu = 0$ if and only if $\delta_J \mcs \nu = 0$ for every $J \subset [n]$.
\end{thnumber}
\end{corollary}

\begin{proof}
This follows from Proposition~\ref{prop - sign meas decomp signs rn} and the fact that $\delta_J \mcs \nu = \delta_J \mc \nu$ for every $J \subset [n]$.
\end{proof}

\begin{proposition}
\label{prop - prod delta j and k}
Let $J, K \subset [n]$. Then
\[
\delta_J \mc \delta_K \,=\, \delta_J \, \one \{ J = K \}
\]
\end{proposition}

\begin{proof}
\begin{eqnarray}
\delta_J \mc \delta_K & = & 2^{-2n} \sum_{s, t \in \sS^n} \sigma_J(s) \, \sigma_K(t) \, \delta_{st} \nonumber\\
 & = & 2^{-2n} \sum_{s, r \in \sS^n} \sigma_J(s) \, \sigma_K(sr) \, \delta_{r} \nonumber\\
 & = & 2^{-2n} \sum_{s, r \in \sS^n} \sigma_J(s) \, \sigma_K(s) \, \sigma_K(r) \, \delta_{r} \nonumber\\
 & = & \delta_J \, \one \{ J = K \} \label{eq - sigma difference}
\end{eqnarray}
where \eqref{eq - sigma difference} follows from Remark~\ref{rema - sum sigma j sigma k}.
\end{proof}

\begin{lemma}
\label{le - density conv}
Let $\mu, \nu \in \sM(\R^n)$ and $J \subset [n]$. Then
\begin{thnumber}
\item \label{le - density conv 1} $\sigma_J \odot (\mu \mc \nu) \,=\, (\sigma_J \odot \mu) \mc (\sigma_J \odot \nu)$
\item \label{le - density conv 2} $\delta_{\O} \mc (\sigma_J \odot \mu) \,=\, \sigma_J \odot (\delta_J \mc \mu)$
\end{thnumber}
\end{lemma}

\begin{proof}
The proof of statement~\eqref{le - density conv 1} is straightforward.

Statement~\eqref{le - density conv 2} follows from statement~\eqref{le - density conv 1} since
\[
\sigma_J \odot (\delta_J \mc \mu) \,=\, (\sigma_J \odot \delta_J) \mc (\sigma_J \odot \mu) \,=\, \delta_{\O} \mc (\sigma_J \odot \mu)
\]
\end{proof}

\begin{remark}
Lemma~\ref{le - density conv}~\eqref{le - density conv 2} shows that the right-hand side of this identity is unconditional as the left-hand side has this property. Moreover this statement is trivial in the case $J = \O$.
\end{remark}

\begin{corollary}
\label{le - density conv sphere}
Let $\mu, \nu \in \sM^{\rm int}(\R^n)$ and $J \subset [n]$. Then
\[
\sigma_J \odot (\mu \mcs \nu) \,=\, (\sigma_J \odot \mu) \mcs (\sigma_J \odot \nu)
\]
\end{corollary}

\begin{proof}
This follows from Corollary~\ref{coro - multconv dens poshom} and
Lemma~\ref{le - density conv}~\eqref{le - density conv 1}.
\end{proof}

\begin{remark}
\label{rema - refl dens}
Let $\mu \in \sM(\R^n)$ and $E, J \subset [n]$. Then
\[
T_E (\sigma_J \odot \mu) \,=\, (-1)^{|E \cap J|} \, \sigma_J \odot T_E(\mu)
\]
\end{remark}

\subsection{Index transformation}

For $E \subset [n]$ and $\sF_e, \sF_o \subset \sP_n$ we define $\sJ(E; \sF_e, \sF_o)$ to be the family of $J \subset E$ such that $|J \cap F|$ is even for every $F \in \sF_e$ and odd for every $F \in \sF_o$. We use the short notation
\[
\sJ(\sF_e, \sF_o) = \sJ([n]; \sF_e, \sF_o)
\]
Moreover we agree that whenever $\sF_e$ or $\sF_o$ is a singleton, then we may omit the corresponding curly brackets; for example we write $\sJ(F, \sF_o)$ instead of $\sJ(\{F\}, \sF_o)$. For definiteness we agree that $\O$ in the second argument means $\sF_o = \O$ and not $\sF_o = \{\O\}$, the latter forming inevitably a non-proper generating pair. With respect to the first argument such a distinction is immaterial.

\begin{remark}
For $E, F \subset [n]$ we have
\begin{eqnarray*}
\sJ(E; \O, \O) & = & \sP(E),\\
\sJ(E; [n], \O) \;=\; \sJ(E; E, \O) & = & \big\{ J \subset E \,:\, |J| \;\, {\rm is \; even} \,\big\},\\
\sJ(E; \sP_n, \O) \;=\; \sJ(E; \sP(E), \O) & = & \{ \O \},\\
\sJ(E; \O, [n]) \;=\; \sJ(E; \O, E) & = & \big\{ J \subset E \,:\, |J| \;\, {\rm is \; odd} \,\big\},\\
\sJ(E; F, \O) & = & \sP(E) \setminus \sJ(E; \O, F)
\end{eqnarray*}
\end{remark}

\begin{remark}
\label{rema - j e properties}
Let $E \subset [n]$ and $\sF_e, \sF_o \subset \sP_n$. Then
\begin{thnumber}
\item $\O \in \sJ(E; \sF_e, \sF_o)$ if and only if $\sF_o = \O$.
\item For every $J \in \sJ(E; \sF_e, \sF_o)$ it follows from Remark~\ref{rema - trafo deltaj te} that
\[
\delta_{E, J} \in \sM^{\rm int}(A_E; \sF_e; \sF_o)
\]
and from Remark~\ref{rema - pf deltaej} that
\[
P_G(\delta_{E, J}) \in \sM^{\rm int}(A_{G \cap E}; \sF_e, \sF_o),\quad G \subset [n]
\]
\item For $E \subset F \subset [n]$ we have $\sJ(E; \sF_e, \sF_o) \subset \sJ(F; \sF_e, \sF_o)$.
\end{thnumber}
\end{remark}

\begin{proposition}
\label{prop - cardinality even odd}
Let $A_1, \ldots, A_m$ be finite sets where $m \geq 1$. Set $A = A_1 \Delta \ldots \Delta A_m$ and
\[
N \,=\, \text{\#} \{ \, i \in [m] : |A_i| \,\; {\rm is \; odd} \, \}
\]
Then $|A|$ is even (odd) if $N$ is even (odd).
\end{proposition}

\begin{proof}
This is straightforward for $m = 2$, and can be shown iteratively for general $m$.
\end{proof}

\begin{corollary}
\label{coro - prod even odd}
Let $m, n \geq 1$ and $J, F_1, \ldots, F_m \subset [n]$. Set $F = F_1 \Delta \ldots \Delta F_m$ and
\[
N \,=\, \text{\#} \{ \, i \in [m] : |J \cap F_i| \,\; {\rm is \; odd} \, \}
\]
Then $|J \cap F|$ is even (odd) if $N$ is even (odd).
\end{corollary}

\begin{proof}
This is a consequence of the distribution law \eqref{eq - distr law sets} and Proposition~\ref{prop - cardinality even odd}.
\end{proof}

\begin{theorem}
\label{th - deltaj zero arb e sym}
Let $\nu \in \sM(\R^n)$ and $F \subset [n]$. Then:
\begin{thnumber}
\item $T_F(\nu) = \nu$ if and only if $\, \delta_J \mc \nu = 0$ for every $J \in \sJ(\O, F)$.
\item $T_F(\nu) = -\nu$ if and only if $\, \delta_J \mc \nu = 0$ for every $J \in \sJ(F, \O)$.
\end{thnumber}
\end{theorem}

\begin{proof}
Note that
\[
T_F(\delta_J \mc \nu) \,=\, T_F(\delta_J) \mc \nu \,=\, (-1)^{|J \cap F|} \, \delta_J \mc \nu
\]
where the first equality is a consequence of Proposition~\ref{prop - refl conv rn} and the second one follows from Remark~\ref{rema - trafo deltaj te}. Now the two claims follow from Proposition~\ref{prop - sign meas decomp signs rn}.
\end{proof}

\begin{corollary}
\label{coro - deltaj zero even odd}
Let $\nu \in \sM(\R^n)$ and $\sF_e, \sF_o \subset \sP_n$. Then $\nu \in \sM(\R^n; \sF_e, \sF_o)$ if and only if $\, \delta_J \mc \nu = 0$ for every $J \subset [n]$ such that $J \notin \sJ(\sF_e, \sF_o)$ (equivalently: such that $|J \cap F|$ is odd for some $F \in \sF_e$ or even for some $F \in \sF_o$).
\end{corollary}

\begin{proof}
This is a consequence of Theorem~\ref{th - deltaj zero arb e sym}.
\end{proof}

\begin{corollary}
\label{coro - deltaj zero even odd sphere}
Let $\nu \in \sM(S^{n-1})$ and $\sF_e, \sF_o \subset \sP_n$. Then $\nu \in \sM(S^{n-1}; \sF_e, \sF_o)$ if and only if $\, \delta_J \mcs \nu = 0$ for every $J \subset [n]$ such that $J \notin \sJ(\sF_e, \sF_o)$ (equivalently: such that $|J \cap F|$ is odd for some $F \in \sF_e$ or even for some $F \in \sF_o$).
\end{corollary}

\begin{proof}
This follows from Corollary~\ref{coro - deltaj zero even odd} and the fact that $\delta_J \mcs \nu = \delta_J \mc \nu$ for every $J \subset [n]$ since $\nu \in \sM(S^{n-1})$.
\end{proof}

Of the following proposition, statement~\eqref{prop - j gen group even odd 1} is used in the proofs of the supplements in Theorems \ref{th - measures univ general even odd ae} and \ref{th - measures univ general even odd be}, and statement~\eqref{prop - j gen group even odd 2} is referred to in Remark~\ref{rema - non proper conditions}.

\begin{proposition}
\label{prop - j gen group even odd}
Let $\sF_e,\sF_o \subset \sP_n$ and $\sF = \sF_e \cup \sF_o$. Further let $\sG$ be the group generated by $\sF$ and let $\sG_e$ and $\sG_o$ be its even and odd parts, respectively. Then:

\begin{thnumber}
\item \label{prop - j gen group even odd 1} $\sJ(\sF_e, \sF_o) = \sJ(\sG_e, \sG_o)$
\item \label{prop - j gen group even odd 2} $(\sF_e, \sF_o)$ is not proper if and only if $\sJ(\sF_e, \sF_o) = \O$.
\end{thnumber}
\end{proposition}

\begin{proof}
We first prove \eqref{prop - j gen group even odd 1}. The case $\sF = \O$ is clear since both sides are equal to $\sP_n$. Suppose $\sF \neq \O$. Let $J \in \sJ(\sF_e, \sF_o)$ and $G \in \sG$. By Proposition~\ref{prop - f and g generated} there are $m \geq 1$ and $F_1, \ldots, F_m \in \sF$ such that $G = F_1 \Delta \ldots \Delta F_m$. Denote by $N_o$ the number of sets from $\sF_o$ in this representation. Defining $N$ as in Corollary~\ref{coro - prod even odd} we have $N_o = N$. If $G$ is in $\sG_e$ ($\sG_o$), then $N_o$ is even (odd), and thus $|G \cap J|$ is even (odd) by Corollary~\ref{coro - prod even odd}. Therefore $J \in \sJ(\sG_e, \sG_o)$. Consequently $\sJ(\sF_e, \sF_o) \subset \sJ(\sG_e, \sG_o)$. The converse inclusion is clear.

We next show \eqref{prop - j gen group even odd 2}. Suppose $(\sF_e, \sF_o)$ is not proper and let $G \in \sG_e \cap \sG_o$. Assume there exists $J \in \sJ(\sF_e, \sF_o)$. It follows from \eqref{prop - j gen group even odd 1} that $J \in \sJ(\sG_e, \sG_o)$. Hence $|G \cap J|$ must be both even and odd, a contradiction. Conversely assume that $(\sF_e, \sF_o)$ is proper. Then there exists $\mu \in \sM(\R^n; \sG_e, \sG_o)$ different from zero by Proposition~\ref{prop - sym pair symmetrisation}. By Corollary~\ref{coro - deltaj zero even odd} we have $\delta_J \mc \mu = 0$ for every $J \subset [n]$ with $J \notin \sJ(\sG_e, \sG_o)$. Since $\mu \neq 0$, it follows from Proposition~\ref{prop - sign meas decomp signs rn} that there is $J \in \sJ(\sG_e, \sG_o)$ such that $\delta_J \mc \mu \neq 0$.
\end{proof}

The following lemma is required for Lemmas \ref{le - univ gen even odd gen supp} and \ref{le - univ gen sym gen supp even odd sphere}.

\begin{lemma}
\label{le - deltaj even odd class}
Let $\sF_e, \sF_o \subset \sP_n$ and $J \in \sJ(\sF_e, \sF_o)$. Then
\begin{thnumber}
\item \label{le - deltaj even odd class 1} $0 \,\neq\, \delta_J \mc \sigma_0^n \,\in\, \sM^{\rm int}(A_{[n]}; \sF_e, \sF_o)$
\item \label{le - deltaj even odd class 2} $0 \,\neq\, \delta_J \mcs \sigma_0^n \,\in\, \sM(B_{[n]}; \sF_e, \sF_o)$
\end{thnumber}
\end{lemma}

\begin{proof}
Statement~\eqref{le - deltaj even odd class 1} follows from Remark~\ref{rema - j e properties} and Proposition~\ref{prop - mult conv even odd support}.

Statement~\eqref{le - deltaj even odd class 2} is a consequence of \eqref{le - deltaj even odd class 1} and Propositions \ref{prop - proj S decomp E} and~\ref{rema - mult refl commute}.
\end{proof}

The following lemma is needed for the supplements in Theorems \ref{th - measures univ general even odd ae} and \ref{th - measures univ general even odd be}.

\begin{lemma}
\label{le - j s g impl}
Let $E \subset [n]$, and let $(\sG_e,\sG_o)$ be a proper symmetry pair. Then:
\begin{thnumber}
\item \label{le - j s g impl even} $E \in \sG_e$ if and only if $\sJ(\sG_e, \sG_o) \subset \sJ(E, \O)$.
\item \label{le - j s g impl odd} $E \in \sG_o$ if and only if $\sJ(\sG_e, \sG_o) \subset \sJ(\O, E)$.
\end{thnumber}
\end{lemma}

\begin{proof}
Necessity in \eqref{le - j s g impl even} and \eqref{le - j s g impl odd} follows directly from the definition of~$\sJ$.

To show sufficiency let $\rho =  M_{\sG_e, \sG_o}(\delta_{1_n})$. Then $\rho \in \sM(\R^n; \sG_e, \sG_o)$ by Proposition~\ref{prop - sym pair symmetrisation}~\eqref{prop - sym pair symmetrisation 1}. Therefore $\delta_J \mc \rho = 0$ holds for $J \in [n] \setminus \sJ(\sG_e, \sG_o)$ by Corollary~\ref{coro - deltaj zero even odd}. To see sufficiency in \eqref{le - j s g impl even} assume that $\sJ(\sG_e, \sG_o) \subset \sJ(E, \O)$. It follows that $\rho$ is even under $T_E$ again by Corollary~\ref{coro - deltaj zero even odd}. Now Proposition~\ref{prop - sym pair symmetrisation}~\eqref{prop - sym pair symmetrisation 4} implies that $E \in \sG_e$. Sufficiency in \eqref{le - j s g impl odd} is shown in a similar way.
\end{proof}

As a consequence of Lemma~\ref{le - j s g impl} the following result holds, which is applied in Example~\ref{ex - deltaej universal} and also interesting by itself.

\begin{theorem}
\label{th - j trafo group}
Let $\sG \subset \sP_n$ be a group. Then $\sJ(\sG, \O)$ is a group and $\sJ(\sJ(\sG, \O), \O) = \sG$.
\end{theorem}

\begin{proof}
To see the first claim let $A, B \in \sJ(\sG, \O)$ and $G \in \sG$. Since $|A \cap G|$ and $|B \cap G|$ are even, also $|(A \Delta B) \cap G|$ is even.

To see the second claim note that clearly $\sG \subset \sJ(\sJ(\sG, \O), \O)$. For the converse inclusion let $E \subset [n]$ such that $E \notin \sG$. Then there is $J \in \sJ(\sG, \O)$ such that $J \notin \sJ(E, \O)$ by Lemma~\ref{le - j s g impl}. Hence $|E \cap J|$ is odd. Consequently $E \notin \sJ(\sJ(\sG, \O), \O)$.
\end{proof}

\section{Lifting}
\label{sec - lifting}

We introduce a transformation of measures that is related to the so-called lifting, see~\cite[Section~2.2]{mos02}. Given $\mu \in \sM^{\rm int}(\R^n)$ where $n \geq 1$ we define a measure on $\R^{n+1}$
\[
L(\mu) = P_S M_{\rm sym} (\delta_1 \otimes \mu)
\]
where $\delta_1$ is a measure on $\R$ that consists of an atom at $1$ with mass~$1$. $L$ is a composition of maps. The first step, $\mu \mapsto \delta_1 \otimes \mu$, is the proper "lifting", which is linear. The first factor is considered as $0$-th coordinate in $\R^{n+1}$. The integrability of $\delta_1 \otimes \mu$ follows from Proposition~\ref{prop - integr prod}. Then the operators $M_{\rm sym}$ and $P_S$ are applied, in arbitrary order due to Remark~\ref{rema - prop sym measures}. Thus $L$ is linear and $L(\mu)$ is always origin-symmetric. $L(\mu)$ is always concentrated on
\[
S^n_0 \,\defeq\, S^n \setminus H_{[n]} \,=\, \{ u \in S^n \,:\, u_0 \neq 0 \}
\]
where $S^n$ denotes the unit sphere in~$\R^{n+1}$ and $[n] = \{1, \ldots, n\}$ as before. Further we set
\[
S^n_+ \,=\, \{ u \in S^n \,:\, u_0 > 0 \}
\]

\begin{lemma}
\label{lemma - finite meas bij}
The map $\mu \mapsto P_S(\delta_1 \otimes \mu)$ is a bijection from $\sM^{\rm int}(\R^n)$ to $\sM(S^n_+)$.
\end{lemma}

\begin{proof}
First note that
\[
1 \,\leq\, \|(1,x)\|_2 \,\leq\, 1 + \|x\|_2,\quad x \in \R^n
\]
Therefore we may define the map
\[
h_1 : \sM^{\rm int}(\R^n) \rightarrow \sM(\R^n),\quad h_1(\mu)(A) = \int_A \mu(dx) \, \|(1,x)\|_2,\quad A \in \mathcal{B}(\R^n)
\]
Moreover we have
\[
0 \,<\, \|(1,x)\|^{-1}_2 \,\leq\, 1,\quad x \in \R^n
\]
Hence the inverse map is well-defined as
\[
h_1^{-1} : \sM(\R^n) \rightarrow \sM^{\rm int}(\R^n),\quad h_1^{-1}(\rho)(A) = \int_A \rho(dx) \, \|(1,x)\|^{-1}_2,\quad A \in \mathcal{B}(\R^n)
\]
Now we define another map by
\[
h_2 : \R^n \rightarrow S^n_+,\quad h_2(x) = \frac{(1, x)}{\|(1, x)\|_2}
\]
which is clearly bijective. This induces a bijective map between measures,
\[
h_2 : \sM(\R^n) \rightarrow \sM(S^n_+)
\]
It remains to note that
\[
h_2 \circ h_1 (\mu) = P_S(\delta_1 \otimes \mu),\quad \mu \in \sM^{\rm int}(\R^n)
\]
\end{proof}

\begin{proposition}
\label{prop - bij meas sphere rn}
$L$ is a bijection from $\sM^{\rm int}(\R^n)$ to $\sM_{\rm sym}(S^n_0)$.
\end{proposition}

\begin{proof}
Use Lemma~\ref{lemma - finite meas bij} and note that $M_{\rm sym}$ is a bijection from $\sM(S^n_+)$ to $\sM_{\rm sym}(S^n_0)$.
\end{proof}

\begin{remark}
\label{rema - L commutes refl}
Let $E \subset [n]$ and $\mu \in \sM^{\rm int}(\R^n)$. Then
\[
T_E \, L (\mu) \,=\, L \, T_E (\mu)
\]
Consequently $\mu$ is even (odd) under $T_E$ if and only if $L(\mu)$ is even (odd) under $T_E$; here sufficiency follows from the injectivity of~$L$. Furthermore $\mu$ is unconditional if and only if $L(\mu)$ is unconditional.
\end{remark}

\begin{proposition}
\label{prop - bij meas sphere rn unc}
$L$ is a bijection from $\sM^{\rm int}_{\rm unc}(\R^n)$ to $\sM_{\rm unc}(S^n_0)$.
\end{proposition}

\begin{proof}
This is a consequence of Proposition~\ref{prop - bij meas sphere rn} and Remark~\ref{rema - L commutes refl}.
\end{proof}

For $E \subset [n]$ we define the "lifted" set
\[
E_L = \{0\} \cup E
\]
In particular, $[n]_L = \{0\} \cup [n]$. Some more properties of $L$ are collected in the following proposition, in particular its relationship with decomposition, multiplicative convolution, and orthogonal projections.

\begin{proposition}
\label{prop - lifting}
Let $\mu, \nu \in \sM^{\rm int}(\R^n)$ and $E \subset [n]$.

\begin{thnumber}
\item \label{prop - lifting zero} $\mu = 0$ if and only if $L(\mu) = 0$.
\item \label{prop - lifting decomp} $L(\mu_E) = (L(\mu))_{E_L}$
\item \label{prop - lifting order} $\mu$ is of order $E$ if and only if $L(\mu)$ is of order $E_L$.
\item \label{prop - lifting degree} $\mu$ is of degree $k$ if and only if $L(\mu)$ is of degree $k + 1$.
\item \label{prop - lifting conv} $L(\mu \mc \nu) = L(\mu) \mcs L(\nu)$
\item \label{prop - lifting proj} $L(P_E(\mu)) = P^S_{E_L}(L(\mu))$
\end{thnumber}

\end{proposition}

\begin{proof}
Statement~\eqref{prop - lifting zero} follows from the linearity and injectivity of $L$, see Proposition~\ref{prop - bij meas sphere rn}.

Statement~\eqref{prop - lifting decomp} follows from Propositions~\ref{prop - proj S decomp E} and \ref{rema - mult refl commute}.

Statement~\eqref{prop - lifting order} now follows from \eqref{prop - lifting zero} and \eqref{prop - lifting decomp}.

Statement~\eqref{prop - lifting degree} is a consequence of \eqref{prop - lifting order}.

To see \eqref{prop - lifting conv} note that
\begin{eqnarray}
L(\mu \mc \nu) & = & P_S M_{\rm sym} (\delta_1 \otimes (\mu \mc \nu)) \nonumber\\
 & = & P_S M_{\rm sym} \big( (\delta_1 \otimes \mu) \mc (\delta_1 \otimes \nu) \big) \label{eq - distr lift prod}\\
 & = & P_S \big[ \big( M_{\rm sym}(\delta_1 \otimes \mu) \big) \mc \big( M_{\rm sym}(\delta_1 \otimes \nu) \big) \big] \nonumber\\
 & = & P_S \big( L(\mu) \mc L(\nu) \big) \nonumber\\
 & = & L(\mu) \mcs L(\nu) \nonumber
\end{eqnarray}
where \eqref{eq - distr lift prod} follows from Proposition~\ref{prop - multconv prod}.

We next show \eqref{prop - lifting proj}. Note that
\begin{eqnarray*}
L(P_E(\mu)) & = & P_S M_{\rm sym} (\delta_1 \otimes P_E(\mu)) \\
 & = & P_S M_{\rm sym} P_{E_L}(\delta_1 \otimes \mu) \\
 & = & P_S P_{E_L} M_{\rm sym} (\delta_1 \otimes \mu) \\
 & = & P^S_{E_L}(L(\mu))
\end{eqnarray*}

\end{proof}

\begin{remark}
Proposition~\ref{prop - lifting}~\eqref{prop - lifting decomp} yields the decomposition
\[
L(\mu) = \sum_{E \subset [n]} (L(\mu))_{E_L}
\]
\end{remark}

We now turn to the "lifting" of generating pairs and symmetry pairs. The role of the following definitions becomes clear in Lemma~\ref{le - bij gen supp gen even odd}. For $\sE, \sF \subset \sP_n$ we introduce the notation
\begin{eqnarray*}
\sE_L & = & \big\{ E_L : E \in \sE \big\},\\
\sF^0 & = & \sF \cup \big\{ [n]_L \big\}
\end{eqnarray*}
For $E \subset [n]$ we define the map
\[
\sL_E : \sP^2(E) \longrightarrow \sP^2(E_L),\quad \sL_E(\sF) = \big\{ F,\, F \,\Delta\, E_L \,:\, F \in \sF \big\} 
\]
For $G \in \sL_E(\sF)$ we have $G \in \sF$ if and only if $0 \notin G$. Moreover clearly
\[
\sL_E(\sF_1 \cup \sF_2) \,=\, \sL_E(\sF_1) \cup \sL_E(\sF_2),\quad \sF_1, \sF_2 \subset \sP(E)
\]
If $\sF = \O$, then also $\sL_E(\sF) = \O$. Further we abbreviate $\sL = \sL_{[n]}$. For $\sF \subset \sP_n$ and $E \subset [n]$ we have
\begin{equation}
\label{eq - lifting restr comu}
\sL_E(\sF|_E) \,=\, \sL(\sF)|_{E_L}
\end{equation}

The same notation is used for $\sL_E$ acting on pairs $\sP^2(E) \times \sP^2(E)$, which is defined as
\[
\sL_E(\sF_e, \sF_o) \,=\, (\sL_E(\sF_e), \sL_E(\sF_o))
\]

\begin{remark}
\label{rema - lifting sym pairs}
Let $\sF_e, \sF_o \subset \sP_n$. Set $\sF = \sF_e \cup \sF_o$ and let $\sG$ be the generated group. It is clear that $\sL(\sG)$ is a group and, if $\sF \neq \O$, then $\sL(\sF)$ generates $\sL(\sG)$. Recall that in the context of multiple reflections the assumption that $\sF \neq \O$ is no restriction. Moreover we have, provided again that $\sF \neq \O$,
\begin{equation}
\label{eq - lifting generation sym pair}
\gamma \circ \sL \, (\sF_e, \sF_o) \,=\, \sL \circ \gamma \, (\sF_e, \sF_o) \,=\, \sL (\sG_e, \sG_o) \,=\, \gamma (\sF^0_e, \sF_o)
\end{equation}
The first equality says that lifting commutes with group generation and identification of even and odd parts. It also shows that a lifted symmetry pair is always a symmetry pair. The last equation shows that the lifted symmetry pair can be generated by a smaller generating pair than $\sL \, (\sF_e, \sF_o)$, namely $(\sF^0_e, \sF_o)$.  Finally notice that $(\sF_e, \sF_o)$ is proper if and only if $(\sF^0_e, \sF_o)$ is proper, and this is the case if and only if $\sL \, (\sF_e, \sF_o)$ is proper; to see this observe that all three pairs are proper in the case $\sF = \O$ and use~\eqref{eq - lifting generation sym pair} in the case $\sF \neq \O$.
\end{remark}

\begin{remark}
\label{rema - lifting restriction}
Lifting commutes with restriction in the following sense. Let $E \subset [n]$ and $\sF_e, \sF_o \subset \sP_n$. It follows from \eqref{eq - lifting restr comu} that
\[
\sL_E \circ \rho_E \, (\sF_e, \sF_o) \,=\, \rho_{E_L} \circ \sL \, (\sF_e, \sF_o)
\]
\end{remark}

\begin{lemma}
\label{le - proper lifting}
Let $E \subset [n]$ and $\sF_e, \sF_o \subset \sP_n$. Then $(\sF_e |_E, \sF_o |_E)$ is proper if and only if $(\sF_e^0 |_{E_L}, \sF_o |_{E_L})$ is proper.
\end{lemma}

\begin{proof}
The statement is clear in the case  $\sF_e \cup \sF_o = \O$. Now suppose $\sF_e \cup \sF_o \neq \O$. Then
\begin{eqnarray}
\gamma \circ \sL_E \circ \rho_E \, (\sF_e, \sF_o) & = & \gamma \circ \rho_{E_L} \circ \sL \, (\sF_e, \sF_o) \label{eq - proper lifting 1}\\
 & = & \gamma \circ \gamma \circ \rho_{E_L} \circ \sL \, (\sF_e, \sF_o) \nonumber\\
 & = & \gamma \circ \rho_{E_L} \circ \gamma \circ \sL \, (\sF_e, \sF_o) \label{eq - proper lifting 2}\\
 & = & \gamma \circ \rho_{E_L} \circ \gamma \, (\sF_e^0, \sF_o) \label{eq - proper lifting 3}\\
 & = & \gamma \circ \gamma \circ \rho_{E_L} \, (\sF_e^0, \sF_o) \label{eq - proper lifting 4}\\
 & = & \gamma \circ \rho_{E_L} \, (\sF_e^0, \sF_o) \nonumber
\end{eqnarray}
where \eqref{eq - proper lifting 1} follows from Remark~\ref{rema - lifting restriction}, \eqref{eq - proper lifting 2} and \eqref{eq - proper lifting 4} from Remark~\eqref{rema - proj even odd}, and \eqref{eq - proper lifting 3} from \eqref{eq - lifting generation sym pair}.

Now the last statement in Remark~\eqref{rema - lifting sym pairs} yields that $\rho_E (\sF_e, \sF_o)$ is proper if and only if $\sL_E \circ \rho_E \, (\sF_e, \sF_o)$ is proper. Now the preceding equations show that this is the case if and only if $\rho_{E_L} (\sF_e^0, \sF_o)$ is proper.
\end{proof}

\begin{lemma}
\label{le - bij gen supp gen even odd}
Let $\sE, \sF_e, \sF_o \subset \sP_n$. Then $L$ is a bijective map between $\sM^{\rm int}(\sE; \sF_e, \sF_o)$ and $\sM^{\rm sph}(\sE_L; \sF_e^0, \sF_o)$.
\end{lemma}

\begin{proof}
First note that $L$ is a bijective map between $\sM^{\rm int}(A_E)$ and
\[
\sM(B_{E_L}; \{ [n]_L \}, \O) \,=\, \sM_{\rm sym}(B_{E_L})
\]
for every $E \subset [n]$ by Propositions \ref{prop - bij meas sphere rn} and \ref{prop - lifting}~\eqref{prop - lifting decomp}. Consequently $L$ is a bijective map between $\sM^{\rm int}(\sE)$ and $\sM^{\rm sph}_{\rm sym}(\sE_L)$ for every $\sE \subset \sP_n$. Now Remark~\ref{rema - L commutes refl} implies the claim.
\end{proof}

\section{Uniqueness of measures}
\label{sec - uniqueness}

\subsection{Uniqueness of measures on \texorpdfstring{$\R^n$}{Rn}}

For $x, \alpha \in \R^n$ we define
\[
[x]^\alpha \,=\, \prod_{j=1}^n |x_j|^{\alpha(j)}
\]
Recall the definition of $M^n$ in \eqref{def - m}.

\begin{lemma}
\label{le - abs mom bound on m}
For $x \in \R^n$ and $\alpha \in M^n$ we have
\[
[x]^\alpha \,\leq\, \max \{1, \|x\|_2\} \,\leq\, 1 + \|x\|_2
\]
\end{lemma}

\begin{proof}
\begin{eqnarray*}
 [x]^\alpha & = & \prod_{j=1}^n |x_j|^{\alpha(j)} \nonumber\\
 & \leq & \prod_{j=1}^n \|x\|_2^{\alpha(j)} \nonumber\\
 & = & \|x\|_2^{\|\alpha\|_1} \nonumber\\
 & \leq & \max \{1, \|x\|_2 \}^{\|\alpha\|_1}\nonumber\\
 & \leq & \max \{1, \|x\|_2 \}
\end{eqnarray*}
\end{proof}

\begin{remark}
\label{re - ex abs mom on m}
As a consequence of Lemma~\ref{le - abs mom bound on m} for every $\mu \in \sM^{\rm int}(\R^n)$ and $\alpha \in M^n$ the integral $\int_{\R^n} \mu(dx) \, [x]^\alpha$
exists and is finite.
\end{remark}

The following theorem is well known (see e.g.\ \cite[Lemma~7]{kab:sch:haan09}). A proof is provided here for the convenience of the reader. It uses a standard argument involving only the identity theorem of complex analysis for one complex variable (see for example \cite[Satz~12]{janich} and the uniqueness of the Fourier transform on~$\R^n$.

\begin{lemma}
\label{le - uniqueness laplace}
Let $\mu_1$ and $\mu_2$ be two non-negative measures on $\R^n$. If there is an open $U \subset \R^n$ such that
\[
\int\limits_{\R^n} \exp(\langle x, s \rangle) \,\mu_1(dx) \,=\, \int\limits_{\R^n} \exp(\langle x, s \rangle) \,\mu_2(dx) \,<\, \infty,\quad s \in U,
\]
then $\mu_1 = \mu_2$.
\end{lemma}

\begin{proof}
For $j \in [n]$ choose $a_j, b_j \in \R$ such that $Q \subset U$ where $Q = Q_1 \times \ldots \times Q_n$ and $Q_j = (a_j, b_j)$. Set
\[
G_j = \{ z \in \C \,:\, \real(z) \in Q_j \},\quad j \in [n],
\]
and $G = G_1 \times \ldots \times G_n$. Now assume the stated assumptions and set
\[
f_k(z) = \int\limits_{\R^n} \exp(\langle x, z \rangle) \,\mu_k(dx),\quad z \in G,\; k = 1, 2,
\]
which are then well-defined and finite. Further $f_1(z) = f_2(z)$ for all $z \in Q$. On $G$ the functions $f_k$ are analytic in $z_j$ for each $j \in [n]$ while keeping the other complex coordinates fixed, see \cite[Satz~IV.5.8]{elstrodt}. Now fix $z_j \in Q_j$ for $j \in \{ 2, \ldots n \}$ and note that $f_1(z) = f_2(z)$ holds for all $z_1 \in Q_1$, and therefore for all $z_1 \in G_1$ by the identity theorem. Next let $l \in \{ 2, \ldots, n \}$. Fix $z_j \in G_j$ for $j \in \{ 1, \ldots, l-1 \}$ and, in case $l < n$, $z_j \in Q_j$ for $j \in \{ l+1, \ldots, n \}$. If $f_1(z) = f_2(z)$ holds for $z_l \in Q_l$, it also holds for $z_l \in G_l$ again by identity theorem. Therefore $f_1(z) = f_2(z)$ holds for $z \in G$.

Now choose $s \in Q$. Define the measures $\nu_k(dx) = \mu_k(dx) \exp(\langle x, s \rangle)$ for $k = 1, 2$. Then
\[
\int\limits_{\R^n} \exp(i \langle x, t \rangle) \,\nu_1(dx) \,=\, 
\int\limits_{\R^n} \exp(i \langle x, t \rangle) \,\nu_2(dx),\quad t \in \R^n
\]
The uniqueness of the Fourier transform yields $\nu_1 = \nu_2$. Therefore $\mu_1 = \mu_2$.
\end{proof}

\begin{lemma}
\label{le - kabluchko measures}
Let $\mu_1$ and $\mu_2$ be two non-negative measures on $(0,\infty)^n$. If there is an open $U \subset \R^n$ such that
\begin{displaymath}
\int\limits_{(0,\infty)^n} [x]^\alpha \,\mu_1(dx) \,=\, \int\limits_{(0,\infty)^n} [x]^\alpha \,\mu_2(dx) \,<\, \infty,\quad \alpha \in U,
\end{displaymath}
then $\mu_1 = \mu_2$.
\end{lemma}

\begin{proof}
We define the componentwise transformation
\[
g : (0,\infty)^n \longrightarrow \R^n,\quad g(x) = (\log x_1, \ldots, \log x_n)
\]
Note that $g$ is a homeomorphism. From the assumption it follows that
\[
\int\limits_{(0,\infty)^n} [x]^\alpha \,\mu_k(dx) \,=\, \int\limits_{\R^n} \exp(\langle y, \alpha \rangle) \,g(\mu_k)(dx) \,<\, \infty,\quad \alpha \in U,\; k = 1, 2 
\]
Now Lemma~\ref{le - uniqueness laplace} implies that $g(\mu_1) = g(\mu_2)$. Therefore $\mu_1 = \mu_2$.
\end{proof}

\begin{corollary}
\label{coro - mu zero first conclusions}
Let $U \subset M^n$ such that $U$ is open in~$\R^n$. Then:
\begin{thnumber}
\item \label{coro - mu zero first conclusions 1} For $\mu \in \sM^{\rm int}((0,\infty)^n)$, $\mu = 0$ if and only if $\int_{(0,\infty)^n} \mu(dx) \, [x]^\alpha = 0$ for all $\alpha \in U$.
\item \label{coro - mu zero first conclusions 2} For $\mu \in \sM^{\rm int}_{\rm unc}(A_{[n]})$, $\mu = 0$ if and only if $\int_{A_{[n]}} \mu(dx) \, [x]^\alpha = 0$ for all $\alpha \in U$.
\item \label{coro - mu zero first conclusions 3} For $\mu \in \sM^{\rm int}_{\rm unc}(\R^n)$ with $\mu \geq 0$, $\mu = 0$ holds if and only if $\int_{\R^n} \mu(dx) \, [x]^\alpha = 0$ for all $\alpha \in M^n$.
\end{thnumber}
\end{corollary}

\begin{proof}
Statement~\eqref{coro - mu zero first conclusions 1} follows from
Lemma~\ref{le - kabluchko measures}.

To see statement~\eqref{coro - mu zero first conclusions 2} use \eqref{coro - mu zero first conclusions 1} and the map $M_{\rm unc}$ (or more precisely the bijection $\widetilde{M}$ in Proposition~\ref{prop - bijection unc positive}).

To show \eqref{coro - mu zero first conclusions 3} set
\[
M^E \,=\, \big\{ x \in M^n \;:\; \forall j \notin E \;\; x_j = 0 \big\},\quad E \subset [n]
\]
The condition implies that
\[
\int_{A_E} \mu_E(dx) \, [x]^\alpha = 0,\quad \alpha \in M^E
\]
Applying~\eqref{coro - mu zero first conclusions 2} to the space $\R^E$ yields $\mu_E = 0$.
\end{proof}

\subsection{Uniqueness of measures on \texorpdfstring{$S^{n-1}$}{Sn-1}}

In this section we derive a result for the sphere similar to that of Lemma~\ref{le - kabluchko measures} for~$\R^n$. It is based on the latter lemma and the procedure of "lifting" introduced in Section~\ref{sec - lifting}. Denote
\[
\Delta^n = \left\{ x \in \R^n_+ \,:\, \|x\|_1 = 1 \right\}
\]
A subset of $\Delta^n$ is called {\em relatively open} if it is open with respect to the topology that the standard topology on $\R^n$ induces on the subset $\Delta^n$.

\begin{lemma}
\label{le - ordern sphere mom inj}
Let $\mu \in \sM_{\rm unc}(B_{[n]})$ where $n \geq 2$. If there is a relatively open $U \subset \Delta^n$ such that
\[
\int_{S^{n-1}} [v]^{\alpha} \,\mu(dv) \,=\, 0,\quad \alpha \in U,
\]
then $\mu = 0$.
\end{lemma}

\begin{proof}
We may enumerate coordinates starting from zero. Setting $m = n - 1$, the basic space then consists of points $(x_0, x_1,\ldots, x_m) \in \R^{m+1}$ where $m \geq 1$. Hence $\mu$ is of order $[m]_L$ or zero, and $\mu \in \sM_{\rm unc}(S^m_0)$. Since $L$ is a bijection from $\sM^{\rm int}_{\rm unc}(\R^m)$ to $\sM_{\rm unc}(S^m_0)$ by Proposition~\ref{prop - bij meas sphere rn unc}, there is $\rho \in \sM^{\rm int}_{\rm unc}(\R^m)$ such that $L(\rho) = \mu$. In fact it follows from Proposition~\ref{prop - lifting}~\eqref{prop - lifting order} that $\rho \in \sM^{\rm int}_{\rm unc}(A_{[m]})$. For $\alpha \in \Delta^{m + 1}$ we have
\begin{eqnarray}
\int_{S^m} \mu(dv_0, dv) \, [v]^\alpha & = & \int_{S^m} L(\rho)(dv_0, dv) \, [v]^\alpha \label{eq - first moment sphere} \\
 & = & \int_{S^m} P_S(\delta_1 \otimes \rho)(dv_0, dv) \, [v]^\alpha \nonumber\\
 & = & \int_{\R^{m+1}} \delta_1 \otimes \rho \, (dy_0, dy) \, \prod_{j=0}^m |y_j|^{\alpha(j)} \label{eq - first moment transform}\\
 & = & \int_{\R^m} \rho(dy) \, \prod_{j=1}^m |y_j|^{\alpha(j)}, \label{eq - first moment rm}
\end{eqnarray}
where Proposition~\ref{prop - multconv poshom} is used in~\eqref{eq - first moment transform}. By the assumptions the left-hand side of \eqref{eq - first moment sphere} is zero for $(\alpha_0, \alpha_1, \ldots, \alpha_m) \in U$. Therefore the right-hand side of \eqref{eq - first moment rm} is zero for $(\alpha_1, \ldots, \alpha_m) \in P_{[m]}(U)$. Since $P_{[m]}(U)$ contains an open subset of $\R^m$, Corollary~\ref{coro - mu zero first conclusions}~\eqref{coro - mu zero first conclusions 2} yields $\rho = 0$. Hence $\mu = 0$.
\end{proof}

\subsection{Uniqueness of symmetry decomposition}

For $\nu \in \sM^{\rm int}(A_{[n]})$ we define the function
\[
g(\nu; \,\cdot\,) : M^n \longrightarrow \R,\quad g(\nu; \alpha) = \int_{A_{[n]}} \nu(dx) \, [x]^\alpha
\]
Notice that $g$ is continuous in $\alpha$ on $M^n$ for fixed $\nu$ by bounded convergence and Lemma~\ref{le - abs mom bound on m}. Moreover $g$ is a linear function in~$\nu$. Furthermore we define
\begin{eqnarray*}
D_{\nu} & = & \big\{ \alpha \in M^n : g(\nu; \alpha) \neq 0 \big\} \\
D^{\rm sph}_{\nu} & = & D_{\nu} \cap \Delta^n
\end{eqnarray*}

\begin{remark}
\label{rema - g simple zeros}
Let $\nu \in \sM^{\rm int}_{\rm unc}(A_{[n]})$. Corollary~\ref{coro - mu zero first conclusions} implies that $\nu = 0$ if and only if $g(\nu; \alpha) = 0$ for all $\alpha \in M^n$. This is the case if there is $U \subset M^n$ such that $U$ is open in $\R^n$ and $g(\nu; \alpha) = 0$ for all $\alpha \in U$. Moreover $\nu \neq 0$ if and only if there is $\alpha_0 \in M^n$ such that $g(\nu; \alpha_0) \neq 0$, and in this case $D_{\nu}$ is open and dense in $M^n$.
\end{remark}

\begin{remark}
\label{rema - prod g two measures}
For $\mu, \nu \in \sM^{\rm int}(A_{[n]})$ we have
\begin{equation}
\label{eq - g prod}
g(\mu \mc \nu; \alpha) = g(\mu; \alpha) \, g(\nu; \alpha),\quad \alpha \in M^n
\end{equation}
Therefore
\begin{align*}
& D_{\mu \mc \nu} = D_{\mu} \cap D_{\nu}, & D^{\rm sph}_{\mu \mc \nu} = D^{\rm sph}_{\mu} \cap D^{\rm sph}_{\nu}
\end{align*}
\end{remark}

\begin{remark}
For $\mu \in \sM^{\rm int}(A_{[n]})$ and $\alpha \in \Delta^n$ the integrand of $g(\mu; \alpha)$ is positively one-homogeneous in its integration variable, and consequently
\begin{equation}
\label{eq - g radial proj}
g(P_S(\mu); \alpha) = g(\mu; \alpha),\quad \alpha \in \Delta^n
\end{equation}
\end{remark}

\begin{remark}
\label{rema - g simple zeros sphere}
Let $\nu \in \sM_{\rm unc}(B_{[n]})$. Lemma~\ref{le - ordern sphere mom inj} implies that $\nu = 0$ if and only if $g(\nu; \alpha) = 0$ for all $\alpha \in \Delta^n$. This is the case if there is a relatively open subset $U$ of $\Delta^n$ such that $g(\nu; \alpha) = 0$ for all $\alpha \in U$. Moreover $\nu \neq 0$ if and only if there is $\alpha_0 \in \Delta^n$ such that $g(\nu; \alpha_0) \neq 0$, and in this case $D^{\rm sph}_{\nu}$ is open and dense in $\Delta^n$.
\end{remark}

\begin{remark}
For $\mu, \nu \in \sM(B_{[n]})$ it follows from \eqref{eq - g prod} and \eqref{eq - g radial proj} that 
\[
g(\mu \mcs \nu; \alpha) \,=\, g(\mu; \alpha) \, g(\nu; \alpha),\quad \alpha \in \Delta^n
\]
\end{remark}

\begin{proposition}
\label{prop - prod an zero}
Let $\mu, \nu \in \sM^{\rm int}_{\rm unc}(A_{[n]})$. Then $\mu \mc \nu \neq 0$ if and only if $\mu \neq 0$ and $\nu \neq 0$.
\end{proposition}

\begin{proof}
Suppose $\mu \neq 0$ and $\nu \neq 0$. Then $D_{\mu}$ and $D_{\nu}$ are open and dense in $M^n$. Thus $D_{\mu \mc \nu}$ is open and dense in $M^n$, implying $\mu \mc \nu \neq 0$.

Conversely, if either $\mu = 0$ or $\nu = 0$, then clearly $\mu \mc \nu = 0$.
\end{proof}

\begin{proposition}
\label{prop - prod bn zero}
Let $\mu, \nu \in \sM_{\rm unc}(B_{[n]})$. Then $\mu \mcs \nu \neq 0$ if and only if $\mu \neq 0$ and $\nu \neq 0$.
\end{proposition}

\begin{proof}
Suppose $\mu \neq 0$ and $\nu \neq 0$. Then $D_{\mu}^{\rm sph}$ and $D_{\nu}^{\rm sph}$ are open and dense in $\Delta^n$. Thus $D_{\mu \mc \nu}^{\rm sph}$ is open and dense in $\Delta^n$. Hence $\mu \mcs \nu \neq 0$.

Conversely, if either $\mu = 0$ or $\nu = 0$, then clearly $\mu \mcs \nu = 0$.
\end{proof}

\begin{proposition}
\label{prop - deltaj mu nu zero}
Let $\mu, \nu \in \sM^{\rm int}(A_{[n]})$ and $J \subset [n]$. Then $\delta_J \mc \mu \mc \nu \neq 0$ holds if and only if $\delta_J \mc \mu \neq 0$ and $\delta_J \mc \nu \neq 0$.
\end{proposition}

\begin{proof}
First note that
\begin{eqnarray}
\sigma_J \odot \big( \delta_J \mc (\mu \mc \nu) \big) & = & \sigma_J \odot \big( (\delta_J \mc \mu) \mc (\delta_J \mc \nu) \big) \label{eq - dens prod deltaj}\\
 & = & \big(\sigma_J \odot (\delta_J \mc \mu) \big) \mc \big(\sigma_J \odot (\delta_J \mc \nu) \big) \nonumber
\end{eqnarray}
Next notice that
\begin{eqnarray*}
\delta_J \mc \mu \neq 0 & \Longleftrightarrow & \sigma_J \odot (\delta_J \mc \mu) \neq 0, \\
\delta_J \mc \nu \neq 0 & \Longleftrightarrow & \sigma_J \odot (\delta_J \mc \nu) \neq 0, \\
\delta_J \mc \mu \mc \nu \neq 0 & \Longleftrightarrow & \sigma_J \odot \big( \delta_J \mc \mu \mc \nu \big) \neq 0
\end{eqnarray*}
By Lemma~\ref{le - density conv}~\eqref{le - density conv 2} the left-hand side and both factors on the right-hand side in \eqref{eq - dens prod deltaj} are unconditional. Now the assertion follows from Proposition~\ref{prop - prod an zero}.
\end{proof}

Notice that in Proposition~\ref{prop - deltaj mu nu zero} integrability of $\mu$ and $\nu$ is required for the application of Proposition~\ref{prop - prod an zero} because the function $g$ should be well-defined.

\begin{proposition}
\label{prop - deltaj mu nu zero sphere}
Let $\mu, \nu \in \sM(B_{[n]})$ and $J \subset [n]$. Then $\delta_J \mcs \mu \mcs \nu \neq 0$ holds if and only if $\delta_J \mcs \mu \neq 0$ and $\delta_J \mcs \nu \neq 0$.
\end{proposition}

\begin{proof}
First note that
\begin{eqnarray}
\sigma_J \odot \big( \delta_J \mcs (\mu \mcs \nu) \big) & = & \sigma_J \odot \big( (\delta_J \mcs \mu) \mcs (\delta_J \mcs \nu) \big) \label{eq - dens prod deltaj sphere}\\
 & = & \big(\sigma_J \odot (\delta_J \mcs \mu) \big) \mcs \big(\sigma_J \odot (\delta_J \mcs \nu) \big) \nonumber
\end{eqnarray}
Next notice that 
\begin{eqnarray*}
\delta_J \mcs \mu \neq 0 & \Longleftrightarrow & \sigma_J \odot (\delta_J \mcs \mu) \neq 0,\\
\delta_J \mcs \nu \neq 0 & \Longleftrightarrow & \sigma_J \odot (\delta_J \mcs \nu) \neq 0,\\
\delta_J \mcs \mu \mcs \nu \neq 0 & \Longleftrightarrow & \sigma_J \odot \big( \delta_J \mcs \mu \mcs \nu \big) \neq 0
\end{eqnarray*}
By Lemma~\ref{le - density conv}~\eqref{le - density conv 2} the left-hand side and both factors on the right-hand side in \eqref{eq - dens prod deltaj sphere} are unconditional. Now the assertion follows from Proposition~\ref{prop - prod bn zero}.
\end{proof}

\section{Universality on \texorpdfstring{$\R^n$}{Rn}}
\label{sec - universality Rn}

\subsection{General case}

Let $\nu \in \sM(\R^n)$ and $\scrA \subset \sM(\R^n)$. Then $\nu$ is called {\em universal on} $\scrA$ if, for every $\mu \in \scrA$, $\nu \mc \mu = 0$ implies $\mu = 0$.

\begin{proposition}
Let $\sE, \sF_e, \sF_o \subset \sP_n$ such that $\sE$ is closed under finite intersections, and for $j \in \{1, 2\}$ let $\nu_j \in \sM(\sE)$. If $\nu_1$ and $\nu_2$ are universal on $\sM(\sE; \sF_e, \sF_o)$, then $\nu_1 \mc \nu_2$ is universal on $\sM(\sE; \sF_e, \sF_o)$.
\end{proposition}

\begin{proof}
Assume that $\nu_1$ and $\nu_2$ are universal on $\sM(\sE; \sF_e, \sF_o)$. Let $\mu \in \sM(\sE; \sF_e, \sF_o)$ such that $\nu_1 \mc \nu_2 \mc \mu = 0$. By Proposition~\ref{prop - mult conv even odd support} $\nu_2 \mc \mu \in \sM(\sE; \sF_e, \sF_o)$, and therefore $\nu_2 \mc \mu = 0$ by the assumption. Hence also $\nu_2 = 0$.
\end{proof}

\begin{theorem}
\label{th - measures univ general even odd ae}
Let $\nu \in \sM^{\rm int}(A_{[n]})$ and $\sF_e, \sF_o \subset \sP_n$ such that $(\sF_e, \sF_o)$ is a proper generating pair. Then $\nu$ is universal on $\sM^{\rm int}(A_{[n]}; \sF_e, \sF_o)$ if and only if $\delta_J \mc \nu \neq 0$ for every $J \in \sJ(\sF_e, \sF_o)$. In this case, if $\nu$ is even (odd) under $T_E$ for some $E \subset [n]$, then $E \in \sG_e$ ($E \in \sG_o$) where $\sG$ is the group generated by~$\sF_e \cup \sF_o$.
\end{theorem}

\begin{proof}
{\em Sufficiency.} Let $\mu \in \sM^{\rm int}(A_{[n]}; \sF_e, \sF_o)$ such that $\nu \mc \mu = 0$. Let $J \subset [n]$. If $J \in \sJ(\sF_e, \sF_o)$, then  Proposition~\ref{prop - deltaj mu nu zero} implies that $\delta_J \mc \mu = 0$ since $\delta_J \mc \nu \neq 0$ by the assumptions. If $J \notin \sJ(\sF_e, \sF_o)$, then Corollary~\ref{coro - deltaj zero even odd} yields $\delta_J \mc \mu = 0$. Now it follows from Proposition~\ref{prop - sign meas decomp signs rn}~\eqref{prop - sign meas decomp signs rn decomp} that $\mu = 0$.

{\em Necessity.} Let $J \in \sJ(\sF_e, \sF_o)$. Then $\delta_J \in \sM^{\rm int}(A_{[n]}; \sF_e, \sF_o)$ by Remark~\ref{rema - j e properties} and clearly $\delta_J \neq 0$. If $\nu$ is universal on $\sM^{\rm int}(A_{[n]}; \sF_e, \sF_o)$, then $\delta_J \mc \nu \neq 0$.

{\em Supplement.} We know from Proposition~\ref{prop - j gen group even odd}~\eqref{prop - j gen group even odd 1} that $\sJ(\sF_e, \sF_o) = \sJ(\sG_e, \sG_o)$. Now suppose that $\delta_J \mc \nu \neq 0$ holds for every $J \in \sJ(\sG_e, \sG_o)$ and that $\nu$ is even (odd) under $T_E$. By Theorem~\ref{th - deltaj zero arb e sym} we have $\delta_J \mc \nu = 0$ for every $J \subset [n]$ such that $J \in \sJ(\O, E)$ ($J \in \sJ(E, \O)$). It follows that $\sJ(\sG_e, \sG_o) \subset \sJ(E, \O)$ ($\sJ(\sG_e, \sG_o) \subset \sJ(\O, E)$). Now Lemma~\ref{le - j s g impl} yields $E \in \sG_e$ ($E \in \sG_o$).
\end{proof}

Note that the last statement in Theorem~\ref{th - measures univ general even odd ae} says that $\nu$ cannot be universal on $\sM^{\rm int}(A_{[n]}; \sG_e, \sG_o)$ if it has itself reflection symmetries, even or odd, other than those contained in $\sG_e$ or $\sG_o$, respectively. An alternative proof of this fact that uses Proposition~\ref{prop - sym pair symmetrisation} directly is the following.

\begin{proof}[Second proof of supplement in Theorem~\ref{th - measures univ general even odd ae}]
Suppose that $\nu$ is universal on $\sM^{\rm int}(A_{[n]}; \sF_e, \sF_o)$ and that it is even (odd) under $T_E$. Set $\rho =  M_{\sF_e, \sF_o}(\delta_{1_n})$. By Proposition~\ref{prop - sym pair symmetrisation} we have $\rho \in \sM^{\rm int}(A_{[n]}; \sG_e, \sG_o)$. Clearly also $T_E(\rho) \in \sM^{\rm int}(A_{[n]}; \sG_e, \sG_o)$. In the "even" case we have
\[
\rho \mc \nu \,=\, \rho \mc T_E(\nu) \,=\, T_E(\rho) \mc \nu
\]
It follows that $T_E(\rho) = \rho$. Thus $E \in \sG_e$ by Proposition~\ref{prop - sym pair symmetrisation}~\eqref{prop - sym pair symmetrisation 4}. The argument for the "odd" case is similar.
\end{proof}

\begin{proposition}
\label{prop - nu proj univ even odd}
Let $\sF_e, \sF_o \subset \sP_n$ and $E, G \subset [n]$ with $G \subset E$. Further let $\scrA \subset \sM(A_G; \sF_e, \sF_o)$ and $\nu \in \sM(\R^n)$. Then $\nu$ is universal on~$\scrA$ if and only if $P_E(\nu)$ is universal on~$\scrA$.
\end{proposition}

\begin{proof}
This follows from the fact that
\[
P_E(\nu) \mc \mu \,=\, \nu \mc P_E(\mu) \,=\, \nu \mc \mu,\quad \mu \in \sM(A_G; \sF_e, \sF_o)
\]
\end{proof}

\begin{lemma}
\label{le - univ gen even odd gen supp}
Let $\nu \in \sM^{\rm int}(\R^n)$ and $\sF_e, \sF_o \subset \sP_n$ such that $(\sF_e, \sF_o)$ is a proper generating pair. If $\nu$ is universal on $\sM^{\rm int}(A_{[n]}; \sF_e, \sF_o)$, then also $\nu_{[n]}$ is universal on $\sM^{\rm int}(A_{[n]}; \sF_e, \sF_o)$.
\end{lemma}

\begin{proof}
Let $J \in \sJ(\sF_e, \sF_o)$. Then Corollary~\ref{coro - sigma0 n term}~\eqref{coro - sigma0 n term rn} implies that
\begin{eqnarray*}
\nu \mc \delta_J \mc \sigma_0^n & = & (\nu \mc \delta_J)_{[n]} \mc \sigma_0^n \\
 & = & \nu_{[n]} \mc \delta_J \mc \sigma_0^n
\end{eqnarray*}
From Lemma~\ref{le - deltaj even odd class}~\eqref{le - deltaj even odd class 1} we know that
\[
0 \,\neq\, \delta_J \mc \sigma_0^n \,\in\, \sM^{\rm int}(A_{[n]}; \sF_e, \sF_o)
\]
Hence, if $\nu$ is universal on $\sM^{\rm int}(A_{[n]}; \sF_e, \sF_o)$,  then $\nu_{[n]} \mc \delta_J \mc \sigma_0^n \neq 0$, and thus $\nu_{[n]} \mc \delta_J \neq 0$. Consequently $\nu_{[n]}$ is universal on $\sM^{\rm int}(A_{[n]}; \sF_e, \sF_o)$ by Theorem~\ref{th - measures univ general even odd ae}.
\end{proof}

\begin{theorem}
\label{th - gen even odd supp char}
Let $\sE, \sF_e, \sF_o \subset \sP_n$ and $\nu \in \sM^{\rm int}(\R^n)$. Then $\nu$ is universal on $\sM^{\rm int}(\sE; \sF_e, \sF_o)$ if and only if
\begin{equation}
\label{th - gen even odd supp char gen 1}
\delta_J \mc R_E P_E(\nu) \neq 0,\quad E \in \sE,\; J \in \sJ(E; \sF_e, \sF_o)
\end{equation}
In this case, if $( \sF_e |_E, \sF_o |_E)$ is proper and $R_E P_E(\nu)$ is even (odd) under $T_G$ for some $E \in \sE$ and $G \subset E$, then $G \in \sG_e |_E$ ($G \in \sG_o |_E$) where $\sG$ is the group generated by $\sF_e \cup \sF_o$.
\end{theorem}

\begin{remark}
In the case $[n] \in \sE$, if $(\sF_e, \sF_o)$ is proper and $\nu$ is universal on $\sM^{\rm int}(\sE; \sF_e, \sF_o)$ and even (odd) under $T_G$ for some $G \subset [n]$, then also $\nu_{[n]}$ is even (odd) under $T_G$ and the second statement in Theorem~\ref{th - gen even odd supp char} says that $G \in \sG_e$ ($G \in \sG_o$).
\end{remark}

\begin{remark}
\label{rema - non proper conditions}
Notice that in Theorem~\ref{th - gen even odd supp char} if $(\sF_e |_E, \sF_o |_E)$ is not proper for some $E \in \sE$, then $E$ can be omitted in $\sE$ because then
\[
\sM^{\rm int}(\sE; \sF_e, \sF_o) = \sM^{\rm int}(\sE \setminus \{E\}; \sF_e, \sF_o)
\]
Note also that for such $E$ there are no conditions \eqref{th - gen even odd supp char gen 1} as $\sJ(E; \sF_e, \sF_o) = \O$ by Proposition~\ref{prop - j gen group even odd}~\eqref{prop - j gen group even odd 2}.
\end{remark}

\begin{proof}[Proof of Theorem~\ref{th - gen even odd supp char}]
{\em Sufficiency.} Let $\mu \in \sM^{\rm int}(\sE; \sF_e, \sF_o)$ such that $\nu \mc \mu = 0$. If $[n] \notin \sE$, then $\mu_{[n]} = 0$. If $[n] \in \sE$, condition~\eqref{th - gen even odd supp char gen 1} for $E = [n]$ says that $\delta_J \mc \nu_{[n]} \neq 0$ for every $J \in \sJ(\sF_e, \sF_o)$. Since multiple reflections commute with $R_{[n]}$, we have $\mu_{[n]} \in \sM^{\rm int}(A_{[n]}; \sF_e, \sF_o)$. Theorem~\ref{th - measures univ general even odd ae} implies that $\nu_{[n]}$ is universal on~$\sM^{\rm int}(A_{[n]}; \sF_e, \sF_o)$. Now note that
\[
\nu_{[n]} \mc \mu_{[n]} \,=\, (\nu \mc \mu)_{[n]} \,=\, 0
\]
Therefore $\mu_{[n]} = 0$.

Now let $k \in \{0, 1, \ldots, n - 1\}$ and assume that $\mu_M = 0$ for all $M \subset [n]$ with $|M| > k$. Let $E \in \sE$ such that $|E| = k$. We have
\begin{eqnarray}
0 & = & (\nu \mc \mu)_E \nonumber\\
 & = & (\nu \mc \mu_E)_E \label{eq - even odd gen supp 1}\\
 & = & \big( P_E(\nu) \mc \mu_E \big)_E \label{eq - even odd gen supp 2}\\
 & = & R_E P_E(\nu) \mc \mu_E \label{eq - even odd gen supp 3}
\end{eqnarray}
where \eqref{eq - even odd gen supp 1} and \eqref{eq - even odd gen supp 3} follow from Corollary~\ref{cor - multconv decomp}, and \eqref{eq - even odd gen supp 2} from Proposition~\ref{prop - multconv proj}. By the assumption the first factor in \eqref{eq - even odd gen supp 3} satisfies
\begin{equation}
\label{eq - delta diff 0 crit}
\delta_{E, J} \mc R_E P_E(\nu) \,=\, \delta_J \mc R_E P_E(\nu) \,\neq\, 0,\quad J \in \sJ(E; \sF_e, \sF_o)
\end{equation}
Theorem~\ref{th - measures univ general even odd ae} implies that $R_E P_E(\nu)$ is universal on~$\sM^{\rm int}(A_E; \sF_e, \sF_o)$. Moreover we know that $\mu_E \in \sM^{\rm int}(A_E; \sF_e, \sF_o)$. Consequently $\mu_E = 0$.

{\em Necessity.} Assume that $\nu$ is universal on $\sM^{\rm int}(\sE; \sF_e, \sF_o)$ and let $E \in \sE$. By Proposition~\ref{prop - nu proj univ even odd} $P_E(\nu)$ is universal on $\sM^{\rm int}(A_E; \sF_e, \sF_o)$. Lemma~\ref{le - univ gen even odd gen supp} implies that also $R_E P_E(\nu)$ is universal on $\sM^{\rm int}(A_E; \sF_e, \sF_o)$. Therefore \eqref{eq - delta diff 0 crit} holds by Theorem~\ref{th - measures univ general even odd ae}.

{\em Supplement.} Let $E \in \sE$. Under the stated conditions $R_E P_E(\nu)$ is universal on
\[
\sM^{\rm int}(A_E; \sF_e, \sF_o) \,=\, \sM^{\rm int}(A_E; \sF_e |_E, \sF_o |_E)
\]
Recall from Remark~\ref{rema - proj even odd} that $\sG_e |_E$ and $\sG_o |_E$ are the even and odd parts, respectively, of the group generated by the pair $(\sF_e |_E, \sF_o |_E)$. Now it follows from the last statement of Theorem~\ref{th - measures univ general even odd ae} that if $R_E P_E(\nu)$ is even (odd) under $T_G$ for some $G \subset E$, then $G \in \sG_e |_E$ ($G \in \sG_o |_E$).
\end{proof}

\begin{corollary}
\label{coro - measures gen even odd univ rn char}
Let $\sF_e, \sF_o \subset \sP_n$ and $\nu \in \sM^{\rm int}(\R^n)$. Then:
\begin{thnumber}
\item \label{coro - measures gen even odd univ rn char an} $\nu$ is universal on $\sM^{\rm int}(A_{[n]}; \sF_e, \sF_o)$ if and only if
\begin{equation*}
\label{coro - measures gen even odd univ rn char an 1}
\delta_J \mc \nu_{[n]} \neq 0,\quad J \in \sJ(\sF_e, \sF_o)
\end{equation*}
\item \label{coro - measures gen even odd univ rn char gen} $\nu$ is universal on $\sM^{\rm int}(\R^n; \sF_e, \sF_o)$ if and only if
\begin{equation}
\label{coro - measures gen even odd univ rn char gen 1}
\delta_J \mc R_E P_E(\nu) \neq 0,\quad E \subset [n],\; J \in \sJ(E; \sF_e, \sF_o)
\end{equation}
\item \label{coro - measures gen even odd univ rn char order-n} For $\nu$ of order $[n]$, $\nu$ is universal on $\sM^{\rm int}(\R^n; \sF_e, \sF_o)$ if and only if
\[
\delta_{J, J} \mc \nu \neq 0,\quad J \in \sJ(\sF_e, \sF_o)
\]
In this case, if $(\sF_e |_E, \sF_o |_E)$ is proper and $P_E(\nu)$ is even (odd) under $T_G$ for some $E \subset [n]$ and $G \subset E$, then $G \in \sG_e |_E$ ($G \in \sG_o |_E$) where $\sG$ is the group generated by $\sF_e \cup \sF_o$.
\end{thnumber}
\end{corollary}

\begin{proof}
Statements~\eqref{coro - measures gen even odd univ rn char an} and \eqref{coro - measures gen even odd univ rn char gen} follow from Theorem~\ref{th - gen even odd supp char} by setting $\sE = \{ [n] \}$ and $\sE = \sP_n$, respectively.

To see \eqref{coro - measures gen even odd univ rn char order-n} assume that $\nu$ is of order $[n]$. In this case for $J \subset E \subset [n]$ we have
\begin{equation*}
\label{eq - nu order n int gen even odd}
\delta_{E, J} \mc R_E P_E(\nu) \,=\, \delta_{E, J} \mc P_E(\nu) \,=\, P_E(\delta_{E, J}) \mc \nu \,=\, \delta_{E, J} \mc \nu
\end{equation*}
and
\begin{equation}
\label{eq - coro nu order n int gen even odd proj j}
P_J(\delta_{E, J} \mc \nu) \,=\, P_J(\delta_{E, J}) \mc \nu \,=\, \delta_{J, J} \mc \nu
\end{equation}
It follows from \eqref{coro - measures gen even odd univ rn char gen} that $\nu$ is universal on $\sM^{\rm int}(\R^n; \sF_e, \sF_o)$ if and only if
\begin{equation}
\label{eq - coro delta e j gen even odd}
\delta_{E, J} \mc \nu \neq 0,\quad E \subset [n],\; J \in \sJ(E; \sF_e, \sF_o)
\end{equation}
Specifically \eqref{eq - coro delta e j gen even odd} implies that $\delta_{J, J} \mc \nu \neq 0$ for every $J \in \sJ(\sF_e, \sF_o)$. On the other hand it follows from~\eqref{eq - coro nu order n int gen even odd proj j} that if $\delta_{E, J} \mc \nu = 0$ for some $J$ and $E$, then $\delta_{J, J} \mc \nu = 0$. The second statement in~\eqref{coro - measures gen even odd univ rn char order-n} follows directly from the second statement of Theorem~\ref{th - gen even odd supp char}.
\end{proof}

\begin{example}
\label{ex - deltaej universal}
Let $\nu \in \sM^{\rm int}(\R^n)$ and $(\sG_e, \sG_o)$ be a proper symmetry pair. It follows from Corollary~\ref{coro - measures gen even odd univ rn char}~\eqref{coro - measures gen even odd univ rn char an} that $\nu$ can only be universal on $\sM^{\rm int}(A_{[n]}; \sG_e, \sG_o)$ if it has degree~$n$. Now let $J \subset E \subset [n]$. Assume that $\delta_{E, J}$ is universal on $\sM^{\rm int}(A_{[n]}; \sG_e, \O)$. Then necessarily $E = [n]$ as $\delta_{E, J}$ must have degree $n$. Moreover Proposition~\ref{prop - prod delta j and k} and Corollary~\ref{coro - measures gen even odd univ rn char}~\eqref{coro - measures gen even odd univ rn char an} imply that $\sJ(\sG_e, \O)$ is a singleton, viz.~$\{J\}$, and thus, since it is a group by Theorem~\ref{th - j trafo group}, we must have $\sJ(\sG_e, \O) = \{ \O \}$. Consequently $J = \O$. Further it follows from Theorem~\ref{th - j trafo group} that in this case $\sG_e = \sP_n$.
\end{example}

\subsection{Specific symmetries}

Recall that unconditionality corresponds to $(\sG_e, \sG_o) = (\sP_n, \O)$. In this case we have $\sJ(\sG_e, \sG_o) = \{ \O \}$, and therefore unconditional measures do have only a single term in their symmetry decomposition (see Proposition~\ref{prop - sign meas decomp signs rn}~\eqref{prop - sign meas decomp signs rn decomp}) by Corollary~\ref{coro - deltaj zero even odd}.

\begin{theorem}[Unconditional measures]
\label{th - measures unc univ char}
Let $\nu \in \sM^{\rm int}(\R^n)$. Then:
\begin{thnumber}
\item \label{th - measures unc univ char ae} $\nu$ is universal on $\sM^{\rm int}_{\rm unc}(A_{[n]})$ if and only if $M_{\rm unc} (\nu)$ has degree $n$.
\item \label{th - measures unc univ char gen} $\nu$ is universal on $\sM^{\rm int}_{\rm unc}(\R^n)$ if and only if $R_E P_E M_{\rm unc}(\nu) \neq 0$ for every $E \subset [n]$.
\item \label{th - measures unc univ char order-n} For $\nu$ of order $[n]$, $\nu$ is universal on $\sM^{\rm int}_{\rm unc}(\R^n)$ if and only if $\nu(\R^n) \neq 0$.
\item \label{th - measures unc univ char pos} For non-negative $\nu$, $\nu$ is universal on $\sM^{\rm int}_{\rm unc}(\R^n)$ if and only if it has degree~$n$.
\end{thnumber}
\end{theorem}

\begin{proof}
Each of the statements \eqref{th - measures unc univ char ae} to \eqref{th - measures unc univ char order-n} follows from the corresponding statement in Corollary~\ref{coro - measures gen even odd univ rn char}. Finally note that \eqref{th - measures unc univ char gen} implies \eqref{th - measures unc univ char pos}.
\end{proof}

\begin{remark}
Notice that the conditions in the case of signed measures \eqref{th - measures unc univ char gen} are much stronger than in the case of non-negative measures~\eqref{th - measures unc univ char pos}.
\end{remark}

\begin{corollary}[Measures on $\R^n_+$]
\label{coro - univ rn plus}
Let $\nu \in \sM^{\rm int}(\R^n)$. If $\nu$ is universal on $\sM^{\rm int}_{\rm unc}(\R^n)$, then $\nu$ is universal on $\sM^{\rm int}(\R^n_+)$.
\end{corollary}

\begin{proof}
Let $\mu \in \sM^{\rm int}(\R^n_+)$ such that $\nu \mc \mu = 0$. Then
\[
0 \,=\, M_{\rm unc}(\nu \mc \mu) \,=\, \nu \mc M_{\rm unc}(\mu) 
\]
By the assumption we obtain $M_{\rm unc}(\mu) = 0$, and Proposition~\ref{prop - bijection unc positive} implies $\mu = 0$.
\end{proof}

Next we consider the origin-symmetric case, that is $\sG_e = \{ \O, [n] \}$ and $\sG_o = \O$, and hence
\[
\sJ(\sG_e, \sG_o) = \{ J \subset [n] \,:\, |J| \;\, {\rm is \; even} \,\}
\]
Hence the symmetry decomposition of origin-symmetric measures contains only non-zero terms such that $|J|$ is even.

\begin{theorem}[Symmetric measures]
\label{th - measures sym univ rn char}
Let $\nu \in \sM^{\rm int}(\R^n)$. Then:
\begin{thnumber}
\item \label{th - measures sym univ rn char ae} $\nu$ is universal on $\sM^{\rm int}_{\rm sym}(A_{[n]})$ if and only if
\[
\delta_J \mc \nu_{[n]} \neq 0,\quad J \subset [n],\; |J| \; {\it even}
\]
\item \label{th - measures sym univ rn char gen} $\nu$ is universal on $\sM^{\rm int}_{\rm sym}(\R^n)$ if and only if
\[
\delta_J \mc R_E P_E(\nu) \neq 0,\quad J \subset E \subset [n],\; |J| \; {\it even}
\]
\item \label{th - measures sym univ rn char order-n} For $\nu$ of order $[n]$, $\nu$ is universal on $\sM^{\rm int}_{\rm sym}(\R^n)$ if and only if
\[
\delta_{J, J} \mc \nu \neq 0,\quad J \subset [n],\; |J| \; {\it even}
\]
In this case, if $P_E(\nu)$ is even under $T_G$ for some $E \subset [n]$ and $G \subset E$, then $G = \O$ or $G = E$.
\end{thnumber}
\end{theorem}

\begin{proof}
Each of the statements \eqref{th - measures sym univ rn char ae} to \eqref{th - measures sym univ rn char order-n} follows from the corresponding statement in Corollary~\ref{coro - measures gen even odd univ rn char}.
\end{proof}

Another interesting special case is $\sF_e = \sF_o = \O$, and thus $\sJ(\O, \O) = \sP_n$. This corresponds to not imposing any symmetry under multiple reflections. In this case each term in the symmetry decomposition of a measure may be non-zero.

\begin{theorem}[No reflection symmetry]
\label{th - measures general univ rn char}
Let $\nu \in \sM^{\rm int}(\R^n)$. Then:
\begin{thnumber}
\item \label{th - measures general univ rn char ae} $\nu$ is universal on $\sM^{\rm int}(A_{[n]})$ if and only if
\[
\delta_J \mc \nu_{[n]} \neq 0,\quad J \subset [n]
\]
\item \label{th - measures general univ rn char gen} $\nu$ is universal on $\sM^{\rm int}(\R^n)$ if and only if
\[
\delta_J \mc R_E P_E(\nu) \neq 0,\quad J \subset E \subset [n]
\]
\item \label{th - measures general univ rn char order-n} For $\nu$ of order $[n]$, $\nu$ is universal on $\sM^{\rm int}(\R^n)$ if and only if
\[
\delta_{J, J} \mc \nu \neq 0,\quad J \subset [n]
\]
In this case, if $P_E(\nu)$ is even under $T_G$ for some $E \subset [n]$ and $G \subset E$, then $G = \O$.
\end{thnumber}
\end{theorem}

\begin{proof}
Each of the statements \eqref{th - measures general univ rn char ae} to \eqref{th - measures general univ rn char order-n} follows from the corresponding statement in Corollary~\ref{coro - measures gen even odd univ rn char}.
\end{proof}

\begin{example}
\label{ex - sigma0n univ}
Since $\delta_J \mc \sigma_0^n \neq 0$ holds for every $J \subset [n]$ by Lemma~\ref{le - deltaj even odd class}~\eqref{le - deltaj even odd class 1}, the condition of Theorem~\ref{th - measures general univ rn char}~\eqref{th - measures general univ rn char ae} is satisfied for $\nu = \sigma_0^n$. Therefore $\sigma_0^n$ is universal on $\sM^{\rm int}(A_{[n]})$.
\end{example}

At last we consider the anti-symmetric case, that is $(\sG_e, \sG_o) = (\{ \O \}, \{ [n] \})$. In this case
\[
\sJ(\sG_e, \sG_o) = \{ J \subset [n] \,:\, |J| \;\, {\rm is \; odd} \,\}
\]
Therefore according to Corollary~\ref{coro - deltaj zero even odd} the symmetry decomposition (see Proposition~\ref{prop - sign meas decomp signs rn}~\eqref{prop - sign meas decomp signs rn decomp}) of anti-symmetric measures contains only non-zero terms such that $|J|$ is odd.

\begin{theorem}[Anti-symmetric measures]
\label{th - measures asym univ rn char}
Let $\nu \in \sM^{\rm int}(\R^n)$. Then:
\begin{thnumber}
\item \label{th - measures asym univ rn char ae} $\nu$ is universal on $\sM^{\rm int}_{\rm asym}(A_{[n]})$ if and only if
\[
\delta_J \mc \nu_{[n]} \neq 0,\quad J \subset [n],\; |J| \; {\it odd}
\]
\item \label{th - measures asym univ rn char gen} $\nu$ is universal on $\sM^{\rm int}_{\rm asym}(\R^n)$ if and only if
\[
\delta_J \mc R_E P_E(\nu) \neq 0,\quad J \subset E \subset [n],\; |J| \; {\it odd}
\]
\item \label{th - measures asym univ rn char order-n} For $\nu$ of order $[n]$, $\nu$ is universal on $\sM^{\rm int}_{\rm asym}(\R^n)$ if and only if
\[
\delta_{J, J} \mc \nu \neq 0,\quad J \subset [n],\; |J| \; {\it odd}
\]
In this case, if $P_E(\nu)$ is even (odd) under $T_G$ for some $E \subset [n]$ and $G \subset E$, then $G = \O$ ($G = E$).
\end{thnumber}
\end{theorem}

\begin{proof}
Each of the statements \eqref{th - measures asym univ rn char ae} to \eqref{th - measures asym univ rn char order-n} follows from the corresponding statement in Corollary~\ref{coro - measures gen even odd univ rn char}.
\end{proof}

\section{Universality on the sphere}
\label{sec - universality sphere}

\subsection{General case}

Let $\nu \in \sM(S^{n-1})$ and $\scrA \subset \sM(S^{n-1})$. Then $\nu$ is called {\em spherically universal on} $\scrA$ if, for every $\mu \in \scrA$, $\nu \mcs \mu = 0$ implies $\mu = 0$.

\begin{proposition}
Let $\sE, \sF_e, \sF_o \subset \sP_n$ such that $\O \notin \sE$ and $\sE \cup \{\O\}$ is closed under finite intersections, and for $j \in \{1, 2\}$ let $\nu_j \in \sM^{\rm sph}(\sE)$. If $\nu_1$ and $\nu_2$ are spherically universal on $\sM^{\rm sph}(\sE; \sF_e, \sF_o)$, then $\nu_1 \mcs \nu_2$ is spherically universal on $\sM^{\rm sph}(\sE; \sF_e, \sF_o)$.
\end{proposition}

\begin{proof}
Assume that $\nu_1$ and $\nu_2$ are spherically universal on $\sM^{\rm sph}(\sE; \sF_e, \sF_o)$. Let $\mu \in \sM^{\rm sph}(\sE; \sF_e, \sF_o)$ such that $\nu_1 \mcs \nu_2 \mcs \mu = 0$. By Proposition~\ref{prop - mult conv even odd support sphere} $\nu_2 \mcs \mu \in \sM^{\rm sph}(\sE; \sF_e, \sF_o)$, and therefore $\nu_2 \mcs \mu = 0$ by the assumption. Hence also $\nu_2 = 0$.
\end{proof}

\begin{theorem}
\label{th - measures univ general even odd be}
Let $\nu \in \sM(B_{[n]})$ and $\sF_e, \sF_o \subset \sP_n$ such that $(\sF_e, \sF_o)$ is a proper generating pair. Then $\nu$ is spherically universal on $\sM(B_{[n]}; \sF_e, \sF_o)$ if and only if $\delta_J \mcs \nu \neq 0$ for every $J \in \sJ(\sF_e, \sF_o)$. In this case, if $\nu$ is even (odd) under $T_E$ for some $E \subset [n]$, then $E \in \sG_e$ ($E \in \sG_o$) where $\sG$ is the group generated by~$\sF_e \cup \sF_o$.
\end{theorem}

\begin{proof}
{\em Sufficiency.} Let $\mu \in \sM(B_{[n]}; \sF_e, \sF_o)$. Assume $\nu \mcs \mu = 0$. Let $J \subset [n]$. If $J \in \sJ(\sF_e, \sF_o)$, then Proposition~\ref{prop - deltaj mu nu zero sphere} implies that $\delta_J \mcs \mu = 0$ since $\delta_J \mcs \nu \neq 0$ by the assumptions. If $J \notin \sJ(\sF_e, \sF_o)$, it follows from Corollary~\ref{coro - deltaj zero even odd sphere} that $\delta_J \mcs \mu = 0$. Now it follows from Corollary~\ref{coro - sign meas bn decomp signs}~\eqref{coro - sign meas bn decomp signs zero} that $\mu = 0$.

{\em Necessity.} Let $J \in \sJ(\sF_e, \sF_o)$. Then $P_S(\delta_J) \in \sM(B_{[n]}; \sF_e, \sF_o)$ by Remark~\ref{rema - j e properties}, and clearly $P_S(\delta_J) \neq 0$. If $\nu$ is spherically universal on $\sM(B_{[n]}; \sF_e, \sF_o)$, then
\[
\delta_J \mcs \nu \,=\, P_S(\delta_J) \mcs \nu \,\neq\, 0
\]

{\em Supplement.} From Proposition~\ref{prop - j gen group even odd}~\eqref{prop - j gen group even odd 1} we know that $\sJ(\sF_e, \sF_o) = \sJ(\sG_e, \sG_o)$. Now suppose that $\delta_J \mcs \nu \neq 0$ holds for every $J \in \sJ(\sG_e, \sG_o)$ and that $\nu$ is even (odd) under $T_E$ where $E \subset [n]$. By Theorem~\ref{th - deltaj zero arb e sym} we have $\delta_J \mcs \nu = 0$ for every $J \in \sJ(\O, E)$ ($J \in \sJ(E, \O)$). It follows that $\sJ(\sG_e, \sG_o) \subset \sJ(E, \O)$ ($\sJ(\sG_e, \sG_o) \subset \sJ(\O, E)$). Since $(\sG_e, \sG_o)$ is proper, Lemma~\ref{le - j s g impl} yields $E \in \sG_e$ ($E \in \sG_o$).
\end{proof}

\begin{proposition}
\label{prop - nu proj univ even odd sphere}
Let $\sF_e, \sF_o \subset \sP_n$ and $E, G \subset [n]$ with $\O \neq G \subset E$. Further let $\nu \in \sM(S^{n-1})$. Then $\nu$ is spherically universal on $\sM(B_G; \sF_e, \sF_o)$ if and only if $P^S_E(\nu)$ is spherically universal on $\sM(B_G; \sF_e, \sF_o)$.
\end{proposition}

\begin{proof}
This follows from the fact that $P^S_E(\nu) \mcs \mu = \nu \mcs \mu$ for every $\mu \in \sM(B_G; \sF_e, \sF_o)$.
\end{proof}

\begin{lemma}
\label{le - univ gen sym gen supp even odd sphere}
Let $\nu \in \sM(S^{n-1})$ and $\sF_e, \sF_o \subset \sP_n$ such that $(\sF_e, \sF_o)$ is a proper generating pair. If $\nu$ is spherically universal on $\sM(B_{[n]}; \sF_e, \sF_o)$, then also $\nu_{[n]}$ is spherically universal on $\sM(B_{[n]}; \sF_e, \sF_o)$.
\end{lemma}

\begin{proof}
Let $J \in \sJ(\sF_e, \sF_o)$. Then 
\begin{eqnarray*}
\nu \mcs \delta_J \mcs \sigma_0^n & = & (\nu \mcs \delta_J)_{[n]} \mcs \sigma_0^n \\
 & = & \nu_{[n]} \mcs \delta_J \mcs \sigma_0^n
\end{eqnarray*}
by Corollary~\ref{coro - sigma0 n term}. Lemma~\ref{le - deltaj even odd class}~\eqref{le - deltaj even odd class 2} says that
\[
0 \,\neq\, \delta_J \mcs \sigma_0^n \,\in\, \sM(B_{[n]}; \sF_e, \sF_o)
\]
Hence, if $\nu$ is spherically universal on $\sM(B_{[n]}; \sF_e, \sF_o)$, then $\nu_{[n]} \mcs \delta_J \mcs \sigma_0^n \neq 0$, and thus $\nu_{[n]} \mcs \delta_J \neq 0$. Consequently $\nu_{[n]}$ is spherically universal on $\sM(B_{[n]}; \sF_e, \sF_o)$ by Theorem~\ref{th - measures univ general even odd be}.
\end{proof}

\begin{theorem}
\label{th - gen even odd supp char sphere}
Let $\sE, \sF_e, \sF_o \subset \sP_n$ with $\O \notin \sE$, and $\nu \in \sM(S^{n-1})$. Then $\nu$ is spherically universal on $\sM^{\rm sph}(\sE; \sF_e, \sF_o)$ if and only if
\begin{equation}
\label{th - gen even odd supp char gen sphere 1}
\delta_J \mcs R_E P^S_E(\nu) \neq 0,\quad E \in \sE,\; J \in \sJ(E; \sF_e, \sF_o)
\end{equation}
In this case, if $( \sF_e |_E, \sF_o |_E)$ is proper and $R_E P^S_E(\nu)$ is even (odd) under $T_G$ for some $E \in \sE$ and $G \subset E$, then $G \in \sG_e |_E$ ($G \in \sG_o |_E$) where $\sG$ is the group generated by $\sF$.
\end{theorem}

\begin{remark}
In the case $[n] \in \sE$, if $(\sF_e, \sF_o)$ is proper and $\nu$ is universal on $\sM^{\rm sph}(\sE; \sF_e, \sF_o)$ and even (odd) under $T_G$ for some $G \subset [n]$, then also $\nu_{[n]}$ is even (odd) under $T_G$ and the second statement in Theorem~\ref{th - gen even odd supp char sphere} says that $G \in \sG_e$ ($G \in \sG_o$).
\end{remark}

Remark~\ref{rema - non proper conditions} applies to spherical universality as well.

\begin{proof}[Proof of Theorem~\ref{th - gen even odd supp char sphere}]
{\em Sufficiency.} Let $\mu \in \sM^{\rm sph}(\sE; \sF_e, \sF_o)$. Assume that $\nu \mcs \mu = 0$. If $[n] \notin \sE$, then $\mu_{[n]} = 0$. If $[n] \in \sE$, condition~\eqref{th - gen even odd supp char gen sphere 1} for $E = [n]$ says that $\delta_J \mcs \nu_{[n]} \neq 0$ for every $J \in \sJ(\sF_e, \sF_o)$. Since multiple reflections commute with $R_{[n]}$, we have $\mu_{[n]} \in \sM(B_{[n]}; \sF_e, \sF_o)$. Theorem~\ref{th - measures univ general even odd be} implies that $\nu_{[n]}$ is spherically universal on~$\sM(B_{[n]}; \sF_e, \sF_o)$. Now note that
\[
\nu_{[n]} \mcs \mu_{[n]} \,=\, (\nu \mcs \mu)_{[n]} \,=\, 0
\]
Therefore $\mu_{[n]} = 0$.

Now let $k \in \{1, 2, \ldots, n - 1\}$ and assume that $\mu_M = 0$ for all $M \subset [n]$ with $|M| > k$. Let $E \in \sE$ such that $|E| = k$. We have
\begin{eqnarray}
0 & = & (\nu \mcs \mu)_E \nonumber\\
 & = & (\nu \mcs \mu_E)_E \label{eq - even odd gen supp sphere 1}\\
 & = & \big( P^S_E(\nu) \mcs \mu_E \big)_E \label{eq - even odd gen supp sphere 2}\\
 & = & R_E P^S_E(\nu) \mcs \mu_E \label{eq - even odd gen supp sphere 3}
\end{eqnarray}
where \eqref{eq - even odd gen supp sphere 1} and \eqref{eq - even odd gen supp sphere 3} follow from Corollary~\ref{cor - decomp star prod}, and \eqref{eq - even odd gen supp sphere 2} from Corollary~\ref{coro - radial proj prod}. By the assumption the first factor in \eqref{eq - even odd gen supp sphere 3} satisfies
\[
\delta_{E, J} \mcs R_E P^S_E(\nu) \,=\, \delta_J \mcs R_E P^S_E(\nu) \,\neq\, 0,\quad J \in \sJ(E; \sF_e, \sF_o)
\]
Theorem~\ref{th - measures univ general even odd be} implies that $R_E P^S_E(\nu)$ is spherically universal on~$\sM(B_E; \sF_e, \sF_o)$. Moreover we know that $\mu_E \in \sM(B_E; \sF_e, \sF_o)$. Consequently $\mu_E = 0$.

{\em Necessity.} Assume that $\nu$ is spherically universal on $\sM^{\rm sph}(\sE; \sF_e, \sF_o)$ and let $E \in \sE$. By Proposition~\ref{prop - nu proj univ even odd sphere} $P^S_E(\nu)$ is spherically universal on $\sM(B_E; \sF_e, \sF_o)$. Lemma~\ref{le - univ gen sym gen supp even odd sphere} implies that also $R_E P^S_E(\nu)$ is spherically universal on $\sM(B_E; \sF_e, \sF_o)$. Now let $J \in \sJ(E; \sF_e, \sF_o)$. Noting that $\delta_{E, J} \in \sM^{\rm int}(A_E; \sF_e, \sF_o)$ and thus $P_S(\delta_{E, J}) \in \sM(B_E; \sF_e, \sF_o)$, we have
\[
\delta_J \mcs R_E P^S_E(\nu) \,=\, \delta_{E, J} \mcs R_E P^S_E(\nu) \,\neq\, 0
\]

{\em Supplement.} Let $E \in \sE$. Under the stated conditions $R_E P^S_E(\nu)$ is spherically universal on
\[
\sM(B_E; \sF_e, \sF_o) \,=\, \sM(B_E; \sF_e |_E, \sF_o |_E)
\]
Recall from Remark~\ref{rema - proj even odd} that $\sG_e |_E$ and $\sG_o |_E$ are the even and odd parts, respectively, of the group generated by $(\sF_e |_E, \sF_o |_E)$. Now it follows from the last statement of Theorem~\ref{th - measures univ general even odd be} that if $R_E P^S_E(\nu)$ is even (odd) under $T_G$ for some $G \subset E$, then $G \in \sG_e |_E$ ($G \in \sG_o |_E$).
\end{proof}

\begin{remark}
Notice a difference in the sufficiency part of the proof of Theorem~\ref{th - gen even odd supp char sphere} to that of Theorem~\ref{th - gen even odd supp char} since the induction over $|E|$ is only for $|E| \geq 1$ here.
\end{remark}

\begin{corollary}
\label{coro - measures gen even odd univ sphere char}
Let $\sF_e, \sF_o \subset \sP_n$ and $\nu \in \sM(S^{n-1})$. Then:
\begin{thnumber}
\item \label{coro - measures gen even odd univ sphere char bn} $\nu$ is spherically universal on $\sM(B_{[n]}; \sF_e, \sF_o)$ if and only if
\begin{equation*}
\label{coro - measures gen even odd univ rn char bn 1}
\delta_J \mcs \nu_{[n]} \neq 0,\quad J \in \sJ(\sF_e, \sF_o)
\end{equation*}
\item \label{coro - measures gen even odd univ sphere char gen} $\nu$ is spherically universal on $\sM(S^{n-1}; \sF_e, \sF_o)$ if and only if
\begin{equation}
\label{coro - measures gen even odd univ sphere char gen 1}
\delta_J \mcs R_E P^S_E(\nu) \neq 0,\quad E \subset [n],\; E \neq \O,\; J \in \sJ(E; \sF_e, \sF_o)
\end{equation}
\item \label{coro - measures gen even odd univ sphere char order-n 1} For $\nu$ of order $[n]$, $\nu$ is spherically universal on $\sM(S^{n-1}; \sF_e, \sF_o)$ if and only if
\begin{equation}
\label{coro - measures gen even odd univ sphere char gen 2}
\delta_{E, J} \mcs \nu \neq 0,\quad E \subset [n],\; E \neq \O,\; J \in \sJ(E; \sF_e, \sF_o)
\end{equation}
\item \label{coro - measures gen even odd univ sphere char order-n 2} For $\nu$ of order $[n]$, $\nu$ is spherically universal on $\sM(S^{n-1}; \sF_e, \sF_o)$ if and only if
both
\begin{equation}
\label{coro - measures gen even odd univ sphere char gen 3}
\delta_{J, J} \mcs \nu \neq 0,\quad J \in \sJ(\sF_e, \sF_o),\; J \neq \O
\end{equation}
and, provided $\sF_o = \O$, also
\begin{equation}
\label{coro - measures gen even odd univ sphere char gen 4}
M_{\rm sym} \, P^S_i(\nu) \neq 0,\quad i \in [n]
\end{equation}
\end{thnumber}
\end{corollary}

\begin{proof}
Statements \eqref{coro - measures gen even odd univ sphere char bn} and \eqref{coro - measures gen even odd univ sphere char gen} follow from Theorem~\ref{th - gen even odd supp char sphere} for the special cases $\sE = \{ [n] \}$ and $\sE = \sP_n \setminus \{\O\}$.

To see \eqref{coro - measures gen even odd univ sphere char order-n 1} assume that $\nu$ is of order $[n]$. In this case for $J \subset E \subset [n]$ with $E \neq \O$ we have
\[
\delta_{E, J} \mcs R_E P^S_E(\nu) \,=\, \delta_{E, J} \mcs P^S_E(\nu) \,=\, P_E(\delta_{E, J}) \mcs \nu \,=\, \delta_{E, J} \mcs \nu
\]

Finally we prove statement~\eqref{coro - measures gen even odd univ sphere char order-n 2}. To this end we show that \eqref{coro - measures gen even odd univ sphere char gen 2} for $J \neq \O$ is equivalent with \eqref{coro - measures gen even odd univ sphere char gen 3} and that \eqref{coro - measures gen even odd univ sphere char gen 2} for $J = \O$ is equivalent with \eqref{coro - measures gen even odd univ sphere char gen 4}. In the first case note that \eqref{coro - measures gen even odd univ sphere char gen 2} for $J \neq \O$ clearly implies \eqref{coro - measures gen even odd univ sphere char gen 3} as the latter is a special case of the former. To see the converse note that, for $J \subset E \subset [n]$ with $J \neq \O$,
\[
P^S_J(\delta_{E, J} \mcs \nu) \,=\, P_J(\delta_{E, J}) \mcs \nu \,=\, \delta_{J, J} \mcs \nu
\]
Hence $\delta_{E, J} \mcs \nu = 0$ implies $\delta_{J, J} \mcs \nu = 0$. Recall that the second case (i.e.\ $J = \O$) occurs if and only if $\sF_o = \O$, see Remark~\ref{rema - j e properties}. Now note that for $E = \{i\}$ with $i \in [n]$ we have
\[
\delta_{\{i\}, \O} \mcs  \nu \,=\, \delta_{\O} \mcs  P^S_i(\nu) \,=\, M_{\rm sym} \, P^S_i(\nu)
\]
Hence \eqref{coro - measures gen even odd univ sphere char gen 2} for $J = \O$ implies \eqref{coro - measures gen even odd univ sphere char gen 4}. Conversely note that, for $i \in E \subset [n]$,
\[
P^S_i (\delta_{E, \O} \mcs \nu) \,=\, P_i(\delta_{E, \O}) \mcs \nu \,=\, \delta_{\{i\}, \O} \mcs \nu
\]
Thus $\delta_{E, \O} \mcs \nu = 0$ implies that $\delta_{\{i\}, \O} \mcs \nu = 0$.
\end{proof}

\begin{example}
\label{ex - deltaej universal sphere}
Similar considerations hold for spherical universality as stated in Example~\ref{ex - deltaej universal} for universality. So let
$\nu \in \sM(S^{n-1})$ and $(\sG_e, \sG_o)$ be a proper symmetry pair. It follows from Corollary~\ref{coro - measures gen even odd univ sphere char}~\eqref{coro - measures gen even odd univ sphere char bn} that $\nu$ can only be spherically universal on $\sM(B_{[n]}; \sG_e, \sG_o)$ if it has degree~$n$. Now let $J \subset E \subset [n]$. If $P_S(\delta_{E, J})$ is spherically universal on $\sM(B_{[n]}; \sG_e, \O)$, then necessarily $E = [n]$, $J = \O$, and $\sG_e = \sP_n$.
\end{example}

\subsection{Specific symmetries}

Again we first consider unconditionality corresponding to the symmetry pair $(\sG_e, \sG_o) = (\sP_n, \O)$, and thus $\sJ(\sG_e, \O) = \{ \O \}$.

\begin{theorem}[Unconditional measures]
\label{th - measures unc univ sphere char}
Let $\nu \in \sM(S^{n-1})$. Then:
\begin{thnumber}
\item \label{th - measures unc univ sphere char ae} $\nu$ is spherically universal on $\sM_{\rm unc}(B_{[n]})$ if and only if $M_{\rm unc} (\nu)$ has degree $n$.
\item \label{th - measures unc univ sphere char gen} $\nu$ is spherically universal on $\sM_{\rm unc}(S^{n-1})$ if and only if $R_E P^S_E M_{\rm unc}(\nu) \neq 0$ for every $E \subset [n]$ with $E \neq \O$.
\item \label{th - measures unc univ sphere char order-n} For $\nu$ of order $[n]$, $\nu$ is spherically universal on $\sM_{\rm unc}(S^{n-1})$ if and only if
\begin{equation}
\label{eq - uncond sphere char dim1}
M_{\rm sym} \, P^S_i(\nu) \neq 0,\quad i \in [n]
\end{equation}
\item \label{th - measures unc univ sphere char pos} For non-negative $\nu$, $\nu$ is spherically universal on $\sM_{\rm unc}(S^{n-1})$ if and only if it has degree~$n$.
\end{thnumber}
\end{theorem}

\begin{proof}
Statements \eqref{th - measures unc univ sphere char ae}, \eqref{th - measures unc univ sphere char gen}, and \eqref{th - measures unc univ sphere char order-n} follow from statements \eqref{coro - measures gen even odd univ sphere char bn}, \eqref{coro - measures gen even odd univ sphere char gen}, and \eqref{coro - measures gen even odd univ sphere char order-n 2} in Corollary~\ref{coro - measures gen even odd univ sphere char}, respectively. Now \eqref{th - measures unc univ sphere char gen} implies \eqref{th - measures unc univ sphere char pos}.
\end{proof}

Observe that, if $\nu$ itself is even (e.g.\ unconditional) in Theorem~\ref{th - measures unc univ char}, then condition~\eqref{eq - uncond sphere char dim1} simplifies to $P^S_i(\nu) \neq 0$ for every $i \in [n]$.

Next we consider the origin-symmetric case, given by $\sG_e = \{ \O, [n] \}$ and $\sG_o = \O$, and therefore
\[
\sJ(\sG_e, \sG_o) = \{ J \subset [n] \,:\, |J| \;\, {\rm is \; even} \,\}
\]

\begin{theorem}[Symmetric measures]
\label{th - measures sym univ sphere char}
Let $\nu \in \sM(S^{n-1})$. Then:
\begin{thnumber}
\item \label{th - measures sym univ sphere char bn} $\nu$ is spherically universal on $\sM_{\rm sym}(B_{[n]})$ if and only if
\[
\delta_J \mcs \nu_{[n]} \neq 0,\quad J \subset [n],\; |J| \; {\it even}
\]
\item \label{th - measures sym univ sphere char gen} $\nu$ is spherically universal on $\sM_{\rm sym}(S^{n-1})$ if and only if
\[
\delta_J \mcs R_E P^S_E(\nu) \neq 0,\quad J \subset E \subset [n],\; E \neq \O,\; |J| \; {\it even}
\]
\item \label{th - measures sym univ sphere char order-n} For $\nu$ of order $[n]$, $\nu$ is spherically universal on $\sM_{\rm sym}(S^{n-1})$ if and only if
both
\[
\delta_{J, J} \mcs \nu \neq 0,\quad J \subset [n],\; J \neq \O,\; |J| \; {\it even}
\]
and
\begin{equation}
\label{eq - sym sphere char dim1}
M_{\rm sym} \, P^S_i(\nu) \neq 0,\quad i \in [n]
\end{equation}
\end{thnumber}
\end{theorem}

\begin{proof}
Statements \eqref{th - measures sym univ sphere char bn}, \eqref{th - measures sym univ sphere char gen}, and \eqref{th - measures sym univ sphere char order-n} follow from statements \eqref{coro - measures gen even odd univ sphere char bn}, \eqref{coro - measures gen even odd univ sphere char gen}, and \eqref{coro - measures gen even odd univ sphere char order-n 2} in Corollary~\ref{coro - measures gen even odd univ sphere char}, respectively.
\end{proof}

If $\nu$ itself is origin-symmetric, then condition~\eqref{eq - sym sphere char dim1} simplifies to $P^S_i(\nu) \neq 0$ for every $i \in [n]$ as in the case of unconditionality.

The case of no reflection symmetry, that is $\sF_e = \sF_o = \O$ and $\sJ(\O, \O) = \sP_n$, is stated in the following.

\begin{theorem}[No reflection symmetry]
\label{th - measures general univ sphere char}
Let $\nu \in \sM(S^{n-1})$.
\begin{thnumber}
\item \label{th - measures general univ sphere char be} $\nu$ is spherically universal on $\sM(B_{[n]})$ if and only if
\[
\delta_J \mcs \nu_{[n]} \neq 0,\quad J \subset [n]
\]
\item \label{th - measures general univ sphere char gen} $\nu$ is spherically universal on $\sM(S^{n-1})$ if and only if
\[
\delta_J \mcs R_E P^S_E(\nu) \neq 0,\quad E \subset [n],\; E \neq \O,\; J \subset E
\]
\item \label{th - measures general univ sphere char order-n} For $\nu$ of order $[n]$, $\nu$ is spherically universal on $\sM(S^{n-1})$ if and only if
both
\[
\delta_{J, J} \mcs \nu \neq 0,\quad J \subset [n],\; J \neq \O,
\]
and
\[
M_{\rm sym} \, P^S_i(\nu) \neq 0,\quad i \in [n]
\]
\end{thnumber}
\end{theorem}

\begin{proof}
Statements \eqref{th - measures general univ sphere char be}, \eqref{th - measures general univ sphere char gen}, and \eqref{th - measures general univ sphere char order-n} follow from statements \eqref{coro - measures gen even odd univ sphere char bn}, \eqref{coro - measures gen even odd univ sphere char gen}, and \eqref{coro - measures gen even odd univ sphere char order-n 2} in Corollary~\ref{coro - measures gen even odd univ sphere char}, respectively.
\end{proof}

\begin{example}
\label{ex - sigma0n univ sphere}
This is very similar to Example~\ref{ex - sigma0n univ}. Since $\delta_J \mcs \sigma_0^n \neq 0$ holds for every $J \subset [n]$ by Lemma~\ref{le - deltaj even odd class}~\eqref{le - deltaj even odd class 2}, the condition of Theorem~\ref{th - measures general univ sphere char}~\eqref{th - measures general univ sphere char be} is satisfied for $\nu = P_S(\sigma_0^n)$. Therefore $P_S(\sigma_0^n)$ is spherically universal on $\sM(B_{[n]})$.
\end{example}

Finally we consider the anti-symmetric case, that is $\sG_e = \{ \O \}$ and $\sG_o = \{ [n] \}$. In this case
\[
\sJ(\sG_e, \sG_o) = \{ J \subset [n] \,:\, |J| \;\, {\rm is \; odd} \,\}
\]

\begin{theorem}[Anti-symmetric measures]
\label{th - measures asym univ sphere char}
Let $\nu \in \sM(S^{n-1})$. Then:
\begin{thnumber}
\item \label{th - measures asym univ sphere char ae} $\nu$ is spherically universal on $\sM_{\rm asym}(B_{[n]})$ if and only if
\[
\delta_J \mcs \nu_{[n]} \neq 0,\quad J \subset [n],\; |J| \; {\it odd}
\]
\item \label{th - measures asym univ sphere char gen} $\nu$ is spherically universal on $\sM_{\rm asym}(S^{n-1})$ if and only if
\[
\delta_J \mcs R_E P^S_E(\nu) \neq 0,\quad E \subset [n],\; E \neq \O,\; J \subset E,\; |J| \; {\it odd}
\]
\item \label{th - measures asym univ sphere char order-n} For $\nu$ of order $[n]$, $\nu$ is spherically universal on $\sM_{\rm asym}(S^{n-1})$ if and only if
\[
\delta_{J, J} \mcs \nu \neq 0,\quad J \subset [n],\; |J| \; {\it odd}
\]
\end{thnumber}
\end{theorem}

\begin{proof}
Statements \eqref{th - measures asym univ sphere char ae}, \eqref{th - measures asym univ sphere char gen}, and \eqref{th - measures asym univ sphere char order-n} follow from statements \eqref{coro - measures gen even odd univ sphere char bn}, \eqref{coro - measures gen even odd univ sphere char gen}, and \eqref{coro - measures gen even odd univ sphere char order-n 2} in Corollary~\ref{coro - measures gen even odd univ sphere char}.
\end{proof}

\section{Application to convex geometry}
\label{sec - application}

The results in the present paper yield a generalisation of some theorems in~\cite{mol:nag:21}. In~\cite[Section~3]{mol:nag:21} the properties of generalised zonoids under diagonal transformations are analysed, in particular the so-called D-universality. First it needs to be clarified how D-universality of a convex body as defined in~\cite{mol:nag:21} is related to spherical universality of a measure on the sphere in the present context.

A convex body $K$ (that is a non-empty, convex, compact subset of~$\R^n$) can analytically be specified by its support function, see for example \cite[Section~1.7.1]{schn2}:
\begin{equation}
\label{eq - def support function}
h(K,x) \,=\, \sup \{ \langle x, y \rangle : y \in K \},\quad x \in \R^n
\end{equation}
A support function is any sublinear function from $\R^n$ to $\R$; in particular it is positively one-homogeneous and therefore determined by its values on~$S^{n-1}$. $K$ is called {\em generalised zonoid}, see for example~\cite[p.~195]{schn2}, if there is $\nu \in \sM(S^{n-1})$ such that 
\begin{equation}
\label{eq - def generating measure}
h(K,u) = \int_{S^{n-1}} \nu(dv) \, |\langle v, u \rangle|,\quad u \in S^{n-1}
\end{equation}
Such $K$ is necessarily origin-symmetric.\footnote{Sometimes translates of such sets are also called generalised zonoids.} $\nu$ is called {\em generating measure} of $K$. $\nu$ can chosen to be symmetric, and then it is unique for given $K$. Note that $\nu = 0$ corresponds to $K = \{ 0 \}$; in this case the generalised zonoid and its generating measure are called {\em trivial}. Not every measure on the sphere yields a generalised zonoid. A symmetric measure $\nu$ on the sphere is a generating measure (for a generalised zonoid) if and only if all of its two-dimensional orthogonal projections are non-negative \cite[Korollar~4.3]{weil82}. Non-negative measures represent zonoids, which are the limits of zonotopes (finite sums of segments) with respect to the Hausdorff metric. For $n = 2$ all origin-symmetric convex bodies are zonoids. For $n \geq 3$ the family of zonoids is closed and nowhere dense with respect to the Hausdorff metric in the family of origin-symmetric convex bodies; however generalised zonoids are dense but not closed \cite[Corollary~3.5.7]{schn2}. Even if $\nu$ is not a generating measure, that is the right-hand side of \eqref{eq - def generating measure} is not a sublinear function in~$u$, that function uniquely specifies $\nu$ provided it is origin-symmetric, see \cite[Theorem~3.5.4]{schn2}.

Now let $K$ be a generalised zonoid in~$\R^n$ with generating measure~$\nu$. For $E \subset [n]$ its orthogonal projection $P_E(K)$ on $H_E$ is generated by the measure $P^S_E(\nu)$ concentrated on $S_E$; this follows from \eqref{eq - def support function}, \eqref{eq - def generating measure}, and Proposition~\ref{prop - multconv poshom}. It is also shown in~\cite[Satz~4.1]{weil82}. Furthermore it follows from \cite[Lemma~3.5.6]{schn2} that the support sets of $K$ in all the coordinate directions are singletons (so-called exposed points) if and only if $M_{\rm sym}(\nu)$ is of order~$[n]$ or $M_{\rm sym}(\nu)$ is trivial.

We now turn more specifically to diagonal transformations of convex bodies and D-universality. Given an origin-symmetric convex body $K$ in $\R^n$ and $\mu \in \sM(S^{n-1})$ the $K$-transform of $\mu$ is defined in~\cite{mol:nag:21} as
\[
T_K \mu (u) \,=\, \int_{S^{n-1}} \mu(dv) \, h(vK, u),\quad u \in S^{n-1},
\]
where the notation \eqref{eq - comp prod vec set} is used. Equivalently a diagonal matrix can be applied to all elements of~$K$, motivating the nomenclature. $K$ is called {\em D-universal} ({\em unconditionally D-universal}) in \cite{mol:nag:21} if $T_K \mu = 0$ implies $\mu = 0$ for every $\mu \in \sM_{\rm sym}(S^{n-1})$ ($\mu \in \sM_{\rm unc}(S^{n-1})$). The following lemma clarifies how these two notions are related to spherical universality:

\begin{lemma}
\label{le - equiv universality gen zon}
Let $K$ be a generalised zonoid in~$\R^n$ with generating measure~$\nu$. Then:
\begin{thnumber}
\item \label{le - equiv universality gen zon sym} $K$ is D-universal if and only if $\nu$ is spherically universal on $\sM_{\rm sym}(S^{n-1})$.
\item \label{le - equiv universality gen zon unc} $K$ is unconditionally D-universal if and only if $\nu$ is spherically universal on $\sM_{\rm unc}(S^{n-1})$.
\end{thnumber}
\end{lemma}

\begin{proof}
Note that for $\mu \in \sM_{\rm sym}(S^{n-1})$ we have
\[
T_K \mu (u) \,=\, \int_{S^{n-1}} \nu \mcs \mu (dv) \, |\langle v, u \rangle|,\quad u \in S^{n-1}
\]
Now the two assertions follow from \cite[Theorem~3.5.4]{schn2}.
\end{proof}

A closer look is now taken at both cases.

\subsection{Unconditionally D-universal convex bodies}

In \cite{mol:nag:21} the unconditional case is considered in Theorem~3.2 and its Corollaries 3.5 and~3.6. In \cite[Theorem~3.2]{mol:nag:21} it is assumed that all the support sets of the generalised zonoid $K$ in the coordinate directions are singletons. It is tacitly assumed that $K$ is not trivial. Hence its generating measure $\nu$ is of order~$[n]$. \cite[Theorem~3.2]{mol:nag:21} can be derived from Theorem~\ref{th - measures unc univ sphere char}~\eqref{th - measures unc univ sphere char order-n} as follows.

\begin{corollary}
\label{coro - univ geom}
Let $K$ be a non-trivial generalised zonoid in~$\R^n$ such that all support sets in the coordinate directions are singletons. Then $K$ is unconditionally D-universal.
\end{corollary}

\begin{proof}
Let $\nu$ be the generating measure of $K$. It follows from the assumptions that $\nu$ is of order $[n]$. By Lemma~\ref{le - equiv universality gen zon}~\eqref{le - equiv universality gen zon unc} it is enough to prove that $\nu$ is spherically universal on $\sM_{\rm unc}(S^{n-1})$. We show that \eqref{eq - uncond sphere char dim1} holds. Assume $P^S_i M_{\rm sym}(\nu)$ is trivial for some $i \in [n]$. Then $P_i(K)$ is trivial. This implies $K \subset H_{[n] \setminus \{i\}}$. It follows from~\cite[Lemma~3.5.6]{schn2} that $M_{\rm sym}(\nu)$ is concentrated on $S_{[n] \setminus \{i\}}$, contradicting the fact that $M_{\rm sym}(\nu)$ is of order~$[n]$.
\end{proof}

Theorem~\ref{th - measures unc univ sphere char}~\eqref{th - measures unc univ sphere char pos} immediately yields that a zonoid is unconditionally universal if its generating (non-negative) measure has degree~$n$. This is a generalisation to \cite[Corollary~3.5]{mol:nag:21} where this could only be shown for $n = 2$. Further generalisations are contained in \eqref{th - measures unc univ sphere char ae} and \eqref{th - measures unc univ sphere char gen} of Theorem~\ref{th - measures unc univ sphere char}, each stating sufficient and necessary conditions on the generating measure. Notice that in both cases $M_{\rm unc} (\nu)$ is supposed to have degree $n$ but may contain terms of orders other than $[n]$, which says in geometric terms that the support sets in the coordinate directions are not necessarily singletons.

\cite[Corollary~3.6]{mol:nag:21} is a statement about universality on $\sM^{\rm int}(\R^n_+)$. There it is required that $\nu$ is a non-negative integrable measure on~$(0,\infty)^n$. Now Corollary~\ref{coro - univ rn plus} imposes weaker conditions on $\nu$, viz.\ $\nu \in \sM^{\rm int}(\R^n)$ and either (1) $\nu$ is of order $[n]$ and $\nu(\R^n) \neq 0$, or (2) $\nu$ has degree $n$ and $\nu \geq 0$. These correspond to conditions \eqref{th - measures unc univ char order-n} and \eqref{th - measures unc univ char pos} of Theorem~\ref{th - measures unc univ char}.

\subsection{D-universal convex bodies}

In \cite{mol:nag:21} the origin-symmetric case is stated in Theorem~3.7. Since generalised zonoids are always symmetric, this is the most general case in that context. \cite[Theorem~3.7]{mol:nag:21} again requires the support sets of the generalised zonoid $K$ in the coordinate directions to be singletons and therefore requires its generating measure $\nu$ to be of order~$[n]$. This situation corresponds to Theorem~\ref{th - measures sym univ sphere char}~\eqref{th - measures sym univ sphere char order-n}. Condition~\eqref{eq - sym sphere char dim1} is automatically satisfied as in the unconditional case by the argument in the proof of Corollary~\ref{coro - univ geom}. The remaining conditions in Theorem~\ref{th - measures sym univ sphere char}~\eqref{th - measures sym univ sphere char order-n} are equivalent to the "asymmetry conditions" (AS) in~\cite{mol:nag:21}. Note that the conditions in \cite[Theorem~3.7]{mol:nag:21} are also necessary, yielding to a full characterisaton in the case in which all support sets in coordinate directions are singletons.

Beyond the results in \cite{mol:nag:21} statements~\eqref{th - measures sym univ sphere char bn} and \eqref{th - measures sym univ sphere char gen} in Theorem~\ref{th - measures sym univ sphere char} each lists a set of sufficient and necessary conditions such that $K$ is universal on the respective family of measures. Notice in particular that $\nu$ may contain terms of order different from~$[n]$, again loosening the singleton support set condition of~\cite{mol:nag:21}.

\subsection{Method of decomposition}

Finally notice that the decomposition of random vectors used in \cite[Lemma~3.4]{mol:nag:21} is the one mentioned in Remark~\ref{rema - coord decomposition random vectors}. This is different from the coordinate decomposition~\eqref{eq - sign meas decomp} used in the present work, and it is not sufficient in order to obtain the more general results obtained here.

\section{Measures of degree \texorpdfstring{$[n]$}{[n]}}
\label{sec - measures of degree n}

The results of this section are not used for the proofs in the rest of the paper. Rather Lemma~\ref{le - combinat id pse}. provides a simple criterion for a measure to have degree~$n$. On $\sM(\R^n)$ we define the linear operator
\[
\widehat{P} = \sum_{E \subset [n]} (-1)^{|E|} \, P_E
\]

\begin{lemma}
\label{le - combinat id pse}
Let $\mu \in \sM(\R^n)$. Then $\mu_{[n]} = 0$ if and only if
\begin{equation}
\label{eq - sum pse id}
\widehat{P}(\mu) \,=\, 0
\end{equation}
\end{lemma}

\begin{remark}
\label{re - combinat id pse}
Note that the term for $E = [n]$ on the left-hand side of~\eqref{eq - sum pse id} yields $(-1)^n \mu$. Thus identity \eqref{eq - sum pse id} allows us to express $\mu$ in terms of its proper orthogonal projections:
\[
\sum_{E \subsetneqq [n]} (-1)^{|E|} \, P_E(\mu) \,=\, (-1)^{n + 1} \mu
\]
\end{remark}

\begin{proof}[Proof of Lemma~\ref{le - combinat id pse}]
{\em Necessity.} Suppose $\mu_{[n]} = 0$. By linearity of the orthogonal projections we may assume that $\mu$ is of order $F \subset [n]$ with $|F| \leq n-1$. We have
\begin{eqnarray}
\widehat{P}(\mu) & = & \sum_{G \subset F} \sum \big\{ (-1)^{|E|} \, P_E(\mu) \,:\, E \subset [n],\; E \cap F = G \big\} \nonumber\\
 & = & \sum_{G \subset F} \left( P_G(\mu) \; \sum \big\{ (-1)^{|E|}  \,:\, E \subset [n],\; E \cap F = G \big\} \right) \label{eq - pse measure id 1} \\
 & = &  \left( \sum_{G \subset F} (-1)^{|G|} \; P_G(\mu) \right) \left( \sum_{k = 0}^{n - |F|} (-1)^k {n - |F| \choose k} \right) \label{eq - pse measure id 2}
\end{eqnarray}
where in \eqref{eq - pse measure id 1} we have applied Proposition~\ref{prop - proj decomp measure}~\eqref{prop - proj decomp measure simpl}. Now notice that the second factor in \eqref{eq - pse measure id 2} is zero.

{\em Sufficiency.} Assume that $\mu$ has degree~$n$. Then $P_{[n]}(\mu)$ has degree $n$. Moreover $P_E(\mu)$ has degree less than $n$ for every $E \subsetneqq [n]$. It follows that the left-hand side of~\eqref{eq - sum pse id} is different from zero.
\end{proof}

\begin{lemma}
\label{le - combinat id pse sphere}
Let $\mu \in \sM(S^{n-1})$. Then $\mu_{[n]} = 0$ if and only if $P_S \widehat{P}(\mu) = 0$.
\end{lemma}

\begin{proof}
Necessity follows from Lemma~\ref{le - combinat id pse}. Sufficiency can be shown in a similar way as in the proof of Lemma~\ref{le - combinat id pse}.
\end{proof}

Finally we establish connections to Corollary~\ref{coro - sigma0 n term}. So let $\mu \in \sM(\R^n)$. It follows from Lemma~\ref{le - combinat id pse} that
\begin{equation}
\label{eq - proj highest order}
\widehat{P}(\mu)  = \widehat{P}(\mu_{[n]})
\end{equation}
Further we know from Corollary~\ref{coro - sigma0 n term}~\eqref{coro - sigma0 n term rn} that
\begin{equation}
\label{eq - sigma0 highest order}
\mu \mc \sigma_0^n = \mu_{[n]} \mc \sigma_0^n
\end{equation}
We show that \eqref{eq - proj highest order} implies \eqref{eq - sigma0 highest order}. To this end first notice that $\widehat{P}(\sigma_0^n) = (-1)^n \sigma_0^n$ by Lemma~\ref{le - proj sigma0}. Assuming $\eqref{eq - proj highest order}$ we have
\begin{eqnarray*}
\mu \mc \sigma_0^n & = & (-1)^n \, \mu \mc \widehat{P}(\sigma_0^n) \\
 & = & (-1)^n \, \widehat{P}(\mu) \mc \sigma_0^n \\
 & = & (-1)^n \, \widehat{P}(\mu_{[n]}) \mc \sigma_0^n \\
 & = & (-1)^n \, \mu_{[n]} \mc \widehat{P}(\sigma_0^n) \\
 & = & \mu_{[n]} \mc \sigma_0^n 
\end{eqnarray*}

\section{Universality and lifting}
\label{sec - universality lifting}

This section is dedicated to the relation between universality and spherical universality. Recall the definitions of the linear map $L$, and the families of sets $\sE_L$ and $\sF^0$ in Section~\ref{sec - lifting}.

\begin{theorem}
\label{th - equiv univ spher univ}
Let $\sE, \sF_e, \sF_o \subset \sP_n$ and $\mu \in \sM^{\rm int}(\R^n)$. Then $\mu$ is universal on $\sM^{\rm int}(\sE; \sF_e, \sF_o)$ if and only if $L(\mu)$ is spherically universal on $\sM^{\rm sph}(\sE_L; \sF_e^0, \sF_o)$.
\end{theorem}

\begin{proof}
{\em Sufficiency.} Suppose $L(\mu)$ is spherically universal on $\sM^{\rm sph}(\sE_L; \sF_e^0, \sF_o)$. Let $\rho \in \sM^{\rm int}(\sE; \sF_e, \sF_o)$ such that $\mu \mc \rho = 0$. Clearly $L(\rho) \in \sM^{\rm sph}(\sE_L; \sF_e^0, \sF_o)$ by Lemma~\ref{le - bij gen supp gen even odd}. Statements \eqref{prop - lifting zero} and \eqref{prop - lifting conv} of Proposition \ref{prop - lifting} yield
\begin{equation}
\label{eq - prod hom zero even odd}
0 \,=\, L(\mu \mc \rho) \,=\, L(\mu) \mcs L(\rho)
\end{equation}
Hence $L(\rho) = 0$. Therefore $\rho = 0$ again by Proposition~\ref{prop - lifting}~\eqref{prop - lifting zero}.

{\em Necessity.} Assume that $\mu$ is universal on $\sM^{\rm int}(\sE; \sF_e, \sF_o)$. Let $\sigma \in \sM^{\rm sph}(\sE_L; \sF_e^0, \sF_o)$ such that $L(\mu) \mcs \sigma = 0$. By Lemma~\ref{le - bij gen supp gen even odd} it is possible to set $\rho = L^{-1}(\sigma)$. Then \eqref{eq - prod hom zero even odd} holds, and Proposition~\ref{prop - lifting}~\eqref{prop - lifting zero} implies $\mu \mc \rho = 0$. Since $\rho \in \sM^{\rm int}(\sE; \sF_e, \sF_o)$, it follows from the assertion that $\rho = 0$.
\end{proof}

\begin{remark}
In Theorem~\ref{th - equiv univ spher univ} an important class of special cases is that for which $\sE = \sP_n$, viz.\ $\mu$ is universal on $\sM^{\rm int}(\R^n; \sF_e, \sF_o)$ if and only if $L(\mu)$ is spherically universal on $\sM(S^n_0; \sF_e^0, \sF_o)$. In particular, $\mu$ is universal on $\sM^{\rm int}(\R^n)$ if and only if $L(\mu)$ is spherically universal on $\sM_{\rm sym}(S^n_0)$. And $\mu$ is universal on $\sM^{\rm int}_{\rm unc}(\R^n)$ if and only if $L(\mu)$ is spherically universal on $\sM_{\rm unc}(S^n_0)$.
\end{remark}

With the help of Theorem~\ref{th - equiv univ spher univ}, Theorem~\ref{th - gen even odd supp char} can be derived as a corollary of Theorem~\ref{th - gen even odd supp char sphere}.

\begin{proof}[Proof of Theorem~\ref{th - gen even odd supp char}]
Let $E \in \sE$ and $J \in \sJ(E_L; \sF_e^0, \sF_o)$. First note that
\begin{eqnarray}
R_{E_L} P^S_{E_L} (L(\nu)) & = & M_{\rm sym} P_S R_{E_L} P_{E_L} (\delta_1 \otimes \nu) \nonumber\\
 & = & M_{\rm sym} P_S \big( \delta_1 \otimes R_E P_E(\nu) \big) \nonumber\\
 & = & L(R_E P_E(\nu)) \label{eq - r p lifting}
\end{eqnarray}
Further
\begin{eqnarray}
\delta_{[n]_L, J} \mcs L(R_E P_E(\nu)) & = & P_S \big( \delta_{[n]_L, J} \mc M_{\rm sym} P_S \big( \delta_1 \otimes R_E P_E(\nu) \big) \big) \nonumber\\
 & = & P_S \big( \delta_{[n]_L, J} \mc \big( \delta_1 \otimes R_E P_E(\nu) \big) \big) \label{th - gen even odd supp char 1} \\
 & = & P_S \sum_{s_0 \in \{-1,1\}} \sum_{s \in \sS^n} \sigma_J(s_0, s) \; \delta_{(s_0, s)} \mc \big( \delta_1 \otimes R_E P_E(\nu) \big)  \nonumber\\
 & = & P_S \sum_{s_0 \in \{-1,1\}} \sum_{s \in \sS^n} \sigma_J(s_0, s) \, \big( \delta_{s_0} \otimes (\delta_s \mc R_E P_E(\nu)) \big) \label{th - gen even odd supp char 3} \\
 & = & P_S \big( \rho \otimes \big( \delta_{[n], J \setminus \{0\}} \mc R_E P_E(\nu) \big) \big) \nonumber
\end{eqnarray}
where
\begin{equation*}
\rho = \left\{ \begin{array}{ll}
\displaystyle \sum_{a \in \{-1,1\}} \sign(a) \, \delta_a,& 0 \in J, \\[1.5em]
\displaystyle \sum_{a \in \{-1,1\}} \delta_a,& 0 \notin J.
\end{array} \right.
\end{equation*}
In particular \eqref{th - gen even odd supp char 1} follows from Remark~\ref{re - e j sym} as for $J \in \sJ(E_L; \sF_e^0, \sF_o)$ $|J|$ is always even since $[n]_L \in \sF_e^0$. Moreover \eqref{th - gen even odd supp char 3} is a consequence of Proposition~\ref{prop - multconv prod}. It follows from Lemma~\ref{lemma - finite meas bij} and from the particular form of $\rho$ that
\[
\delta_{[n]_L, J} \mcs R_{E_L} P^S_{E_L} (L(\nu)) \,\neq\, 0
\]
holds if and only if
\[
\delta_{[n], J \setminus \{0\}} \mc R_E P_E(\nu) \,\neq\, 0
\]
Sufficiency and necessity now follow from Theorems \ref{th - equiv univ spher univ} and \ref{th - gen even odd supp char sphere}.

To show the second assertion let $E \in \sE$, $G \subset E$, and assume that $( \sF_e |_E, \sF_o |_E)$ is proper, that $\nu$ is universal on $\sM^{\rm int}(\sE; \sF_e, \sF_o)$, and that $R_E P_E(\nu)$ is even (odd) under~$T_G$. It follows from Lemma~\ref{le - proper lifting} that $(\sF_e^0 |_{E_L}, \sF_o |_{E_L})$ is proper. Moreover Theorem~\ref{th - equiv univ spher univ} implies that $L(\nu)$ is spherically universal on $\sM^{\rm sph}(\sE_L; \sF_e^0, \sF_o)$. Further $L(R_E P_E(\nu))$ is even (odd) under~$T_G$ by Remark~\ref{rema - L commutes refl}. It follows from Theorem~\ref{th - gen even odd supp char sphere}, equation~\eqref{eq - r p lifting} and the fact that $E_L \in \sE_L$ and $G \subset E_L$ that $G$ is in the even (odd) part of \[
\rho_{E_L} \circ \gamma \,(\sF_e^0, \sF_o) \,=\, \rho_{E_L} \circ \sL \,(\sG_e, \sG_o) 
\]
Hence, in the "even" case,
\[
G \,\in\, \rho_{E_L} \circ \sL \,(\sG_e) \,=\, \sL_E \circ \rho_E \,(\sG_e) 
\]
Since $0 \notin G$, we have $G \in \rho_E (\sG_e)$. Similarly, in the "odd" case, we obtain $G \in \rho_E (\sG_o)$.
\end{proof}

\begin{remark}
The preceding proof shows that the characterisation of universality in $\R^n$ as stated in Theorem~\ref{th - gen even odd supp char} can be considered as a special case of the universality on the sphere in $n+1$ dimensions.
\end{remark}


\begin{thebibliography}{10}

\bibitem{cohn}
D.~L.~Cohn.
\newblock {\em Measure theory}.
\newblock Springer, New York, 2nd edition, 2013.

\bibitem{dam:mik:ros:14}
E.~Damek, T.~Mikosch, J.~Rosi\'{n}ski, and G.~Samorodnitsky.
\newblock General inverse problems for regular variation.
\newblock {\em J. Appl. Probab.}, 51A(Celebrating 50 Years of The Applied Probability Trust):229--248, 2014.

\bibitem{dav:sha:87}
P.~L.~Davies, D.~N.~Shanbhag.
\newblock A generalization of a theorem of Deny with applications in characterization theory.
\newblock {\em Quart. J. Math. Oxford}, 38:13--34, 1987.

\bibitem{den60}
J.~Deny.
\newblock Sur l'\'{e}quation de convolution $\mu = \mu \mc \sigma$ [About the evolution equation $\mu = \mu \mc \sigma$.]
\newblock {\em Seminaire Brelot-Choquet-Deny} (Th\'{e}orie du Potentiel), 4e ann\'{e}e, 1959/60, no.~5, p.1--11.

\bibitem{elstrodt}
J.~Elstrodt.
\newblock {\em Mass- und Integrationstheorie.} [Measure and integration theory.]
\newblock Springer, Berlin, 6th edition, 2009.

\bibitem{grillet}
P.~A.~Grillet.
\newblock {\em Abstract algebra}.
\newblock Springer, New York, 2nd edition, 2007.

\bibitem{jac:mic:ros:sam09}
M.~Jacobson, T.~Mikosch, J.~Rosi\'{n}ski, and G.~Samorodnitsky.
\newblock Inverse problems for regular variation of linear filters, a cancellation property for $\sigma$-finite measures and identification of stable laws.
\newblock {\em Ann. Probab.}, 19:210--242, 2009.

\bibitem{janich}
K.~J\"{a}nich.
\newblock {\em Funktionentheorie.} [Complex analysis.]
\newblock Springer, Berlin, 6th edition, 2004.

\bibitem{kab:sch:haan09}
Z.~Kabluchko, M.~Schlather, and L.~de~Haan.
\newblock Stationary max-stable fields associated to negative definite functions.
\newblock {\em Ann. Probab.}, 37:2042--2065, 2009.

\bibitem{mol:nag:21}
I.~Molchanov and F.~Nagel.
\newblock Diagonal Minkowski classes, zonoid equivalence, and stable laws.
\newblock {\em Communications in Contemporary Mathematics}, 23(02):1950091, 2021.

\bibitem{mos02}
K.~Mosler.
\newblock {\em Multivariate Dispersion, Central Regions and Depth. The Lift Zonoid Approach}, volume 165 of {\em Lect. Notes Statist.}
\newblock Springer, Berlin, 2002.

\bibitem{sam:taq94}
G.~Samorodnitsky and M.~S. Taqqu.
\newblock {\em Stable non-{Gaussian} Random Processes}.
\newblock Chapman \&\ Hall, New York, 1994.

\bibitem{schn2}
R.~Schneider.
\newblock {\em Convex Bodies. The {Brunn--Minkowski} Theory}.
\newblock Cambridge University Press, Cambridge, 2nd edition, 2014.

\bibitem{weil82}
W.~Weil.
\newblock Zonoide und verwandte Klassen konvexer K\"orper. [Zonoids and related classes of convex bodies.]
\newblock {\em Mh. Math.}, 94:73--84, 1982.

\bibitem{werner}
D.~Werner.
\newblock {\em Funktionalanalysis}. [Functional analysis.]
\newblock Springer, Berlin, 6th edition, 2007.

\end{thebibliography}
\end{document}